\documentclass[12pt,leqno]{article}

\usepackage[english]{babel}
\usepackage{amssymb}
\usepackage[mathscr]{eucal}
\usepackage{amsmath,amssymb,latexsym,theorem,bbm}
\usepackage{color,url}
\usepackage{appendix}
\usepackage{graphicx}
\newcommand{\NN}{\mathbb{N}}

\newcommand{\RR}{\mathbb{R}}

\newcommand{\bB}{{\boldsymbol{B}}}

\newcommand{\bc}{{\boldsymbol{c}}}
\newcommand{\bC}{{\boldsymbol{C}}}

\newcommand{\bD}{{\boldsymbol{D}}}

\newcommand{\Bf}{{\boldsymbol{f}}}

\newcommand{\tBf}{{\widetilde{\Bf}}}

\newcommand{\bG}{{\boldsymbol{G}}}
\newcommand{\tbG}{{\widetilde{\bG}}}
\newcommand{\bh}{{\boldsymbol{h}}}

\newcommand{\tbh}{\widetilde{\bh}}
\newcommand{\tildeh}{\widetilde{h}}
\newcommand{\oh}{{\overline{h}}}

\newcommand{\bI}{{\boldsymbol{I}}}
\newcommand{\tI}{\widetilde{I}}

\newcommand{\bM}{{\boldsymbol{M}}}

\newcommand{\bQ}{{\boldsymbol{Q}}}
\newcommand{\bR}{{\boldsymbol{R}}}
\newcommand{\bS}{{\boldsymbol{S}}}

\newcommand{\bV}{{\boldsymbol{V}}}

\newcommand{\bx}{{\boldsymbol{x}}}

\newcommand{\bz}{{\boldsymbol{z}}}

\newcommand{\bU}{{\boldsymbol{U}}}

\newcommand{\bpsi}{{\boldsymbol{\psi}}}

\newcommand{\bxi}{{\boldsymbol{\xi}}}
\newcommand{\bzero}{{\boldsymbol{0}}}

\newcommand{\cB}{{\mathcal B}}

\newcommand{\cC}{{\mathcal C}}
\newcommand{\tbC}{\widetilde{\bC}}
\newcommand{\ttbC}{\widetilde{\widetilde{\bC}}}
\newcommand{\cD}{{\mathcal D}}

\newcommand{\cF}{{\mathcal F}}

\newcommand{\cN}{{\mathcal N}}
\newcommand{\cP}{{\mathcal P}}

\newcommand{\cS}{{\mathcal S}}

\newcommand{\cW}{{\mathcal W}}

\newcommand{\dd}{\mathrm{d}}
\newcommand{\ee}{\mathrm{e}}
\newcommand{\ii}{\mathrm{i}}

\newcommand{\cont}{\mathrm{cont}}

\renewcommand{\tau}{\mathcal{T}}

\DeclareMathOperator*{\argmax}{arg\,max}

\newcommand{\EE}{\operatorname{\mathbb{E}}}
\newcommand{\PP}{\operatorname{\mathbb{P}}}

\newcommand{\ha}{\widehat{a}}
\newcommand{\hb}{\widehat{b}}

\newcommand{\htheta}{\widehat{\theta}}

\newcommand{\hbSigma}{\widehat{\Sigma}}

\newcommand{\hbpsi}{\widehat{\bpsi}}

\newcommand{\ttheta}{\widetilde{\theta}}

\newcommand{\btpsi}{{\widetilde{\bpsi}}}
\newcommand{\tmu}{\widetilde{\mu}}
\newcommand{\hkappa}{\widehat{\kappa}}
\newcommand{\hmu}{\widehat{\mu}}

\newcommand{\tS}{\widetilde{S}}

\newcommand{\tY}{\widetilde{Y}}

\newcommand{\hY}{\widehat{Y}}
\newcommand{\teta}{\widetilde{\eta}}
\newcommand{\tkappa}{\widetilde{\kappa}}

\newcommand{\vare}{\varepsilon}

\renewcommand{\mid}{\,|\,}

\renewcommand{\leq}{\leqslant}
\renewcommand{\geq}{\geqslant}

\newcommand{\stoch}{\stackrel{\PP}{\longrightarrow}}
\newcommand{\distr}{\stackrel{\cD}{\longrightarrow}}
\newcommand{\distre}{\stackrel{\cD}{=}}

\newcommand{\as}{\stackrel{{\mathrm{a.s.}}}{\longrightarrow}}

\newcommand{\bbone}{\mathbbm{1}}

\newcommand{\proofend}{\hfill\mbox{$\Box$}}

\numberwithin{equation}{section}

\theorembodyfont{\em}
\newtheorem{Lem}{Lemma}[section]
\newtheorem{Thm}[Lem]{Theorem}
\newtheorem{Pro}[Lem]{Proposition}

\theorembodyfont{\rm}
\newtheorem{Rem}[Lem]{Remark}

\begin{document}

\selectlanguage{english}
\pagenumbering{arabic}

\begin{center}
 {\bfseries\Large
  Asymptotic behavior of maximum likelihood estimators\\[2mm]
  for a jump-type Heston model}\\[5mm]
 {\sc\large
 M\'aty\'as $\text{Barczy}^{*,\diamond}$, \ Mohamed $\text{Ben Alaya}^{**}$,}\\[2mm]
 {\sc\large
  Ahmed $\text{Kebaier}^{***}$ \ and \ Gyula $\text{Pap}^{****}$}
\end{center}

\vskip0.2cm

\noindent
 * MTA-SZTE Analysis and Stochastics Research Group,
   Bolyai Institute, University of Szeged,
   Aradi v\'ertan\'uk tere 1, H--6720 Szeged, Hungary.

\noindent
 ** Laboratoire De Math\'ematiques Rapha\"el Salem, UMR 6085, Universit\'e De Rouen,
   Avenue de L'Universit\'e Technop\^ole du Madrillet, 76801 Saint-Etienne-Du-Rouvray, France.

\noindent
 *** Universit\'e Paris 13, Sorbonne Paris Cit\'e, LAGA, CNRS (UMR 7539),
   Villetaneuse, France.

\noindent
 **** Bolyai Institute, University of Szeged,
     Aradi v\'ertan\'uk tere 1, H--6720 Szeged, Hungary.

\noindent e--mails: barczy@math.u-szeged.hu (M. Barczy), \\
\phantom{e--mails:\,} mohamed.ben-alaya@univ-rouen.fr (M. Ben Alaya), \\
\phantom{e--mails:\,} kebaier@math.univ-paris13.fr (A. Kebaier), \\
\phantom{e--mails:\,} papgy@math.u-szeged.hu (G. Pap).

\noindent $\diamond$ Corresponding author

\renewcommand{\thefootnote}{}
\footnote{\textit{2010 Mathematics Subject Classifications\/}:
          60H10, 91G70, 60F05, 62F12.}
\footnote{\textit{Key words and phrases\/}:
 jump-type Heston model, maximum likelihood estimator.}
\vspace*{0.2cm}
\footnote{This research is supported by Laboratory of Excellence MME-DII,  Grant no.\ ANR11-LBX-0023-01
 (\texttt{http://labex-mme-dii.u-cergy.fr/}).
M\'aty\'as Barczy was supported between September 2016 and January 2017 by the
 ''Magyar \'Allami E\"otv\"os \"Oszt\"ond\'{\i}j 2016'' Grant no.\ 75141 funded by the Tempus
 Public Foundation, and from September 2017 by the J\'anos Bolyai Research Scholarship of the
 Hungarian Academy of Sciences.
Ahmed Kebaier benefited from the support of the chair
"Risques Financiers", Fondation du Risque.}

\vspace*{-10mm}

\begin{abstract}
We study asymptotic properties of maximum likelihood estimators of drift parameters
 for a jump-type Heston model based on continuous time observations,
 where the jump process can be any purely non-Gaussian L\'evy process of not necessarily
 bounded variation with a L\'evy measure concentrated on \ $(-1, \infty)$.
\ We prove strong consistency and asymptotic normality for all admissible
 parameter values except one, where we show only weak consistency and mixed normal
 (but non-normal) asymptotic behavior.
It turns out that the volatility of the price process is a measurable function of the price
 process.
We also present some numerical illustrations to confirm our results.
\end{abstract}

\section{Introduction}

Parameter estimation, especially studying asymptotic properties of maximum likelihood
 estimator (MLE) of drift parameters for Cox--Ingersoll--Ross (CIR) and Heston models
 is an active area of research mainly due to the wide range of applications of these
 models in financial mathematics.

The present paper gives a new contribution to the theory of asymptotic properties
 of MLE for jump-type Heston models based on continuous time observations.
Concerning related works, due to the vast literature on parameter estimation for Heston
 models, we will restrict ourselves to mention only papers that investigate the very same
 types of questions.
For a detailed and recent survey on parameter estimation for Heston models in general,
 see the Introduction of Barczy and Pap \cite{BarPap}.

Overbeck \cite{Ove} studied MLE of the drift parameters of the first coordinate process
 of a (diffusion type) Heston model (see \eqref{Heston_SDE}) based on continuous time observations,
 which is nothing else but a CIR process, also called square root process or Feller process.
Ben-Alaya and Kebaier \cite{BenKeb1}, \cite{BenKeb2} made a progress in MLE for the CIR
 process, giving explicit forms of joint Laplace transforms of the building blocks of
 this MLE as well.

Barczy and Pap \cite{BarPap} considered a Heston model
 \begin{align}\label{Heston_SDE}
  \begin{cases}
   \dd Y_t = (a - b Y_t) \, \dd t + \sigma_1 \sqrt{Y_t} \, \dd W_t , \\
   \dd X_t = (\alpha - \beta Y_t) \, \dd t
             + \sigma_2  \sqrt{Y_t}
               \bigl(\varrho \, \dd W_t
                     + \sqrt{1 - \varrho^2} \, \dd B_t\bigr) ,
  \end{cases} \qquad t \in [0, \infty) ,
 \end{align}
 where \ $a, \sigma_1, \sigma_2 \in (0, \infty)$, \ $b, \alpha, \beta \in \RR$,
 \ $\varrho \in (-1, 1)$ \ and \ $(W_t, B_t)_{t\in[0,\infty)}$ \ is a 2-dimensional
 standard Wiener process.
Here \ $(X_t)_{t\in[0,\infty)}$ \ is the log-price process of an asset,
 \ $(Y_t)_{t\in[0,\infty)}$ \ is its stochastic volatility (or instantaneous variance),
 \ $\sigma_1 \in (0, \infty)$ \ is the so-called volatility of the volatility, and
 \ $\varrho \in (-1, 1)$ \ is the correlation between the driving standard Wiener
 processes \ $(W_t)_{t\in[0,\infty)}$ \ and
 \ $(\varrho W_t + \sqrt{1 - \varrho^2} B_t)_{t\in[0,\infty)}$.
\ The MLE of the drift parameters \ $(a, b, \alpha, \beta)$ \ and its asymptotic behavior
 have been investigated based on continuous time observations \ $(Y_t,X_t)_{t\in[0,T]}$
 \ with \ $T \in (0, \infty)$ \ for all admissible parameter values
 (according to \ $b > 0$, \ $b = 0$, \ and \ $b < 0$).
\ It turned out that, for all \ $t\in[0,T]$, \ $Y_t$ \ is a measurable function of \ $(X_t)_{t\in[0,T]}$,
 \ hence, for the calculation of the MLE in question, one does not need the sample \ $(Y_t)_{t\in[0,T]}$.

The original Heston model (see Heston \cite{Hes}) takes the form
 \begin{align}\label{original_Heston_SDE}
  \begin{cases}
   \dd Y_t = \kappa (\theta - Y_t) \, \dd t + \sigma \sqrt{Y_t} \, \dd W_t , \\
   \dd S_t = \mu S_t \, \dd t
             + S_t \sqrt{Y_t}
               \bigl(\varrho \, \dd W_t
                     + \sqrt{1 - \varrho^2} \, \dd B_t\bigr) ,
  \end{cases} \qquad t \in [0, \infty) ,
 \end{align}
 where \ $(S_t)_{t\in[0,\infty)}$ \ is the price process of an asset, \ $\mu \in \RR$
 \ is the rate of return of the asset, \ $\theta \in (0, \infty)$ \ is the so-called
 long variance (long run average price variance, i.e., the limit of \ $\EE(Y_t)$ \ as
 \ $t \to \infty$), \ $\kappa \in (0, \infty)$ \ is the rate at which
 \ $(Y_t)_{t\in[0,\infty)}$ \ reverts to \ $\theta$, \ and \ $\sigma \in (0, \infty)$ \ is
 the so-called volatility of the volatility.
We call the attention that there are two differences between the models
 \eqref{Heston_SDE} and \eqref{original_Heston_SDE}.
Namely, in \eqref{original_Heston_SDE} the coefficient \ $\kappa$ \ can only be positive,
 while in \eqref{Heston_SDE} the corresponding coefficient \ $b$ \ can be an arbitrary
 real number.
In other words, the first coordinate process in \eqref{Heston_SDE} can be subcritical,
 critical or supercritical (according to \ $b > 0$, \ $b = 0$, \ and \ $b < 0$), but in
 \eqref{original_Heston_SDE} it can only be subcritical (since \ $\kappa > 0$).
\ Moreover, the second coordinate process in \eqref{original_Heston_SDE} is the price
 process, while in \eqref{Heston_SDE} it is the log-price process.

In this paper we study a jump-type Heston model (also called a stochastic volatility with
 jumps model, SVJ model)
 \begin{align}\label{j_Heston_SDE}
  \begin{cases}
   \dd Y_t = \kappa (\theta - Y_t) \, \dd t + \sigma \sqrt{Y_t} \, \dd W_t , \\
   \dd S_t = \mu S_t \, \dd t
             + S_t \sqrt{Y_t}
               \bigl(\varrho \, \dd W_t
                     + \sqrt{1 - \varrho^2} \, \dd B_t\bigr)
             + S_{t-} \, \dd L_t ,
  \end{cases} \qquad t \in [0, \infty) ,
 \end{align}
 where \ $(L_t)_{t\in[0,\infty)}$ \ is a purely non-Gaussian L\'evy process
 independent of \ $(W_t, B_t)_{t\in[0,\infty)}$ \ with L\'evy--Khintchine
 representation
 \begin{equation}\label{LK}
  \EE(\ee^{\ii u L_1})
  = \exp\left\{\ii \gamma u
               + \int_{-1}^\infty
                   (\ee^{\ii uz} - 1 - \ii u z \bbone_{(-1,1]}(z)) \, m(\dd z)\right\} ,
  \qquad u \in \RR ,
 \end{equation}
 where \ $\gamma \in \RR$ \ and \ $m$ \ is a L\'evy measure concentrated on
 \ $(-1, \infty)$ \ with \ $m(\{0\}) = 0$.
\ Here, let us recall that the L\'evy process \ $L$ \ has finite variation on each interval \ $[0, t]$,
 \ $t \in [0, \infty)$, \ if and only if \ $\int_{-1}^1 |z| \, m(\dd z) < \infty$, \ see, e.g.,
 Sato \cite[Theorem 21.9]{Sat}.
We point out that the assumption \ $\PP(Y_0\in[0,\infty), S_0\in(0,\infty))=1$ \ and
 the assumption in question on the support of the L\'evy measure \ $m$ \ assure
 that \ $\PP(\text{$S_t \in (0, \infty)$ \ for all $t \in [0, \infty)$}) = 1$
 \ (see Proposition \ref{Pro_j_Heston}), so the process \ $S$ \ can be used for modeling
 prices in a financial market.
From the point of view of financial mathematics, a natural question may occur
 concerning the model \eqref{j_Heston_SDE}.
Namely, is the drift coefficient of the second SDE in \eqref{j_Heston_SDE} well-adjusted
 in the sense that the discounted price process forms a martingale under some suitable
 equivalent martingale measure?
We renounce to consider this question, we just note that one may have to choose the parameter
 \ $\mu$ \ in an appropriate way to assure this property.
In Lamberton and Lapeyre \cite[Section 7]{LamLap} one can find a detailed discussion of the same type of question
 for a jump-type Black-Scholes model, where the jumps of the log-price process is modeled by a compound
 Poisson process.
They derived a necessary and sufficient condition for the drift coefficient of the underlying SDE
 in terms of the discounting factor and the parameters of the compound Poisson process in question
 in order that the discounted price process is a martingale, see \cite[page 146]{LamLap}.
For a good survey on jump-type Heston models, pricing and hedging in these models,
 see Runggaldier \cite{Run}.
In fact, the model \eqref{j_Heston_SDE} is quite popular in finance with the special choice
 of the L\'evy process \ $L$ \ as a compound Poisson process.
Namely, let
 \begin{align}\label{L_special}
  L_t := \sum_{i=1}^{\pi_t} (\ee^{J_i} - 1) , \qquad t \in [0, \infty) ,
 \end{align}
 where \ $(\pi_t)_{t\in[0,\infty)}$ \ is a Poisson process with intensity
 \ $\lambda \in (0, \infty)$,
 \ $(J_i)_{i\in\NN}$ \ is a sequence of independent identically distributed random variables
 having no atom at zero (i.e., \ $\PP(J_1=0)=0$), \ and being independent of \ $\pi$ \ as well.
We also suppose that \ $\pi$, \ $(J_i)_{i\in\NN}$, \ $W$ \ and \ $B$ \ are independent.
One can interpret \ $J$ \ as the jump size of the logarithm of the asset price.
 Then
 \begin{align*}
  \EE(\ee^{\ii u L_1})
  &= \exp\left\{\lambda \int_{-1}^\infty (\ee^{\ii uz} - 1) \, m(\dd z)\right\} \\
  &= \exp\left\{\ii u \lambda \int_{-1}^1 z \, m(\dd z)
               + \lambda
                 \int_{-1}^\infty
                  (\ee^{\ii uz} - 1 - \ii u z \bbone_{[-1,1]}(z)) \, m(\dd z)\right\} ,
  \qquad u \in \RR ,
 \end{align*}
 has the form \eqref{LK} with \ $m$ \ being the distribution of
 \ $\ee^{J_1} - 1$ \ and \ $\gamma = \lambda \int_{-1}^1 z \, m(\dd z)$.
 \ Moreover,
 \ $S_t$ \ takes the form
 \begin{align}\label{HESCPP}
  S_t
  = S_0
    \exp\biggl\{\int_0^t \Bigl(\mu - \frac{1}{2} Y_u\Bigr) \, \dd u
                + \int_0^t \sqrt{Y_u} \, (\varrho\, \dd  W_u + \sqrt{1 - \varrho^2}\, \dd  B_u)
                + \sum_{i=1}^{\pi_t} J_i\biggr\}
 \end{align}
 for \ $t \in [0, \infty)$, \ see \eqref{Solutions}.
We note that the SDE \eqref{j_Heston_SDE} with the L\'evy process \ $L$ \ given in
 \eqref{L_special} has been studied, e.g., by Bates \cite[equation (1)]{Bat},
 Bakshi et al. \cite[equations (1) and (2) with \ $R\equiv0$]{BakCaoChe},
 by Broadie and Kaya \cite[equations (30)-(31)]{BroKay}
 (where a factor \ $S_{t-}$ \ is missing from the last term of equation (30)),
 by Runggaldier \cite[Remark 3.1 with \ $\lambda_t \equiv \lambda$]{Run} and by Sun et al. \cite[equation (1) with $J_v=0$]{SunLiuGuo}.
 Bates \cite{Bat}, Bakshi et al. \cite{BakCaoChe} and Broadie and Kaya \cite{BroKay} have chosen the common distribution
 of \ $J$ \ as a normal distribution.
Bakshi et al. \cite{BakCaoChe} used this model for studying (European style) S\&P 500 options, e.g., they derived
 a practically implementable closed-form pricing formula.
Broadie and Kaya \cite{BroKay} gave an exact simulation algorithm for this model, further, they considered the pricing
  of forward start options in this model.
Sun et al. \cite{SunLiuGuo} have chosen the common distribution of \ $J$ \ as a normal distribution,
 a one-sided exponential distribution or a two-sided distribution, and they applied the Fourier-cosine
 series expansion method for pricing vanilla options under these jump-type Heston models.

The aim of this paper is to study the MLE of the parameter \ $\bpsi := (\theta, \kappa, \mu)$ \
 for the model \eqref{j_Heston_SDE} based on continuous time observations
 \ $(Y_t, S_t)_{t\in[0,T]}$ \ with \ $T \in (0, \infty)$, \ starting the process
 \ $(Y, S)$ \ from some deterministic initial value \ $(y_0, s_0) \in (0,\infty)^2$ \
 supposing that \ $\sigma$, \ $\varrho$, \ $\gamma$ \ and the L\'evy measure \ $m$ \ are known.
Here we stress that under these assumptions, the underlying statistical space
 corresponding to the parameters
 \ $(\kappa, \theta, \mu) \in (0, \infty)^2 \times \RR$
 \ is identifiable, however it would not be true for the statistical space
 corresponding to the parameters
 \ $(\kappa, \theta, \mu, \gamma) \in (0,\infty)^2 \times \RR^2$.
We call the attention that the MLE in question contains stochastic integrals
 with respect to \ $L$.
\ We prove that, for all \ $t\in[0,T]$, \ $L_t$ \ is a measurable function (i.e., a statistic)
 of \ $(S_t)_{t\in[0,T]}$, \ by providing a sequence of measurable functions of
 \ $(S_t)_{t\in[0,T]}$ \ converging in probability to \ $L_t$, \ see Remark \ref{observationL}
(note that this sequence depends on \ $\gamma$ \ and \ $m$ \ as well).
Further, it turns out that \ $Y_t$ \ for all \ $t \in [0, T]$, \ and the parameters \ $\sigma$
 \ and \ $\varrho$ \ are also measurable functions of \ $(S_t)_{t\in[0,T]}$, \ see Remarks
 \ref{observation} and \ref{Thm_MLE_cons_sigma_rho}, respectively.
Hence, for the calculation of the MLE in question, one needs only the sample
 \ $(S_t)_{t\in[0,T]}$, \ the parameter \ $\gamma$ \ and the L\'evy measure \ $m$ \ ($\gamma$
 \ and \ $m$ are needed for the reconstruction of \ $(L_t)_{t\in[0,T]}$).
\ Though we do not need to estimate the parameters \ $\sigma$ \ and \ $\varrho$,
 \ it is worth mentioning that the market microstructure effects may cause serious damage to
 the approximation of \ $\sigma$ \ and \ $\varrho$ \ given in Remark
 \ref{Thm_MLE_cons_sigma_rho} and to the MLE of \ $(\theta,\kappa,\mu)$ \ in case of high-frequency
 observations as in Zhang et al.\ \cite{ZMA}.
This type of question can be another interesting topic for future research.

The paper is organized as follows.
In Section \ref{Prel}, we prove that the SDE \eqref{j_Heston_SDE} has a pathwise
 unique strong solution (under some appropriate conditions), see Proposition
 \ref{Pro_j_Heston}, we recall a result  about the existence of a unique stationary
 distribution and ergodicity for the process \ $(Y_t)_{t\in[0,\infty)}$ \ given by
 the first equation in \eqref{j_Heston_SDE}, see Theorem \ref{Ergodicity}.
In Proposition \ref{Grigelionis}, we derive a Grigelionis representation for the process
 \ $(S_t)_{t\in[0,\infty)}$.
\ Further, we prove that for all \ $t\in[0,T]$, \ $L_t$ \ and \ $Y_t$ \ are measurable functions
 of \ $(S_t)_{t\in[0,T]}$, \ and we justify
 why we do not estimate the parameters \ $\sigma$ \ and \ $\varrho$, \ see Remarks \ref{observationL},
 \ref{observation} and \ref{Thm_MLE_cons_sigma_rho}.
Section \ref{section_EUMLE} is devoted to study the existence and uniqueness of the
 MLE \ $(\htheta_T, \hkappa_T, \hmu_T)$ \ of \ $(\theta, \kappa, \mu)$ \ based on
 observations \ $(Y_t, S_t)_{t\in[0,T]}$ \ with \ $T \in (0, \infty)$.
\ In Proposition \ref{Pro_MLE}, under appropriate conditions,
 we prove the unique existence of \ $(\htheta_T, \hkappa_T, \hmu_T)$, \ and we derive
 an explicit formula for it as well, see \eqref{MLE_coordinatewise}.
In Remark \ref{Luschgy}, we describe the connection with the so called score vector due to S{\o}rensen \cite{SorM}
 and the estimating equation due to Luschgy \cite{Lus2}, \cite{Lus} leading to the same estimator.
In Section \ref{section_CMLE}, we prove that the MLE of \ $(\theta, \kappa, \mu)$
 \ is strongly consistent if \ $\theta, \kappa \in (0, \infty)$ \ with
 \ $\theta\kappa \in \bigl(\frac{\sigma^2}{2},\infty\bigr)$, \ and weakly consistent
 if \ $\theta,\kappa \in (0, \infty)$ \ with \ $\theta\kappa = \frac{\sigma^2}{2}$,
 \ see Theorem \ref{Thm_MLE_cons} and Remark \ref{Rem_weak_const_a_szigma2},
 respectively.
Section \ref{section_AMLE} is devoted to investigate the asymptotic behaviour of the
 MLE of \ $(\theta, \kappa, \mu)$.
\ In Theorem \ref{Thm_MLE}, provided that \ $\theta,\kappa \in(0, \infty)$ \ with
 \ $\theta\kappa \in \bigl(\frac{\sigma^2}{2},\infty\bigr)$, \ we show that the MLE
 of \ $(\theta, \kappa, \mu)$ \ is asymptotically normal with a usual square root
 normalization \ $(T^{1/2})$, \ but as usual, the asymptotic covariance matrix
 depends on the unknown parameters \ $\theta$ \ and \ $\kappa$, \ as well.
To get around this problem, we also replace the normalization \ $T^{1/2}$ \ by a
 random one (depending only on the sample, but not on the parameters \ $\theta$,
 $\kappa$ \ and \ $\mu$) \ with the advantage that the MLE of
 \ $(\theta, \kappa, \mu)$ \ with the random scaling is asymptotically
 $3$-dimensional standard normal.
Theorem \ref{Thm_MLE=} is a counterpart of Theorem \ref{Thm_MLE} in some sense.
Namely, provided that \ $\theta, \kappa \in (0, \infty)$ \ with
 \ $\theta\kappa = \frac{\sigma^2}{2}$, \ we derive two limit theorems for
 the MLE \ $(\htheta_T, \hkappa_T, \hmu_T)$ \ with mixed normal limit distributions.
First, we have a non-random scaling, but for \ $\hmu_T$ \ instead of the usual
 scaling \ $T^{1/2}$ \ we have \ $T$; \ and then we have a random scaling as well.
We point out that, surprisingly, the limit distributions in Theorems \ref{Thm_MLE}
 and \ref{Thm_MLE=} do not depend on \ $L$ \ (roughly speaking, they do not depend on the jump part).
From a practical point of view, a natural question can occur, namely, how one can
 decide whether Theorems  \ref{Thm_MLE} and \ref{Thm_MLE=} can be applied (if yes, then which one),
 since one does not know the product \ $\theta\kappa$ \ of the unknown parameters \ $\theta$ \
 and \ $\kappa$ \ in advance.
To answer this question, one can build up a probe for testing the null hypothesis
 \ $\theta\kappa=\frac{\sigma^2}{2}$ \ against some alternative hypothesis, e.g., \ $\theta\kappa>\frac{\sigma^2}{2}$.
\ In Section \ref{section_simulation} we present some numerical illustrations of our
 limit theorems.
We close the paper with Appendices, where we recall certain sufficient conditions
 for the absolute continuity of probability measures induced by semimartingales together
 with a representation of the Radon--Nikodym derivative (Appendix \ref{App_LR}),
 some limit theorems for continuous local martingales for studying asymptotic behavior of
 \ $(\htheta_T, \hkappa_T, \hmu_T)$ \ (Appendix \ref{App_MCLT}) and a version of
 the continuous mapping theorem (Appendix \ref{CMT}), and we give an explicit formula
 for the non-normal but mixed normal density function of the limit distribution of
 \ $T (\hmu_T - \mu)$ \ as \ $T \to \infty$ \ in Theorem \ref{Thm_MLE=}
 (Appendix \ref{App_density}).

We call the attention that in both cases \ $\theta \kappa > \frac{\sigma^2}{2}$ \ and
 \ $\theta \kappa = \frac{\sigma^2}{2}$, \ the CIR process \ $Y$ \ has a unique
 stationary distribution and is ergodic, nevertheless, in case
 \ $\theta \kappa > \frac{\sigma^2}{2}$ \ the asymptotic limit distribution of the MLE
 of \ $\bpsi =(\theta,\kappa,\mu)$ \ is normal, while in case \ $\theta \kappa = \frac{\sigma^2}{2}$ \ it is
 mixed normal.
The interesting point is that we have an ergodic case with an asymptotically mixed
 normal (but non-normal) limit distribution.
The main difference between the two ergodic cases is that
 \ $\EE\bigl(\frac{1}{Y_\infty}\bigr) < \infty$ \ if
 \ $\theta \kappa > \frac{\sigma^2}{2}$, \ but
 \ $\EE\bigl(\frac{1}{Y_\infty}\bigr) = \infty$ \ if
 \ $\theta \kappa = \frac{\sigma^2}{2}$.

\section{Preliminaries}
\label{Prel}

The next proposition is about the existence and uniqueness of a strong
 solution of the SDE \eqref{j_Heston_SDE}.

\begin{Pro}\label{Pro_j_Heston}
Let \ $(\eta_0, \zeta_0)$ \ be a random vector such that \ $\eta_0$
 \ is independent of \ $(W_t)_{t\in[0,\infty)}$ \ satisfying
 \ $\PP(\eta_0 \in [0, \infty), \, \zeta_0 \in (0, \infty)) = 1$.
\ Then for all \ $\theta, \kappa, \sigma \in (0, \infty)$, \ $\mu \in \RR$
 \ and \ $\varrho \in (-1, 1)$, \ there is a
 (pathwise) unique strong solution \ $(Y_t, S_t)_{t\in[0, \infty)}$ \ of the
 SDE \eqref{j_Heston_SDE} such that
 \ $\PP((Y_0, S_0) = (\eta_0, \zeta_0)) = 1$ \ and
 \ $\PP(\text{$Y_t \in [0, \infty)$ and $S_t \in (0, \infty)$ for all $t \in [0, \infty)$}) = 1$.
\ Further,
 \begin{align}\label{Solutions}
  S_t = S_0
        \exp\biggl\{\int_0^t \Bigl(\mu - \frac{1}{2} Y_u\Bigr) \, \dd u
                    + \int_0^t
                       \sqrt{Y_u} \,
                       (\varrho\, \dd  W_u + \sqrt{1 - \varrho^2}\, \dd  B_u)
                    + L_t\biggr\}
        \prod_{u\in[0,t]} (1 + \Delta L_u) \ee^{-\Delta L_u}
 \end{align}
 for \ $t \in [0, \infty)$, \ where \ $\Delta L_u := L_u - L_{u-}$,
 \ $u \in (0, \infty)$, \ $\Delta L_0 := 0$, \ and the (possibly) infinite product is
 absolutely convergent.
If, in addition, \ $\theta \kappa \in \bigl[\frac{\sigma^2}{2}, \infty\bigr)$ \ and
 \ $\PP(\eta_0 \in (0, \infty)) = 1$, \ then
 \ $\PP(\text{$Y_t \in (0, \infty)$ for all $t \in [0, \infty)$}) = 1$.
\end{Pro}

Note that, due to Sato \cite[Theorem 21.3]{Sat}, for each \ $t \in (0, \infty)$,
 \ the product \ $\prod_{u\in[0,t]} (1 + \Delta L_u) \ee^{-\Delta L_u}$ \ in
 \eqref{Solutions} contains finitely many terms different from 1 if and only if
 \ $m((-1, 1]) < \infty$.

\noindent{\bf Proof of Proposition \ref{Pro_j_Heston}.}
By a theorem due to Yamada and Watanabe (see, e.g., Karatzas and Shreve
 \cite[Proposition 5.2.13]{KarShr}), the strong uniqueness holds for the first
 equation in \eqref{j_Heston_SDE}.
By Ikeda and Watanabe \cite[Example 8.2, page 221]{IkeWat}, there is a
 (pathwise) unique non-negative strong solution \ $(Y_t)_{t\in[0,\infty)}$ \ of
 the first equation in \eqref{j_Heston_SDE} with any initial value \ $\eta_0$
 \ such that \ $\PP(\eta_0 \in [0, \infty)) = 1$.
\ The second equation in \eqref{j_Heston_SDE} can be written in the form
 \[
   \dd S_t = S_{t-} \, \dd L^*_t , \qquad t \in [0, \infty) ,
 \]
 where
 \begin{align}\label{L*}
   L^*_t
   := \mu t
      + \int_0^t
         \sqrt{Y_u} \bigl(\varrho \, \dd W_u + \sqrt{1 - \varrho^2} \, \dd B_u\bigr)
      + L_t , \qquad t \in [0, \infty) ,
 \end{align}
 is a semimartingale, since the process \ $(\sqrt{Y_t})_{t\in[0,\infty)}$ \ has
 continuous sample paths almost surely and hence locally bounded almost surely yielding that
 \ $\bigl(\int_0^t \sqrt{Y_u}
           (\varrho \, \dd W_u + \sqrt{1 - \varrho^2} \, \dd B_u)\bigr)_{t\in[0,\infty)}$
 \ is a square integrable martingale, and since \ $L$ \ is a semimartingale being a
 L\'evy process (see, e.g., Jacod and Shiryaev \cite[Corollary II.4.19]{JSh}).
Using \ $\Delta L^*_t = \Delta L_t$, \ $t \in [0, \infty)$, \ and Theorem 1 in
 Jaschke \cite{Jas}, which is a generalization of the Dol\'eans--Dade exponential formula
 (see, e.g., Jacod and Shiryaev \cite[I.4.61]{JSh}), we obtain
 \begin{align*}
  S_t &= S_0 \exp\left\{ L_t^* - L_0^*
                     - \frac{1}{2}
                       \langle (L^*)^\cont \rangle_t \right\}
         \prod_{u\in[0,t]} (1 + \Delta L^*_u)\ee^{-\Delta L^*_u} \\
      &= S_0 \exp\left\{ \mu t
                     + \int_0^t
                        \sqrt{Y_u} \bigl(\varrho \, \dd W_u + \sqrt{1 - \varrho^2}
                        \, \dd B_u\bigr)
                     + L_t
                     - \frac{1}{2}\int_0^t Y_u \, \dd u \right\}
         \prod_{u\in[0,t]} (1 + \Delta L_u)\ee^{-\Delta L_u} ,
 \end{align*}
 where \ $(\langle (L^*)^\cont \rangle_t)_{t\in[0,\infty)}$ \ denotes the
 (predictable) quadratic variation process of the continuous martingale part
 \ $(L^*)^\cont$ \ of \ $L^*$, \ and the (possibly) infinite product is absolutely
 convergent.
Here we used that
 \ $\langle (L^*)^\cont \rangle_t = \int_0^t Y_u\,\dd u$, \ $t\in[0,\infty)$,
 \ being a consequence of the fact that
 \begin{align}\label{L_cont_mart_part}
  (L^*)^\cont_t
  = \int_0^t
     \sqrt{Y_u} \bigl(\varrho \, \dd W_u + \sqrt{1 - \varrho^2} \, \dd B_u\bigr) , \qquad
  t \in [0, \infty) ,
 \end{align}
 which can be checked as follows.
The L\'evy--It\^{o}'s representation of \ $L$ \ takes the form
 \begin{equation}\label{Levy_Ito}
  \begin{aligned}
   L_t & = \lim_{\delta\downarrow 0}
           \int_{(0,t]} \int_{\{\delta<|z|\leq 1\}}
            z \big(\mu^L(\dd u,\dd z) - \dd u \, m(\dd z)\big)
           + \int_{(0,t]} \int_{\{1<|z|<\infty\}}
              z \,\mu^L(\dd u,\dd z)
           + \gamma t \\
       & =: \int_0^t \int_{\RR}
             h_1(z) \, \big(\mu^L(\dd u,\dd z) - \dd u \, m(\dd z)\big)
            + \int_0^t \int_{\RR} (z - h_1(z)) \, \mu^L(\dd u,\dd z)
             + \gamma t
 \end{aligned}
 \end{equation}
 for \ $t \in [0, \infty)$, \ where
 \ $\mu^L(\dd u,\dd z)
    := \sum_{v\in[0,\infty)}
        \bbone_{\{\Delta L_v\ne 0\}} \vare_{(v,\Delta L_v)}(\dd u,\dd z)$
 \ is the integer-valued Poisson random measure on \ $[0, \infty) \times \RR$
 \ associated with the jumps of the process \ $L$, \ $\vare_{(v,x)}$ \ denotes
 the Dirac measure at the point \ $(v, x) \in [0,\infty) \times \RR$, \ and
 \begin{equation}\label{h1}
  h_1(z) := z \bbone_{[-1,1]}(z) , \qquad z \in \RR ,
 \end{equation}
 is a truncation function, see, e.g., Sato \cite[Theorem 19.2]{Sat}.
The first term in \eqref{Levy_Ito} is a purely discontinuous local martingale,
 see, e.g., Jacod and Shiryaev \cite[Definitions II.1.27]{JSh}.
The second term in \eqref{Levy_Ito} can be written as a finite sum (see, e.g., Sato \cite[Lemma 20.1]{Sat})
 \[
   \int_0^t \int_{\RR} (z - h_1(z)) \, \mu^L(\dd u,\dd z)
   = \sum_{u\in[0,t]} \bbone_{\{|\Delta L_u|>1\}} \Delta L_u , \qquad t \in [0, \infty) ,
 \]
 which is a compound Poisson process with L\'evy-Khintchine representation
  \[
    \EE\left(  \ee^{\ii\theta \sum_{u\in[0,1]} \bbone_{\{|\Delta L_u|>1\}} \Delta L_u}\right)
        = \exp\left\{\int_1^\infty (\ee^{\ii\theta u} - 1)\, m(\dd u) \right\},
        \qquad \theta\in\RR.
  \]
Hence it is a process with finite variation over each finite interval \ $[0, t]$,
 \ $t \in [0, \infty)$, \ see, e.g., Sato \cite[Theorem 21.9]{Sat}.
\ Consequently, we conclude \eqref{L_cont_mart_part}.
An alternative way for deriving \eqref{L_cont_mart_part} is as follows.
Using \eqref{Levy_Ito}, the process \ $L^*$ \ can be written in the form III.2.23 in
 Jacod and Shiryaev \cite{JSh}, and hence, by Jacod and Shiryaev
 \cite[Remarks III.2.28, part 1)]{JSh}, we get \eqref{L_cont_mart_part}.
Thus the (pathwise) unique strong solution \ $(S_t)_{t\in[0,\infty)}$ \ of the second
 equation in \eqref{j_Heston_SDE} is given by \eqref{Solutions}.
Further,
 \[
  \PP(\text{$\Delta L_t \in (-1, \infty)$ \ for all \ $t \in [0, \infty)$}) = 1,
 \]
 since the L\'evy measure \ $m$ \ of \ $L$ \ is concentrated on
 \ $(-1, 0) \cup (0, \infty)$.
\ Using again \ $\Delta L^*_t = \Delta L_t$, \ $t \in [0, \infty)$, \ we obtain
 \ $\PP(\text{$\Delta L^*_t \in (-1, \infty)$ for all $t \in [0, \infty)$}) = 1$,
 \ and hence \ $\PP(\text{$S_t \in (0, \infty)$ for all $t \in [0, \infty)$}) = 1$.
\ Indeed, if \ $S_0 = 1$, \ then this follows, e.g., from Theorem I.4.61 (c)
 in Jacod and Shiryaev \cite{JSh}, hence, in general, this is a consequence of
 formula \eqref{Solutions} and \ $\PP(S_0 \in (0, \infty)) = 1$.

The proof of the last statement can be found, e.g., in Ikeda and Watanabe
 \cite[Chapter IV, Example 8.2]{IkeWat} and in Lamberton and Lapeyre
 \cite[Proposition 6.2.4]{LamLap}.
\proofend

In the sequel \ $\stoch$, \ $\distr$ \ and \ $\as$ \ will denote convergence
 in probability, in distribution and almost surely, respectively.

The following result states the existence of a unique stationary distribution
 and the ergodicity for the process \ $(Y_t)_{t\in[0,\infty)}$ \ given by the
 first equation in \eqref{j_Heston_SDE}, see, e.g., Feller
 \cite{Fel}, Cox et al.\ \cite[Equation (20)]{CoxIngRos}, Li and Ma
 \cite[Theorem 2.6]{LiMa} or Theorem 3.1 with \ $\alpha = 2$ \ and Theorem 4.1
 in Barczy et al.\ \cite{BarDorLiPap2}.

\begin{Thm}\label{Ergodicity}
Let \ $\theta, \kappa, \sigma \in (0, \infty)$.
\ Let \ $(Y_t)_{t\in[0,\infty)}$ \ be the unique strong solution of the first
 equation of the SDE \eqref{j_Heston_SDE} satisfying
 \ $\PP(Y_0 \in [0, \infty)) = 1$.
 \renewcommand{\labelenumi}{{\rm(\roman{enumi})}}
 \begin{enumerate}
  \item
   Then \ $Y_t \distr Y_\infty$ \ as \ $t \to \infty$, \ and the distribution
    of \ $Y_\infty$ \ is given by
    \begin{align}\label{Laplace}
     \EE(\ee^{-\lambda Y_\infty})
     = \Bigl(1 + \frac{\sigma^2}{2\kappa} \lambda \Bigr)^{-\frac{2\theta\kappa}{\sigma^2}} ,
     \qquad \lambda \in [0, \infty) ,
    \end{align}
    i.e., \ $Y_\infty$ \ has Gamma distribution with parameters
    \ $2 \theta \kappa / \sigma^2$ \ and \ $2 \kappa / \sigma^2$, \ hence
    \[
      \EE(Y_\infty^K)
      = \frac{\Gamma\bigl(\frac{2\theta\kappa}{\sigma^2} + K \bigr)}
             {\bigl(\frac{2\kappa}{\sigma^2}\bigr)^K
              \Gamma\bigl(\frac{2\theta\kappa}{\sigma^2}\bigr)} , \qquad
      K \in \Bigl( -\frac{2\theta\kappa}{\sigma^2}, \infty \Bigr) .
    \]
   Especially, \ $\EE(Y_\infty) = \theta$.
   \ Further, if \ $\theta \kappa \in \bigl( \frac{\sigma^2}{2}, \infty \bigr)$,
    \ then
    \ $\EE\bigl(\frac{1}{Y_\infty}\bigr) = \frac{2\kappa}{2\theta\kappa-\sigma^2}$.
 \item
  For all Borel measurable functions \ $f : \RR \to \RR$ \ such that
   \ $\EE(|f(Y_\infty)|) < \infty$, \ we have
   \begin{equation}\label{ergodic}
    \frac{1}{T} \int_0^T f(Y_u) \, \dd u \as \EE(f(Y_\infty)) \qquad
    \text{as \ $T \to \infty$.}
   \end{equation}
\end{enumerate}
\end{Thm}

Note that, by Proposition \ref{Pro_j_Heston}, the process \ $(Y_t, S_t)_{t\in[0,\infty)}$
 \ is a semimartingale, see, e.g., Jacod and Shiryaev \cite[I.4.34]{JSh}.
Now we derive a so-called Grigelionis form for the semimartingale
 \ $(S_t)_{t\in[0,\infty)}$, \ see, e.g., Jacod and Shiryaev
 \cite[III.2.23]{JSh} or Jacod and Protter \cite[Theorem 2.1.2]{JacPro}.

\begin{Pro}\label{Grigelionis}
Let \ $\theta, \kappa, \sigma \in (0, \infty)$, \ $\mu \in \RR$,
 \ $\varrho \in (-1, 1)$.
\ Let \ $(Y_t, S_t)_{t\in[0,\infty)}$ \ be the unique strong solution of the SDE
 \eqref{j_Heston_SDE} satisfying
 \ $\PP(Y_0 \in [0, \infty), \, S_0 \in (0, \infty)) = 1$.
\ Then
 \begin{equation}\label{SLevy_Ito}
  \begin{aligned}
   S_t &= S_0 + (\mu + \gamma) \int_0^t S_u \, \dd u
          + \int_0^t
             \left(\int_\RR (h_1(S_{u-} z)- S_{u-}h_1(z)) \, m(\dd z)\right)
             \dd u \\
       &\quad
          + \int_0^t
             S_u \sqrt{Y_u} \bigl(\varrho \, \dd W_u + \sqrt{1 - \varrho^2} \,
             \dd B_u\bigr) \\
       &\quad
          + \int_0^t \int_\RR
             h_1(S_{u-}z) \, \bigl(\mu^L(\dd u,\dd z) - \dd u \, m(\dd z)\bigr)
          + \int_0^t \int_\RR (S_{u-}z - h_1(S_{u-}z)) \, \mu^L(\dd u,\dd z)
  \end{aligned}
 \end{equation}
 for \ $t \in [0, \infty)$, \ where \ $h_1$ \ is defined in \eqref{h1}.
\end{Pro}

\noindent{\bf Proof.}
Using \eqref{Levy_Ito} and Proposition II.1.30 in Jacod and Shiryaev \cite{JSh}, we obtain
 \begin{align*}
  S_t &= S_0 + (\mu + \gamma) \int_0^t S_u \, \dd u
         + \int_0^t
            S_u \sqrt{Y_u} \bigl(\varrho \, \dd W_u + \sqrt{1 - \varrho^2} \,
            \dd B_u\bigr) \\
      &\quad
         + \int_0^t \int_\RR
            S_{u-} h_1(z) \, \big(\mu^L(\dd u,\dd z) - \dd u \, m(\dd z)\big)
         + \int_0^t \int_\RR S_{u-} (z - h_1(z)) \, \mu^L(\dd u,\dd z)
 \end{align*}
 for \ $t \in [0, \infty)$.
\ In order to prove the statement, it is enough to show
 \begin{gather}\label{G1}
  \int_0^t \int_\RR
   S_{u-} h_1(z) \, \big(\mu^L(\dd u,\dd z) - \dd u \, m(\dd z)\big)
   = I_1 - I_2 , \\
  \int_0^t \int_\RR
   S_{u-} (z - h_1(z)) \, \mu^L(\dd u,\dd z)
   = I_3 + I_4 , \label{G2}
 \end{gather}
 with
 \begin{align*}
  I_1 &:= \int_0^t \int_\RR
           h_1(S_{u-}z) \, \bigl(\mu^L(\dd u,\dd z) - \dd u \, m(\dd z)\bigr) , \\
  I_2 &:= \int_0^t \int_\RR
           (h_1(S_{u-} z)- S_{u-}h_1(z))
           \, \bigl(\mu^L(\dd u,\dd z) - \dd u \, m(\dd z)\bigr) , \\
  I_3 &:= \int_0^t \int_\RR
           (S_{u-}z - h_1(S_{u-}z)) \, \mu^L(\dd u,\dd z) , \\
  I_4 &:= \int_0^t \int_\RR
           (h_1(S_{u-} z)- S_{u-}h_1(z)) \, \mu^L(\dd u,\dd z) ,
 \end{align*}
 and the equality
 \begin{equation}\label{drift}
  I_4 - I_2 = I_5 \qquad \text{with} \qquad
  I_5 := \int_0^t \left(\int_\RR (h_1(S_{u-} z)- S_{u-}h_1(z)) \, m(\dd z)\right) \dd u .
 \end{equation}
For the equations \eqref{G1}, \eqref{G2} and \eqref{drift}, it suffices to check the
 existence of \ $I_2$, \ $I_3$ \ and \ $I_5$.

First note that for every \ $s \in (0, \infty)$ \ we have
 \begin{align}\label{help_ITO_atiras}
   h_1(s z)- sh_1(z)
   = \begin{cases}
      s z \bbone_{\{1<|z|\leq\frac{1}{s}\}}
       & \text{if \ $s \in (0, 1)$, \ $z \in \RR$,} \\
      0 & \text{if \ $s = 1$, \ $z \in \RR$,} \\
      -s z \bbone_{\{\frac{1}{s}<|z|\leq1\}}
       & \text{if \ $s \in (1, \infty)$, \ $z \in \RR$.}
     \end{cases}
 \end{align}
The existence of \ $I_2$ \ will be a consequence of \ $I_2 = I_{2,1} - I_{2,2} - I_{2,3}$
 \ with
 \begin{align*}
  I_{2,1} &:= \int_0^t \int_\RR
               S_{u-} z \bbone_{\{1<|z|\leq\frac{1}{S_{u-}}\}} \bbone_{\{S_{u-}\in(0,1)\}}
               \, \mu^L(\dd u,\dd z) , \\
  I_{2,2} &:= \int_0^t \int_\RR
               S_{u-} z \bbone_{\{1<|z|\leq\frac{1}{S_{u-}}\}} \bbone_{\{S_{u-}\in(0,1)\}}
               \, \dd u \, m(\dd z) , \\
  I_{2,3} &:= \int_0^t \int_\RR
               S_{u-} z \bbone_{\{\frac{1}{S_{u-}}<|z|\leq1\}}
               \bbone_{\{S_{u-}\in(1,\infty)\}}
               \, \bigl(\mu^L(\dd u,\dd z) - \dd u \, m(\dd z)\bigr) .
 \end{align*}
Here we have
 \[
   |I_{2,1}|
   \leq
   \int_0^t \int_\RR
    |S_{u-} z| \bbone_{\{1<|z|\leq\frac{1}{S_{u-}}\}} \bbone_{\{S_{u-}\in(0,1)\}}
    \, \mu^L(\dd u,\dd z)
   \leq \int_0^t \int_\RR \bbone_{\{1<|z|\}} \, \mu^L(\dd u,\dd z)
   < \infty ,
 \]
 see, e.g., Sato \cite[Lemma 20.1]{Sat}.
Moreover,
 \begin{align*}
  |I_{2,2}|
  &\leq
   \int_0^t \int_\RR
    |S_{u-} z| \bbone_{\{1<|z|\leq\frac{1}{S_{u-}}\}} \bbone_{\{S_{u-}\in(0,1)\}}
    \, \dd u \, m(\dd z) \\
  &\leq
   \int_0^t \int_\RR \bbone_{\{1<|z|\}} \, \dd u \, m(\dd z)
   = t m(\{z \in \RR : |z| > 1 \})
   < \infty .
 \end{align*}
Further, the function
 \ $\Omega \times [0, \infty) \times \RR \ni (\omega, t, z) \mapsto h_1(z)$
 \ belongs to \ $G_{\mathrm{loc}}(\mu^L)$, \ see Jacod and Shiryaev
 \cite[Definitions II.1.27, Theorem II.2.34]{JSh}.
We have
 \ $|z \bbone_{\{\frac{1}{S_{u-}}<|z|\leq1\}} \bbone_{\{S_{u-}\in(1,\infty)\}}|
    \leq |h_1(z)|$,
 \ hence, by the definition of \ $G_{\mathrm{loc}}(\mu^L)$, \ the function
 \ $\Omega \times [0, \infty) \times \RR \ni (\omega, t, z)
    \mapsto z \bbone_{\{\frac{1}{S_{u-}}<|z|\leq1\}} \bbone_{\{S_{u-}\in(1,\infty)\}}$
 \ also belongs to \ $G_{\mathrm{loc}}(\mu^L)$.
\ By Jacod and Shiryaev \cite[Proposition II.1.30]{JSh}, we conclude that the function
 \ $\Omega \times [0, \infty) \times \RR \ni (\omega, t, z)
    \mapsto S_{u-} z \bbone_{\{\frac{1}{S_{u-}}<|z|\leq1\}} \bbone_{\{S_{u-}\in(1,\infty)\}}$
 \ also belongs to \ $G_{\mathrm{loc}}(\mu^L)$, \ thus the integral \ $I_{2,3}$ \ exists,
 and hence we obtain the existence of \ $I_2$, \ and hence of \ $I_1$.

Next observe that we have \ $\Delta S_t = S_{t-} \Delta L_t$, \ $t \in [0, \infty)$, \ see,
 e.g., Jacod and Shiryaev \cite[page 60, formula (5)]{JSh}.
Consequently,
 \[
   I_3
   = \int_0^t \int_\RR
      S_{u-} z \bbone_{\{|S_{u-}z|>1\}} \, \mu^L(\dd u,\dd z)
   = \sum_{u\in[0,t]} S_{u-} (\Delta L_u) \bbone_{\{|S_{u-}\Delta L_u|>1\}}
   = \sum_{u\in[0,t]} \Delta S_u \bbone_{\{|\Delta S_u|>1\}}
 \]
 is a finite sum, since the process \ $(S_t)_{t\in[0,\infty)}$ \ admits c\`adl\`ag
 trajectories, hence there can be at most finitely many points \ $u \in [0, t]$ \ at which
 the jump \ $|\Delta S_u|$ \ exceeds 1, see, e.g., Billingsley \cite[page 122]{Bil}.
Thus we obtain the existence of \ $I_3$, \ and hence of \ $I_4$.

Finally, we have
 \begin{align*}
  |I_5|
  &\leq
   \int_0^t
    \left(\int_\RR
           |S_{u-} z| \bbone_{\{1<|z|\leq\frac{1}{S_{u-}}\}} \bbone_{\{S_{u-}\in(0,1)\}}
           \, m(\dd z)\right) \dd u \\
  &\phantom{\leq}
    + \int_0^t
       \left(\int_\RR
              |S_{u-} z| \bbone_{\{\frac{1}{S_{u-}}<|z|\leq1\}}
              \bbone_{\{S_{u-}\in(1,\infty)\}}
              \, m(\dd z)\right) \dd u \\
  &\leq
   \int_0^t
    \left(\int_\RR \bbone_{\{1<|z|\}} \, m(\dd z)\right) \dd u
   + \left(\int_\RR |S_{u-} z|^2 \bbone_{\{|z|\leq1\}} \, m(\dd z)\right) \dd u \\
  &= t m(\{z \in \RR : |z| > 1\})
     + \int_0^t S_{u-}^2 \, \dd u \int_{-1}^1 |z|^2 \, m(\dd z)
   < \infty ,
 \end{align*}
 hence we conclude the existence of \ $I_5$.
\proofend

In the next remark, we show that, for all \ $t\in[0,T]$, \ $L_t$ \ is a measurable function
 of \ $(S_t)_{t\in[0,T]}$ \ depending on the parameter \ $\gamma$ and the L\'evy
  measure \ $m$.

\begin{Rem}\label{observationL}
For all \ $t\in[0,T]$ \ and \ $\delta\in(0,1)$,
 \begin{equation*}
  \begin{aligned}
    \int_{(0,t]} \int_{\{\delta<|z|\leq 1\}}
             z \big(\mu^L(\dd u,\dd z) - \dd u \, m(\dd z)\big)
     & = \sum_{u\in[0,t]} \bbone_{\{ \delta < |\Delta L_u|\leq 1 \}}\Delta L_u
        - \int_{(0,t]} \int_{\{\delta<|z|\leq 1\}} z\,\dd u \, m(\dd z) \\
     & = \sum_{u\in[0,t]}
       \bbone_{\{ \delta < \frac{|\Delta S_u|}{S_{u-}}\leq 1 \}} \frac{\Delta S_u}{S_{u-}}
       - t \int_{\{\delta<|z|\leq 1\}} z\, m(\dd z),
 \end{aligned}
 \end{equation*}
 which is a measurable function of \ $(S_t)_{t\in[0,T]}$.
Similarly, for all \ $t\in[0,T]$,
 \begin{equation*}
  \begin{aligned}
    \int_{(0,t]} \int_{\{1<|z|<\infty\}} z \,\mu^L(\dd u,\dd z)
       = \sum_{u\in[0,t]} \bbone_{\{ |\Delta L_u| > 1\}} \Delta L_u
       = \sum_{u\in[0,t]} \bbone_{\{ \frac{|\Delta S_u|}{S_{u-}}>1 \}} \frac{\Delta S_u}{S_{u-}},
 \end{aligned}
 \end{equation*}
 which is a measurable function of \ $(S_t)_{t\in[0,T]}$ \ as well.
Hence, using \eqref{Levy_Ito}, for all \ $t\in[0,T]$,
 \begin{equation*}
  \sum_{u\in[0,t]} \bbone_{\{ \frac{|\Delta S_u|}{S_{u-}}>\delta \}} \frac{\Delta S_u}{S_{u-}}
      - t \int_{\{\delta<|z|\leq 1\}} z\, m(\dd z)
      + \gamma t
    \stoch L_t \quad\text{as \ $\delta\downarrow 0$,}
 \end{equation*}
 yielding that \ $L_t$ \ is a measurable function of \ $(S_t)_{t\in[0,T]}$.
\ In the special case of
 \begin{equation}\label{L}
  L_t = \sum_{s\in[0,t]} \Delta L_s , \qquad t \in [0, \infty),
 \end{equation}
the above statement readily follows from \ $\Delta L_s = \frac{\Delta S_s}{S_{s-}}$, $s\in[0,\infty)$.
\ Condition \eqref{L} is satisfied if
 \ $\int_{-1}^1 |z| \, m(\dd z) < \infty$ \ and
 \ $\gamma = \int_{-1}^1 z \, m(\dd z)$, \
 since, by \eqref{LK},
 \begin{align*}
   L_t = \int_0^t \int_\RR z \, \mu^L(\dd u,\dd z)
             + t \left(\gamma - \int_{-1}^1 z \, m(\dd z)\right)
       = \sum_{s\in[0,t]} \Delta L_s
              + t \left(\gamma - \int_{-1}^1 z \, m(\dd z)\right)
 \end{align*}
 for \ $t \in [0, \infty)$, \ see Sato \cite[Theorem 19.3]{Sat}.
Recall that, using \eqref{LK}, the L\'evy process
 \ $L$ \ has finite variation on each interval \ $[0, t]$, \ $t \in [0, \infty)$, \ if and
 only if \ $\int_{-1}^1 |z| \, m(\dd z) < \infty$, \ see, e.g., Sato
 \cite[Theorem 21.9]{Sat}.
For example, it is satisfied for a compound Poisson process given in
 \eqref{L_special}, where \ $m$ is a probability measure.
\proofend
\end{Rem}

In the next remark, we show that, for all \ $t\in[0,T]$, \ $Y_t$ \ is a measurable function
 of \ $(S_t)_{t\in[0,T]}$.

\begin{Rem}\label{observation}
Let \ $\theta, \kappa, \sigma \in (0, \infty)$, \ $\mu \in \RR$,
 \ $\varrho \in (-1, 1)$.
\ Let \ $(Y_t, S_t)_{t\in[0,\infty)}$ \ be the unique strong solution of the SDE
 \eqref{j_Heston_SDE} satisfying
 \ $\PP(Y_0 \in [0, \infty), \, S_0 \in (0, \infty)) = 1$.
\ The Grigelionis representation given in Proposition \ref{Grigelionis} implies that the
 continuous martingale part \ $S^\cont$ \ of \ $S$ \ is
 \begin{align}\label{S_cont_mart_part}
   S^\cont_t
    = \int_0^t
       S_u \sqrt{Y_u} \bigl(\varrho \, \dd W_u + \sqrt{1 - \varrho^2} \, \dd B_u\bigr),
 \qquad t \in [0, \infty),
 \end{align}
 see Jacod and Shiryaev \cite[III.2.28 Remarks, part 1)]{JSh}.
Consequently, the (predictable) quadratic variation process of \ $S^\cont$ \ is
 \ $\langle S^\cont \rangle_t = \int_0^t S_u^2 Y_u \, \dd u$, \ $t \in [0, \infty)$.
\ Since
 \[
   \PP(\text{$S_t, \, S_{t-} \in (0, \infty)$ \ for all \ $t \in [0, \infty)$}) = 1
  \]
 with the convention \ $S_{0-} := S_0$ \ (due to Proposition \ref{Pro_j_Heston}),
\ one can apply It\^o's rule to the function
 \ $f(x) = \log(x)$, \ $x \in (0, \infty)$, \ for which \ $f'(x) = 1 / x$,
 \ $f''(x) = - 1 / x^2$, \ $x \in (0, \infty)$, \ and we obtain
 \begin{align}\label{logS_SDE}
  \begin{aligned}
   \log S_T
   &= \log S_0 + \int_0^T \frac{\dd S_u}{S_{u-}}
      - \frac{1}{2} \int_0^T \frac{1}{S_u^2} \, \dd\langle S^\cont \rangle_u
      + \sum_{u\in[0,T]}
         \biggl(\log S_u - \log S_{u-} - \frac{1}{S_{u-}} \, \Delta S_u \biggr) \\
   &= \log S_0 + \mu T
      + \int_0^T \sqrt{Y_u} (\varrho \, \dd W_u + \sqrt{1 - \varrho^2} \, \dd B_u)
      - \frac{1}{2} \int_0^T Y_u \, \dd u + L_T \\
   &\quad
      + \sum_{u\in[0,T]}
         \biggl(\log\frac{S_u}{S_{u-}} + 1 - \frac{S_u}{S_{u-}}\biggr)
  \end{aligned}
 \end{align}
 for \ $T \in [0, \infty)$, \ see, e.g.,
 von Weizs\"acker and Winkler \cite[Theorem 8.4.1]{WeiWin}.
All terms in \eqref{logS_SDE} are well-defined.
In particular, the last term is a process with finite variation over each finite
 interval \ $[0, t]$, \ $t \in [0, \infty)$, \ see, e.g., Sato
 \cite[Lemma 21.8.(iii)]{Sat}.
Taking into account of the L\'evy--It\^{o}'s representation \eqref{Levy_Ito} of \ $L$,
 \ we conclude that the continuous martingale part \ $(\log S)^\cont$ \ of \ $\log S$
 \ is
 \ $(\log S)^\cont_t
    = \int_0^t \sqrt{Y_u} \bigl(\varrho \, \dd W_u + \sqrt{1 - \varrho^2} \, \dd B_u\bigr)$,
 \ $t \in [0, \infty)$.
\ Hence the (predictable) quadratic variation process of \ $(\log S)^\cont$ \ is
 \[
   \langle (\log S)^\cont \rangle_t = \int_0^t Y_u \, \dd u , \qquad t \in [0, \infty) .
 \]
By Theorem I.4.47 a) in Jacod and Shiryaev \cite{JSh},
 \[
   \sum_{i=1}^{\lfloor nt\rfloor} (\log S_{\frac{i}{n}} - \log S_{\frac{i-1}{n}})^2
   \stoch [\log S]_t \qquad \text{as \ $n \to \infty$,} \quad
   t \in [0, \infty) ,
 \]
 where \ $\lfloor x \rfloor$ \ and \ $([\log S]_t)_{t\in[0,\infty)}$ \ denotes the
 integer part of a real number \ $x \in \RR$, \ and the quadratic variation process of
 the semimartingale \ $\log S$, \ respectively.
\ By Theorem I.4.52 in Jacod and Shiryaev \cite{JSh},
 \[
   [\log S]_t = \langle (\log S)^\cont \rangle_t
   + \sum_{u\in[0,t]} (\log S_u - \log S_{u-})^2 , \qquad t \in [0, \infty) .
 \]
Consequently, for all \ $t \in [0, \infty)$, \ we have
 \[
   \sum_{i=1}^{\lfloor nt\rfloor} (\log S_{\frac{i}{n}} - \log S_{\frac{i-1}{n}})^2
   - \sum_{u\in[0,t]} (\log S_u - \log S_{u-})^2
   \stoch \langle (\log S)^\cont \rangle_t \qquad \text{as \ $n \to \infty$.}
 \]
Note that this convergence holds almost surely along a suitable subsequence,
 the members of this sequence are measurable functions of \ $(S_u)_{u\in[0,t]}$,
 \ hence, using Theorems 4.2.2 and 4.2.8 in Dudley \cite{Dud}, we obtain that
 \ $\langle (\log S)^\cont \rangle_t = \int_0^t Y_u \, \dd u$ \ is a measurable
 function (i.e., a statistic) of \ $(S_u)_{u\in[0,t]}$.
\ Moreover,
 \begin{equation}\label{QV}
  \frac{\langle (\log S)^\cont \rangle_{t+h}
        - \langle (\log S)^\cont \rangle_t}
       {h}
  = \frac{1}{h} \int_t^{t+h} Y_s \, \dd s
  \as Y_t  \qquad \text{as \ $h \to 0$,} \quad t \in [0, \infty) ,
 \end{equation}
 since \ $Y$ \ has continuous sample paths almost surely.
Consequently, for all \ $t \in [0, T]$, \ $Y_t$ \ is a measurable function
 (i.e., a statistic) of \ $(S_u)_{u\in[0,T]}$ \ (where for \ $t = T$, \ one may take
 \ $h \uparrow 0$), \ however, we also point out that this measurable function
 remains inexplicit.
\proofend
\end{Rem}

Next we give statistics for the parameters \ $\sigma$ \ and \ $\varrho$ \ using
 continuous time observations \ $(S_t)_{t\in[0,T]}$ \ with some \ $T > 0$.
\ Due to this result we do not consider the estimation of these parameters,
 they are supposed to be known.

\begin{Rem}\label{Thm_MLE_cons_sigma_rho}
Let \ $\theta, \kappa, \sigma \in (0, \infty)$, \ $\mu \in \RR$,
 \ $\varrho \in (-1, 1)$, \ and
 \ $\PP(Y_0 \in [0, \infty), \, S_0 \in (0, \infty)) = 1$.
\ Then for all \ $T > 0$,
 \[
   \begin{bmatrix}
    \sigma^2 & \varrho \sigma \\
    \varrho \sigma & 1
   \end{bmatrix}
   = \frac{1}{\int_0^T Y_s \, \dd s}
     \begin{bmatrix}
      \langle Y \rangle_T & \langle Y, (\log S)^\cont \rangle_T \\
      \langle Y, (\log S)^\cont \rangle_T
       & \langle (\log S)^\cont \rangle_T
     \end{bmatrix}
   =: \hbSigma_T ,
 \]
 where \ $(\langle Y, (\log S)^\cont \rangle_t)_{t\in[0,\infty)}$ \ denotes the
 (predictable) quadratic covariation process of \ $Y$ \ and \ $(\log S)^\cont$,
 \ since, by the SDEs
 \eqref{j_Heston_SDE} and \eqref{logS_SDE},
 \begin{gather*}
  \langle Y \rangle_T
   = \sigma^2 \int_0^T Y_u \, \dd u , \qquad
  \langle (\log S)^\cont \rangle_T
  = \int_0^T Y_u \, \dd u , \qquad
  \langle Y, (\log S)^\cont \rangle_T
  = \varrho \sigma \int_0^T Y_u \, \dd u .
 \end{gather*}
We point out that \ $\PP\bigl(\int_0^T Y_u \, \dd u \in (0, \infty)\bigr) = 1$.
\ Indeed, if \ $\omega \in \Omega$ \ is such that
 \ $[0, T] \ni s \mapsto Y_s(\omega)$ \ is continuous and
 \ $Y_t(\omega) \in [0, \infty)$ \ for all \ $t \in [0, \infty)$, \ then we have
 \ $\int_0^T Y_u(\omega)\,\dd u = 0$ \ if and only if \ $Y_u(\omega) = 0$ \ for
 all \ $u \in [0, T]$.
\ Using the method of the proof of Theorem 3.1 in Barczy et.\ al
 \cite{BarDorLiPap}, we get \ $\PP(\int_0^T Y_u \, \dd u = 0) = 0$, \ as desired.
We note that \ $\hbSigma_T$ \ is a statistic, i.e., there exists a measurable
 function \ $\Xi : D([0,T], \RR) \to \RR^{2\times2}$ \ such that
 \ $\hbSigma_T = \Xi((S_u)_{u\in[0,T]})$, \ where \ $D([0,T], \RR)$ \ denotes the
 space of real-valued c\`adl\`ag functions defined on \ $[0,T]$, \ since
 \begin{equation}\label{sigma_rho}
  \begin{aligned}
   &\frac{1}{\frac{1}{n} \sum_{i=1}^{\lfloor nT\rfloor} Y_{\frac{i-1}{n}}}
    \sum_{i=1}^{\lfloor nT\rfloor}
     \begin{bmatrix}
      Y_{\frac{i}{n}} - Y_{\frac{i-1}{n}} \\
      \log S_{\frac{i}{n}} - \log S_{\frac{i-1}{n}}
     \end{bmatrix}
     \begin{bmatrix}
      Y_{\frac{i}{n}} - Y_{\frac{i-1}{n}} \\
      \log S_{\frac{i}{n}} - \log S_{\frac{i-1}{n}}
     \end{bmatrix}^\top \\
   &- \frac{1}{\frac{1}{n} \sum_{i=1}^{\lfloor nT\rfloor} Y_{\frac{i-1}{n}}}
      \sum_{u\in[0,T]}
      \begin{bmatrix}
       0 & 0 \\
       0 & (\log S_u - \log S_{u-})^2
      \end{bmatrix}
    \stoch
    \hbSigma_T \qquad \text{as \ $n \to \infty$,}
  \end{aligned}
 \end{equation}
 where the convergence in \eqref{sigma_rho} holds almost surely along a suitable
 subsequence, the members of the sequence in
 \eqref{sigma_rho} are measurable functions of \ $(S_u)_{u\in[0,T]}$ \ (due to Remark \ref{observation}), and
 one can use Theorems 4.2.2 and 4.2.8 in Dudley \cite{Dud}.
\ Next we prove \eqref{sigma_rho}.
By Theorems I.4.47 a) and I.4.52 in Jacod and Shiryaev \cite{JSh},
 \begin{gather*}
  \sum_{i=1}^{\lfloor nT\rfloor} (Y_{\frac{i}{n}} - Y_{\frac{i-1}{n}})^2
  \stoch [Y]_T = \langle Y \rangle_T , \\
  \sum_{i=1}^{\lfloor nT\rfloor} (\log S_{\frac{i}{n}} - \log S_{\frac{i-1}{n}})^2
  \stoch [\log S]_T
  = \langle (\log S)^\cont \rangle_T
    + \sum_{u\in[0,T]} (\log S_u - \log S_{u-})^2 , \\
  \sum_{i=1}^{\lfloor nT\rfloor}
   (Y_{\frac{i}{n}} - Y_{\frac{i-1}{n}}) (\log S_{\frac{i}{n}} - \log S_{\frac{i-1}{n}})
  \stoch [Y, \log S]_T = \langle Y, (\log S)^\cont \rangle_T
 \end{gather*}
 as \ $n \to \infty$, \ where \ $([Y, \log S]_t)_{t\in[0,\infty)}$ \ denotes the
 quadratic covariation process of  the semimartingales \ $Y$ \ and \ $\log S$.
\ Consequently,
 \begin{align*}
  &\sum_{i=1}^{\lfloor nT\rfloor}
    \begin{bmatrix}
     Y_{\frac{i}{n}} - Y_{\frac{i-1}{n}} \\
     \log S_{\frac{i}{n}} - \log S_{\frac{i-1}{n}}
    \end{bmatrix}
    \begin{bmatrix}
     Y_{\frac{i}{n}} - Y_{\frac{i-1}{n}} \\
     \log S_{\frac{i}{n}} - \log S_{\frac{i-1}{n}}
    \end{bmatrix}^\top \\
  &- \sum_{u\in[0,T]}
      \begin{bmatrix}
       0 & 0 \\
       0 & (\log S_u - \log S_{u-})^2
      \end{bmatrix}
   \stoch
   \left( \int_0^T Y_u \, \dd u \right)
   \hbSigma_T \qquad \text{as \ $n \to \infty$,}
 \end{align*}
 see, e.g., van der Vaart \cite[Theorem 2.7, part (vi)]{Vaart}.
Moreover,
 \[
   \frac{1}{n} \sum_{i=1}^{\lfloor nT\rfloor} Y_{\frac{i-1}{n}}
   \as \int_0^T Y_u \, \dd u  \qquad \text{as \ $n \to \infty$}
 \]
 since \ $Y$ \ has continuous sample paths almost surely.
Hence \eqref{sigma_rho} follows by Slutsky's lemma.
Finally, we note that the sample size \ $T$ \ is fixed above, and it is enough to
 know any short sample \ $(S_u)_{u\in[0,T]}$ \ to carry out the above calculations.
\proofend
\end{Rem}

\section{Existence and uniqueness of MLE}
\label{section_EUMLE}

From this section, we will consider the jump-type Heston model \eqref{j_Heston_SDE}
 with known \ $\sigma \in (0, \infty)$, \ $\varrho \in (-1, 1)$,
 \ $\gamma \in \RR$, \ L\'evy measure \ $m$, \ and deterministic
 initial value \ $(Y_0, S_0) = (y_0, s_0) \in (0, \infty)^2$, \ and we will consider
 \ $\bpsi := (\theta, \kappa, \mu) \in (0, \infty)^2 \times \RR =: \Psi$ \ as a parameter.

Let \ $\PP_\bpsi$ \ denote the probability measure induced by
 \ $(Y_t, S_t)_{t\in[0,\infty)}$ \ on the measurable space
 \ $(D([0,\infty), \RR^2), \cD([0,\infty), \RR^2))$ \ of \ $\RR^2$-valued c\`adl\`ag
 functions defined on \ $[0, \infty)$ \ endowed with a right continuous
 filtration \ $(\cD_t([0,\infty), \RR^2))_{t\in[0,\infty)}$, \ see Appendix
 \ref{App_LR}.
Further, for all \ $T \in (0, \infty)$, \ let
 \ $\PP_{\bpsi,T} := \PP_{\bpsi\!\!}|_{\cD_T([0,\infty), \RR^2)}$ \ be the
 restriction of \ $\PP_\bpsi$ \ to \ $\cD_T([0,\infty), \RR^2)$.

Let us write the Heston model \eqref{j_Heston_SDE} in the form
 \begin{equation}\label{j_Heston_SDE_matrix}
  \begin{bmatrix} \dd Y_t \\ \dd S_t \end{bmatrix}
  = A(Y_t,S_t) H(\bpsi) \, \dd t
    + \Gamma(Y_t,S_t) \begin{bmatrix} \dd W_t \\ \dd B_t \end{bmatrix}
    + \begin{bmatrix} 0 \\ S_{t-} \, \dd L_t \end{bmatrix} ,
  \qquad t \in [0, \infty) ,
 \end{equation}
 where the functions \ $A : [0, \infty) \times (0, \infty) \to \RR^{2\times3}$,
 \ $\Gamma : [0, \infty) \times (0, \infty) \to \RR^{2\times2}$
 \ and \ $H : \RR^3 \to \RR^3$ \ are defined by
 \begin{gather*}
  A(y, s)
  := \begin{bmatrix}
      1 & -y & 0 \\
      0 & 0 & s
     \end{bmatrix} , \qquad
  \Gamma(y, s)
  := \sqrt{y}
     \begin{bmatrix}
      \sigma & 0 \\
      \varrho s & \sqrt{1 - \varrho^2} s
     \end{bmatrix} , \qquad
  H(x_1, x_2, x_3)
  := \begin{bmatrix}
      x_1 x_2 \\
      x_2 \\
      x_3
     \end{bmatrix}
 \end{gather*}
 for \ $(y, s) \in [0, \infty) \times (0, \infty)$ \ and \ $(x_1, x_2, x_3) \in \RR^3$.
\ Note that \ $H$ \ is bijective on the set
 \ $\RR \times (\RR \setminus \{0\}) \times \RR$ \ having inverse
 \begin{align}\label{H_inverse}
   H^{-1}(y_1, y_2, y_3)
   = \left( \frac{y_1}{y_2}, y_2, y_3 \right) , \qquad
   (y_1, y_2, y_3) \in \RR \times (\RR \setminus \{0\}) \times \RR .
 \end{align}
Let us introduce the function
 \ $\Sigma : [0, \infty) \times (0, \infty) \to \RR^{2\times2}$ \ given by
 \[
   \Sigma(y, s)
   := \Gamma(y, s) \Gamma(y, s)^\top
   = y
     \begin{bmatrix}
      \sigma^2 & \varrho \sigma s \\
      \varrho \sigma s & s^2
     \end{bmatrix} , \qquad
   (y, s) \in [0, \infty) \times (0, \infty) .
 \]
If \ $(y, s) \in (0, \infty)^2$ \ then \ $\Sigma(y, s)$ \ is invertible, namely,
 \begin{equation}\label{bSigma^{-1}}
  \begin{aligned}
   \Sigma(y, s)^{-1}
     = (\Gamma(y, s)^\top)^{-1} \Gamma(y, s)^{-1}
   & = \frac{1}{(1-\varrho^2) \sigma^2 y s^2}
      \begin{bmatrix}
       \sqrt{1 - \varrho^2} s & -\varrho s \\
       0 & \sigma
      \end{bmatrix}
      \begin{bmatrix}
       \sqrt{1 - \varrho^2} s & 0 \\
       -\varrho s & \sigma
      \end{bmatrix} \\
   & = \frac{1}{(1-\varrho^2) \sigma^2 y s^2}
       \begin{bmatrix}
       s^2 & -\varrho\sigma s \\
       -\varrho\sigma s & \sigma^2
      \end{bmatrix} .
  \end{aligned}
 \end{equation}
Further, let
 \begin{gather*}
  \bG_t
  :=\int_0^t A(Y_u, S_u)^\top \Sigma(Y_u, S_u)^{-1} A(Y_u, S_u) \, \dd u , \qquad
  t \in [0, \infty) , \\
  \Bf_t
  :=\int_0^t
     A(Y_{u-}, S_{u-})^\top
     \Sigma(Y_{u-}, S_{u-})^{-1}
    \begin{bmatrix}
     \dd Y_u \\
     \dd S_u - S_{u-} \, \dd L_u
    \end{bmatrix} , \qquad
  t \in [0, \infty) ,
 \end{gather*}
 provided that \ $\PP(\text{$Y_t \in (0, \infty)$ for all $t \in [0, \infty)$}) = 1$, \ which
 holds if \ $\theta \kappa \in \bigl[ \frac{\sigma^2}{2}, \infty \bigr)$.
\ Using \eqref{bSigma^{-1}}, we obtain
 \begin{align}\label{bGt}
  \begin{aligned}
   \bG_t
   &= \int_0^t
       A(Y_u, S_u)^\top (\Gamma(Y_u, S_u)^\top)^{-1}
       \Gamma(Y_u, S_u)^{-1} A(Y_u, S_u) \, \dd u \\
   &= \int_0^t
       (\Gamma(Y_u, S_u)^{-1} A(Y_u, S_u))^\top
       (\Gamma(Y_u, S_u)^{-1} A(Y_u, S_u)) \, \dd u \\
   &= \frac{1}{(1-\varrho^2)\sigma^2}
      \int_0^t
       \frac{1}{Y_u}
       \begin{bmatrix}
        1 & -Y_u & -\varrho \sigma \\
        -Y_u & Y_u^2 & \varrho \sigma Y_u \\
        -\varrho \sigma & \varrho \sigma Y_u & \sigma^2
       \end{bmatrix}
       \dd u , \qquad
   t \in [0, \infty) ,
  \end{aligned}
 \end{align}
 provided that \ $\PP(\text{$Y_t \in (0, \infty)$ for all $t \in [0, \infty)$}) = 1$, \ since
 \begin{equation}\label{bB^{-1}bA}
   \Gamma(y, s)^{-1} A(y, s)
   = \frac{1}{\sigma\sqrt{(1-\varrho^2)y}}
     \begin{bmatrix}
      \sqrt{1-\varrho^2} & -y\sqrt{1-\varrho^2} & 0 \\
      -\varrho & \varrho y & \sigma
      \end{bmatrix} , \qquad (y, s) \in (0, \infty)^2 .
 \end{equation}

The next lemma is about the form of the Radon--Nikodym derivative
 \ $\frac{\dd \PP_{\bpsi,T}}{\dd \PP_{\btpsi,T}}$ \ for certain
 \ $\bpsi, \btpsi \in \Psi$.

\begin{Lem}\label{RN}
Let \ $\bpsi = (\theta, \kappa, \mu) \in \Psi$ \ and
 \ $\btpsi := (\ttheta, \tkappa, \tmu) \in \Psi$ \ with
 \ $\theta \kappa, \ttheta \tkappa \in \bigl[ \frac{\sigma^2}{2}, \infty \bigr)$.
\ Then for all \ $T \in (0, \infty)$, \ the probability measures \ $\PP_{\bpsi,T}$ \ and
 \ $\PP_{\btpsi,T}$ \ are absolutely continuous with respect to each other, and,
 under \ $\PP$,
 \begin{equation}\label{RNformula}
  \log \frac{\dd \PP_{\bpsi,T}}{\dd \PP_{\btpsi,T}}(\tY, \tS)
  = \bigl(H(\bpsi) - H(\btpsi)\bigr)^\top \tBf_T
    - \frac{1}{2}
      \bigl(H(\bpsi) - H(\btpsi)\bigr)^\top \tbG_T \bigl(H(\bpsi) + H(\btpsi)\bigr) ,
 \end{equation}
 where \ $\tY$, \ $\tS$, \ $\tbG$ \ and \ $\tBf$ \ are the processes
 corresponding to the parameter \ $\btpsi$.
\end{Lem}

\noindent{\bf Proof.}
In what follows, we will apply Theorem III.5.34 in Jacod and Shiryaev \cite{JSh}
 (see also Appendix \ref{App_LR}).
We will work on the canonical space \ $(D([0,\infty), \RR^2), \cD([0,\infty), \RR^2))$.
\ Let \ $(\eta_t, \zeta_t)_{t\in[0,\infty)}$ \ denote the canonical process
 \ $(\eta_t, \zeta_t)(\omega) := \omega(t)$, \ $\omega \in D([0,\infty), \RR^2)$,
 \ $t \in [0, \infty)$.
Using \eqref{j_Heston_SDE_matrix} and \eqref{Levy_Ito}, the Heston model
 \eqref{j_Heston_SDE} can be written in the form
 \begin{align*}
  \begin{bmatrix} Y_t \\ S_t \end{bmatrix}
   & = \begin{bmatrix} y_0 \\ s_0 \end{bmatrix}
      + \int_0^t
         \left(A(Y_u, S_u) H(\bpsi)
               + \begin{bmatrix}
                  0 \\
                  \gamma S_u
                 \end{bmatrix}\right)
         \dd u \\
   &\quad
      + \int_0^t \Gamma(Y_u, S_u) \begin{bmatrix} \dd W_u \\ \dd B_u \end{bmatrix}
      + \begin{bmatrix}
         0 \\
         \int_0^t \int_{\RR}
          S_{u-}h_1(z) (\mu^L(\dd u,\dd z) - \dd u \, m(\dd z))
        \end{bmatrix} \\
   &\quad
      + \begin{bmatrix}
         0 \\
         \int_0^t \int_{\RR} S_{u-}(z - h_1(z)) \, \mu^L(\dd u,\dd z)
        \end{bmatrix} ,
      \qquad  t\in[0,\infty).
 \end{align*}
Using Proposition \ref{Grigelionis}, we obtain
 \begin{align*}
  \begin{bmatrix} Y_t \\ S_t \end{bmatrix}
   & = \begin{bmatrix} y_0 \\ s_0 \end{bmatrix}
      + \int_0^t
         \left(A(Y_u, S_u) H(\bpsi)
               + \begin{bmatrix}
                  0 \\
                  \gamma S_u
                  + \int_{\RR} (h_1(S_{u-} z) - S_{u-} h_1(z)) \, m(\dd z)
                 \end{bmatrix}\right)
         \dd u \\
   &\quad
      + \int_0^t \Gamma(Y_u, S_u) \begin{bmatrix} \dd W_u \\ \dd B_u \end{bmatrix}
      + \begin{bmatrix}
         0 \\
         \int_0^t \int_{\RR}
          h_1(S_{u-} z) (\mu^L(\dd u,\dd z) - \dd u \, m(\dd z))
        \end{bmatrix} \\
   &\quad
      + \begin{bmatrix}
         0 \\
         \int_0^t \int_{\RR} (S_{u-} z - h_1(S_{u-} z)) \, \mu^L(\dd u,\dd z)
        \end{bmatrix} ,
      \qquad  t\in[0,\infty),
 \end{align*}
 which is a special case of III.2.23 in Jacod and Shiryaev \cite{JSh}, since
 \begin{align}\label{help_RN_atiras}
  \begin{split}
  \begin{bmatrix} Y_t \\ S_t \end{bmatrix}
   &= \begin{bmatrix} y_0 \\ s_0 \end{bmatrix}
      + \int_0^t
         \left(A(Y_u,S_u) H(\bpsi)
               + \begin{bmatrix}
                  0 \\
                  \gamma S_u
                  + \int_{\RR} (h_1(S_{u-} z) - S_{u-} h_1(z)) \, m(\dd z)
                 \end{bmatrix}\right)
         \dd u \\
   &\quad
      + \int_0^t \Gamma(Y_u,S_u) \begin{bmatrix} \dd W_u \\ \dd B_u \end{bmatrix}
      + \int_0^t \int_{\RR}
         h_2\left(\begin{bmatrix} 0 \\ S_{u-} z \end{bmatrix}\right)
         (\mu^L(\dd u,\dd z) - \dd u \, m(\dd z)) \\
   &\quad
      + \int_0^t \int_{\RR}
         \left(\begin{bmatrix} 0 \\ S_{u-} z \end{bmatrix}
               - h_2\left(\begin{bmatrix} 0 \\ S_{u-} z \end{bmatrix}\right)\right)
         \mu^L(\dd u,\dd z) , \qquad t \in [0, \infty) ,
  \end{split}
 \end{align}
 with the truncation function \ $h_2(\bz) := \bz \bbone_{[-1,1]^2}(\bz)$,
 \ $\bz \in \RR^2$, \ where we used that
 \[
   \begin{bmatrix} 0 \\ h_1(z) \end{bmatrix}
   = h_2\left(\begin{bmatrix} 0 \\ z \end{bmatrix}\right) , \qquad z \in \RR .
 \]
By Proposition \ref{Pro_j_Heston}, the SDE \eqref{j_Heston_SDE} has a pathwise
 unique strong solution (with the given deterministic initial value
 \ $(y_0,s_0)\in(0,\infty)^2$), \ and hence, by Theorem III.2.26 in Jacod and
 Shiryaev \cite{JSh}, under the probability measure \ $\PP_\bpsi$, \ the canonical
 process \ $(\eta_t, \zeta_t)_{t\in[0,\infty)}$ \ is a semimartingale with
 semimartingale characteristics \ $(B^{(\bpsi)}, C, \nu)$ \ associated with the
 truncation function \ $h_2$, \ where
 \begin{align*}
  &B^{(\bpsi)}_t
   =  \int_0^t
         \left(A(\eta_u, \zeta_u) H(\bpsi)
               + \begin{bmatrix}
                  0 \\
                  \gamma \zeta_u
                  + \int_{\RR} (h_1(\zeta_u z) - \zeta_u h_1(z)) \, m(\dd z)
                 \end{bmatrix}\right)
         \dd u , \\
  &C_t = \int_0^t \Gamma(\eta_u, \zeta_u) \Gamma(\eta_u, \zeta_u)^\top \, \dd u
       = \int_0^t \Sigma(\eta_u, \zeta_u) \, \dd u ,
 \end{align*}
 for \ $t \in [0, \infty)$, \ and
 \begin{align*}
  \nu(\dd t, \dd y, \dd z) = K(\eta_t, \zeta_t, \dd y, \dd z) \, \dd t
 \end{align*}
 with the Borel transition kernel \ $K$ \ from \ $[0, \infty)^2 \times \RR^2$
 \ into \ $\RR^2$ \ given by
 \[
   K(y, s, R) := \int_{\RR} \!\bbone_{R\setminus\{\bzero\}}(0, s z) \, m(\dd z)
   \qquad \text{for \ $(y, s) \in [0, \infty)^2$ \ and \ $R \in \cB(\RR^2)$.}
 \]
The aim of the following discussion is to check the set of sufficient conditions
 presented in Appendix \ref{App_LR} (of which the notations will be used) in
 order to have right to apply Theorem III.5.34 in Jacod and Shiryaev \cite{JSh}.
First note that \ $(C_t)_{t\in[0,\infty)}$ \ and \ $\nu(\dd t, \dd y, \dd z)$ \ do not
 depend on the unknown parameter \ $\bpsi$, \ and hence \ $V^{(\bpsi,\btpsi)}$ \ is identically one and then \eqref{GIR1} and \eqref{GIR2} readily hold.
We also have
 \[
   \PP_\bpsi\big(\nu( \{t\}\times \RR^2) = 0 \big)
     = \PP_\bpsi \left( \int_{\{t\}} K(\eta_s, \zeta_s, \RR^2) \,\dd s = 0\right)
     = 1,
     \qquad t \in [0, \infty) , \quad \bpsi \in \Psi .
 \]
Further, \ $(C_t)_{t\in[0,\infty)}$ \ can be represented as
 \ $C_t = \int_0^t c_u \, \dd F_u$, \ $t \in [0, \infty)$, \ where the stochastic
 processes \ $(c_t)_{t\in[0,\infty)}$ \ and \ $(F_t)_{t\in[0,\infty)}$ \ are given by
 \ $c_t := \Sigma(\eta_t, \zeta_t)$, \ $t \in [0, \infty)$, \ and
 \ $F_t = t$, \ $t \in [0, \infty)$.
\ Next, note that, under the condition
 \ $\theta\kappa \in \bigl[\frac{\sigma^2}{2}, \infty\bigr)$, \ we have
 \ $\PP_\bpsi(\text{$(\eta_t, \zeta_t) \in (0, \infty)^2$ for all $t \in [0, \infty)$}) = 1$
 \ (due to Proposition \ref{Pro_j_Heston}), hence, by \eqref{bSigma^{-1}}, for each
 \ $t \in [0, \infty)$, \ the matrix \ $c_t$ \ is invertible \ $\PP_\bpsi$-almost surely.
\ Consequently, for all \ $\bpsi = (\theta, \kappa, \mu) \in \Psi$ \ and
 \ $\btpsi = (\ttheta, \tkappa, \tmu) \in \Psi$ \ with
 \ $\theta \kappa, \ttheta \tkappa \in \bigl[ \frac{\sigma^2}{2}, \infty \bigr)$,
 \begin{align*}
  B^{(\bpsi)}_t - B^{(\btpsi)}_t
  = \int_0^t A(\eta_u, \zeta_u) (H(\bpsi) - H(\btpsi)) \, \dd u
  = \int_0^t c_u \beta^{(\btpsi,\bpsi)}_u \, \dd F_u ,
  \qquad t \in [0, \infty) ,
 \end{align*}
 where the stochastic process \ $(\beta^{(\btpsi,\bpsi)}_t)_{t\in[0,\infty)}$ \ is given by
 \[
   \beta^{(\btpsi,\bpsi)}_t
   = c_t^{-1} A(\eta_t,\zeta_t) (H(\bpsi) - H(\btpsi))
   = \Sigma(\eta_t, \zeta_t)^{-1} A(\eta_t, \zeta_t) (H(\bpsi) - H(\btpsi)) ,
   \qquad t \in [0, \infty) ,
 \]
 which yields \eqref{GIR3}.

Next we check \eqref{GIR4}, i.e.,
 \begin{gather}\label{COND}
  \PP_\bpsi\left(\int_0^t
                 \bigl(\beta^{(\btpsi,\bpsi)}_u\bigr)^\top
                 c_u
                 \beta^{(\btpsi,\bpsi)}_u
                 \, \dd F_u
           < \infty\right)
  = 1,  \qquad t\in[0,\infty).
 \end{gather}
We have
 \begin{align*}
   &\int_0^t (\beta^{(\btpsi,\bpsi)}_s)^\top c_s \beta^{(\btpsi,\bpsi)}_s \, \dd F_s \\
   &\qquad =  (H(\bpsi) - H(\btpsi))^\top
        \int_0^t A(\eta_s,\zeta_s)^\top  (\Sigma(\eta_s, \zeta_s)^{-1})^\top \,
                 A(\eta_s,\zeta_s) \,\dd s\,
           (H(\bpsi) - H(\btpsi)) \\
   & \qquad =   (H(\bpsi) - H(\btpsi))^\top  G_t  (H(\bpsi) - H(\btpsi)),
   \qquad t\in[0,\infty) ,
 \end{align*}
 where \ $G_t$, \ $t \in [0, \infty)$, \ is understood as the original \ $\bG_t$,
 \ $t \in [0, \infty)$, \ replacing \ $(Y, S)$ \ by \ $(\eta, \zeta)$.
\ Since \ $\eta$ \ has continuous sample paths \ $\PP_\bpsi$-almost surely and
 \ $\PP_\bpsi(\eta_t \in (0, \infty), \, \forall \, t \in [0,\infty)) = 1$
 \ (due to \ $\theta\kappa\in[\frac{\sigma^2}{2},\infty)$),
 \ we have \ $\PP_\bpsi(\inf_{t\in[0,T]} \eta_t \in (0, \infty)) = 1$ \ for all
 \ $T \in [0, \infty)$, \ which, together with the \ $\PP_\bpsi$-almost sure
 continuity of \ $\eta$ \ and formula \eqref{bGt}, yield \eqref{COND}.

Next, we check that, under the probability measure \ $\PP_{\bpsi}$, \ local
 uniqueness holds for the martingale problem on the canonical space corresponding
 to the triplet \ $(B^{(\bpsi)}, C, \nu)$ \ with the given initial value
 \ $(y_0, s_0)$ \ with \ $\PP_\bpsi$ \ as its unique solution.
By Proposition \ref{Pro_j_Heston}, the SDE \eqref{j_Heston_SDE} has a pathwise
 unique strong solution (with the given deterministic initial value
 \ $(y_0,s_0)\in(0,\infty)^2$), \ and hence Theorem III.2.26 in Jacod and
 Shiryaev \cite{JSh} yields that the set of all solutions to the martingale
 problem on the canonical space corresponding to \ $(B^{(\bpsi)}, C, \nu)$ \ has
 only one element \ $(\PP_{\bpsi})$ \ yielding the desired local uniqueness.
We also mention that Theorem III.4.29 in Jacod and Shiryaev \cite{JSh} implies that under the probability
 measure \ $\PP_\bpsi$, \ all local martingales have the integral representation property
 relative to \ $(\eta, \zeta)$.

By Theorem III.5.34 in Jacod and Shiryaev \cite{JSh} (see also Appendix \ref{App_LR}),
 \ $\PP_{\bpsi,T}$ \ and \ $\PP_{\btpsi,T}$ \ are equivalent (one can change the roles of
 \ $\bpsi$ \ and \ $\btpsi$), \ and under the probability measure \ $\PP_{\btpsi,T}$, \  we have
 \begin{align*}
  \frac{\dd \PP_{\bpsi,T}}{\dd \PP_{\btpsi,T}}{(\eta, \zeta)}
  = \exp\bigg\{\int_0^T
                \bigl(\beta^{(\btpsi,\bpsi)}_u\bigr)^\top
                \begin{bmatrix}
                 \dd (\eta^\cont)^{(\btpsi)}_u \\
                 \dd (\zeta^\cont)^{(\btpsi)}_u
                \end{bmatrix}
               -\frac{1}{2}
                \int_0^T
                 \bigl(\beta^{(\btpsi,\bpsi)}_u\bigr)^\top
                 c_u
                 \beta^{(\btpsi,\bpsi)}_u
                 \, \dd u
        \bigg\}
 \end{align*}
 for \ $T \in (0, \infty)$, \ where
 \ $((\eta^\cont)^{(\btpsi)}_t, (\zeta^\cont)^{(\btpsi)}_t)_{t\in[0,\infty)}$
 \ denotes the continuous (local) martingale part of
 \ $(\eta_t,\zeta_t)_{t\in[0,\infty)}$ \ under \ $\PP_{\btpsi}$.
\ Using part 1) of Remarks III.2.28 in Jacod and Shiryaev \cite{JSh} and
 \eqref{help_RN_atiras}, the continuous (local) martingale part
 \ $(\tY^\cont_t,\tS^\cont_t)_{t\in[0,\infty)}$ \ of
 \ $(\tY_t,\tS_t)_{t\in[0,\infty)}$ \ takes the form
 \begin{align*}
    \begin{bmatrix} \tY^\cont_t \\ \tS^\cont_t \end{bmatrix}
      = \int_0^t \Gamma(\tY_u,\tS_u)  \begin{bmatrix} \dd W_u \\ \dd B_u \end{bmatrix},
     \qquad t\in[0,\infty),
 \end{align*}
 and, by \eqref{j_Heston_SDE_matrix}, we have
 \begin{equation}\label{clmp}
   \begin{bmatrix} \dd \tY^\cont_t \\ \dd \tS^\cont_t \end{bmatrix}
   = \begin{bmatrix} \dd \tY_t \\ \dd \tS_t \end{bmatrix}
     - A(\tY_t, \tS_t) H(\btpsi) \, \dd t
     - \begin{bmatrix} 0 \\ \tS_{t-} \, \dd L_t \end{bmatrix} , \qquad t \in [0, \infty) .
 \end{equation}
Hence, under \ $\PP$,
 \begin{align*}
  &\log\,\frac{\dd \PP_{\bpsi,T}}{\dd \PP_{\btpsi,T}}{(\tY, \tS)}
  = (H(\bpsi) - H(\btpsi))^\top
     \int_0^T A(\tY_u, \tS_{u-})^\top
      \Sigma(\tY_u, \tS_{u-})^{-1}
      \begin{bmatrix} \dd \tY_u \\ \dd \tS_u - \tS_{u-} \, \dd L_u \end{bmatrix} \\
  &- (H(\bpsi) - H(\btpsi))^\top
       \left( \int_0^T
               A(\tY_u, \tS_u)^\top \Sigma(\tY_u, \tS_u)^{-1}
               A(\tY_u, \tS_u) \, \dd u \right)
       H(\btpsi) \\
  &- \frac{1}{2} (H(\bpsi) - H(\btpsi))^\top
       \left(\int_0^T
              A(\tY_u, \tS_u)^\top \Sigma(\tY_u, \tS_u)^{-1}
              A(\tY_u, \tS_u) \, \dd u \right)
       (H(\bpsi) - H(\btpsi)) \\
  & = (H(\bpsi) - H(\btpsi))^\top \tBf_T - (H(\bpsi) - H(\btpsi))^\top \tbG_T  H(\btpsi) \\
  &\phantom{=\;} -\frac{1}{2}(H(\bpsi) - H(\btpsi))^\top \tbG_T (H(\bpsi) - H(\btpsi)),
 \end{align*}
 which yields the statement.
\proofend

Note that \ $\tBf_T$ \ in Lemma \ref{RN} contains a stochastic integral with respect to \ $L$,
 \ but, by Remark \ref{observationL}, for all \ $t\in[0,T]$, \ $L_t$ \ is a measurable
 function of \ $(\tS_t)_{t\in[0,T]}$ \ (depending on \ $\gamma$ \ and \ $m$).

We point out that we use the condition \ $\theta\kappa\in[\frac{\sigma^2}{2},\infty)$ \
 in the proof of Lemma \ref{RN} to assure the invertibility of
 \ $(c_t)_{t\in[0,\infty)}$.

Next, using Lemma \ref{RN}, by considering \ $\PP_{\btpsi,T}$ \ as a fixed reference measure,
 we will derive an MLE for the parameter \ $\bpsi = (\theta, \kappa, \mu)$ \ based on the
 observations \ $(Y_t, S_t)_{t\in[0,T]}$.
\ Our method for deriving an MLE is one of the known ones in the literature, and it
 turns out that these lead to the same estimator \ $\hbpsi_T$, \ see Remark \ref{Luschgy}.
Let us denote the right hand side of
 \eqref{RNformula} by \ $\Lambda_T(\bpsi, \btpsi)$ \ replacing \ $(\tBf_T, \tbG_T)$ \ by
 \ $(\Bf_T, \bG_T)$.
\ For convenience, first we calculate an MLE \ $\hbpsi_T$ \ of the parameter \ $\bpsi$ \ on the set \ $\RR^3$ \ based on the
 observations \ $(Y_t, S_t)_{t\in[0,T]}$, \ namely,
 \[
   \hbpsi_T := \argmax_{\bpsi\in\RR^3} \Lambda_T(\bpsi, \btpsi) ,
 \]
 which will turn out to be not dependent on \ $\btpsi$.
\ Here the function \ $\Lambda_T$ \ is extended for all
 \ $\bpsi = (\theta, \kappa, \mu) \in \RR^3$ \ in a natural way
 (note that for the calculation of the random matrices \ $\bG_t$,
 \ $t \in [0, \infty)$, \ and the random vectors \ $\Bf_t$, \ $t \in [0, \infty)$,
 \ one does not need to know the parameters \ $\bpsi$ \ or \ $\btpsi$).
\ In Remark \ref{REM_argmax_PSI}, we describe the connection between \ $\hbpsi_T$ \ and an MLE given by
 \ $\argmax_{\bpsi\in\Psi} \Lambda_T(\bpsi, \btpsi)$ \  on the set \ $\Psi$.

\begin{Pro}\label{Pro_MLE}
Let \ $\theta, \kappa \in (0, \infty)$ \ with \ $\theta \kappa \in \bigl[ \frac{\sigma^2}{2}, \infty \bigr)$,
 \ $\mu \in \RR$, \ $\sigma \in (0, \infty)$, \ $\varrho \in (-1, 1)$, \ and \ $(Y_0, S_0) = (y_0, s_0) \in (0, \infty)^2$.
\ Then for all \ $T \in (0, \infty)$, \ there exists a unique MLE \ $\hbpsi_T = (\htheta_T, \hkappa_T, \hmu_T)^\top$ \
 of \ $\bpsi = (\theta, \kappa, \mu)^\top$ \ on the set \ $\RR^3$ \ based on the observations
 \ $(Y_t, S_t)_{t\in[0,T]}$ \ taking the form
  \begin{equation}\label{MLE}
  \hbpsi_T
  = \begin{bmatrix}
     \htheta_T \\ \hkappa_T \\ \hmu_T
    \end{bmatrix}
  = \begin{bmatrix}
     \frac{(\bG_T^{-1} \Bf_T)_1}{(\bG_T^{-1} \Bf_T)_2} \\
     (\bG_T^{-1} \Bf_T)_2 \\
     (\bG_T^{-1} \Bf_T)_3
    \end{bmatrix} ,
 \end{equation}
 provided that \ $\bG_T$ \ is strictly positive definite and
 \ $(\bG_T^{-1} \Bf_T)_2 \ne 0$, \ which hold almost surely.
Further, we have
 \begin{align}\label{MLE_coordinatewise}
 \begin{split}
  \htheta_T
  &= \frac{\int_0^T Y_u \, \dd u \int_0^T \frac{\dd Y_u}{Y_u}
           - T \int_0^T \dd Y_u
           + \varrho \sigma T \int_0^T \frac{\dd S_u - S_{u-} \, \dd L_u}{S_{u-}}
           - \varrho \sigma T^2 \frac{\int_0^T \frac{\dd S_u - S_{u-} \, \dd L_u}{Y_u S_{u-}}}{\int_0^T \frac{\dd u}{Y_u}} }
          {T  \int_0^T \frac{\dd Y_u}{Y_u}
           - \int_0^T \frac{\dd u}{Y_u} \int_0^T \dd Y_u
           + \varrho \sigma \int_0^T \frac{\dd u}{Y_u}
             \int_0^T \frac{\dd S_u - S_{u-} \, \dd L_u}{S_{u-}}
           - \varrho \sigma T \int_0^T \frac{\dd S_u - S_{u-} \, \dd L_u}{Y_u S_{u-}}} , \\[1mm]
  \hkappa_T
  &= \frac{T \int_0^T \frac{\dd Y_u}{Y_u}
           - \int_0^T \frac{\dd u}{Y_u} \int_0^T \dd Y_u
           + \varrho \sigma \int_0^T \frac{\dd u}{Y_u}
             \int_0^T \frac{\dd S_u - S_{u-} \, \dd L_u}{S_{u-}}
           - \varrho \sigma T
             \int_0^T \frac{\dd S_u - S_{u-} \, \dd L_u}{Y_u S_{u-}}}
         {\int_0^T Y_u \, \dd u \int_0^T \frac{\dd u}{Y_u} - T^2 } , \\
  \hmu_T
  &= \frac{\int_0^T \frac{\dd S_u - S_{u-} \, \dd L_u}{Y_u S_{u-}}}
          {\int_0^T \frac{\dd u}{Y_u}} .
  \end{split}
 \end{align}
\end{Pro}

\noindent {\bf Proof.}
The function \ $\Lambda_T$ \ can be written in the form
 \begin{align*}
  \Lambda_T(\bpsi, \btpsi)
  = - \frac{1}{2} H(\bpsi)^\top \bG_T H(\bpsi) + H(\bpsi)^\top \Bf_T - \bc ,
  \qquad \bpsi \in \RR^3 ,
 \end{align*}
 with
 \begin{gather*}
  \bc := H(\btpsi)^\top \Bf_T - \frac{1}{2} H(\btpsi)^\top \bG_T H(\btpsi) ,
 \end{gather*}
 since the symmetry of \ $\Sigma(y, s)^{-1}$, \ $(y, s) \in (0, \infty)^2$,
 \ implies the symmetry of \ $\bG_T$, \ and hence
 \ $H(\btpsi)^\top \bG_T H(\bpsi) = H(\bpsi)^\top \bG_T H(\btpsi)$.
\ The symmetric random matrix \ $\bG_T$ \ is almost surely strictly positive
 definite, since its \ $k \times k$ \ minors, \ $k \in \{1, 2, 3\}$ \ (see, \eqref{bGt}),
 are almost surely positive, namely,
 \begin{gather*}
   \frac{1}{(1-\varrho^2)\sigma^2} \int_0^T \frac{\dd u}{Y_u} > 0 , \qquad
   \frac{1}{((1-\varrho^2)\sigma^2)^2}
   \biggl(\int_0^T Y_u \, \dd u \int_0^T \frac{\dd u}{Y_u} - T^2\biggr) > 0 , \\
   \frac{1}{((1-\varrho^2)\sigma^2)^2}
   \biggl(\int_0^T Y_u \, \dd u \int_0^T \frac{\dd u}{Y_u} - T^2\biggr)
   \int_0^T \frac{\dd u}{Y_u} > 0
 \end{gather*}
 almost surely.
Indeed, \ $\int_0^T \frac{\dd u}{Y_u} > 0$ \ a.s.\ follows from
 \ $\PP(\text{$Y_t \in (0, \infty)$ \ for all \ $t \in [0, \infty)$}) = 1$,
 \ which can be found, e.g., in Lamberton and Lapeyre
 \cite[Proposition 6.2.4]{LamLap} (see also Proposition \ref{Pro_j_Heston}),
 and the proof of \ $\int_0^T Y_u \, \dd u \int_0^T \frac{\dd u}{Y_u} > T^2$ \ a.s.\ is
 given, e.g., in Barczy and Pap \cite[Lemma 3.3]{BarPap}.
Thus the matrix \ $\bG_T$ \ is almost surely invertible, namely,
 \begin{align*}
  \bG_T^{-1}
  &=\frac{1}
         {\bigl(\int_0^T Y_u \, \dd u \int_0^T \frac{\dd u}{Y_u} - T^2\bigr)
          \int_0^T \frac{\dd u}{Y_u}} \\
  &\quad
    \times
    \begin{bmatrix}
     \sigma^2 \int_0^T Y_u \, \dd u \int_0^T \frac{\dd u}{Y_u}
     - \varrho^2 \sigma^2 T^2
      & (1 - \varrho^2) \sigma^2 T \int_0^T \frac{\dd u}{Y_u}
      & \varrho \sigma \int_0^T Y_u \, \dd u \int_0^T \frac{\dd u}{Y_u}
        - \varrho \sigma T^2 \\
     (1 - \varrho^2) \sigma^2 T \int_0^T \frac{\dd u}{Y_u}
      & (1 - \varrho^2) \sigma^2 \bigl(\int_0^T \frac{\dd u}{Y_u}\bigr)^2 & 0 \\
     \varrho \sigma \int_0^T Y_u \, \dd u \int_0^T \frac{\dd u}{Y_u}
     - \varrho \sigma T^2
      & 0 & \int_0^T Y_u \, \dd u \int_0^T \frac{\dd u}{Y_u} - T^2
    \end{bmatrix}
 \end{align*}
 whenever \ $\int_0^T Y_u \, \dd u \int_0^T \frac{\dd u}{Y_u} > T^2$ \ and
 \ $\int_0^T \frac{\dd u}{Y_u} > 0$, \ which hold almost surely.
Provided that \ $\bG_T$ \ is strictly positive definite, we have
 \begin{align*}
  \Lambda_T(\bpsi, \btpsi)
  = - \frac{1}{2} \bigl(H(\bpsi) - \bG_T^{-1} \Bf_T\bigr)^\top \bG_T
       \bigl(H(\bpsi) - \bG_T^{-1} \Bf_T\bigr)
     + \frac{1}{2} \Bf_T^\top \bG_T^{-1} \Bf_T
     - \bc
  \leq \frac{1}{2} \Bf_T^\top \bG_T^{-1} \Bf_T - \bc ,
 \end{align*}
 and equality holds if and only if \ $H(\bpsi) = \bG_T^{-1} \Bf_T$.
\ The aim of the following discussion is to show that the inverse mapping \ $H^{-1}$ \
 given in \eqref{H_inverse} can be applied to
 \ $\bG_T^{-1} \Bf_T =: ((\bG_T^{-1} \Bf_T)_1, (\bG_T^{-1} \Bf_T)_2, (\bG_T^{-1} \Bf_T)_3)$
 \ almost surely, that is, \ $\PP((\bG_T^{-1} \Bf_T)_2 = 0) = 0$.
\ Applying \eqref{bSigma^{-1}}, we obtain
 \[
   \Bf_T
   = \frac{1}{(1-\varrho^2)\sigma^2}
     \begin{bmatrix}
      \int_0^T \frac{\dd Y_u}{Y_u}
       - \varrho \sigma \int_0^T \frac{\dd S_u - S_{u-} \, \dd L_u}{Y_u S_{u-}} \\[1mm]
      - \int_0^T \dd Y_u
       + \varrho \sigma \int_0^T \frac{\dd S_u - S_{u-} \, \dd L_u}{S_{u-}} \\[1mm]
      - \varrho \sigma \int_0^T \frac{\dd Y_u}{Y_u}
       + \sigma^2 \int_0^T \frac{\dd S_u - S_{u-} \, \dd L_u}{Y_u S_{u-}}
     \end{bmatrix} , \qquad T \in (0, \infty) .
 \]
Using the explicit formula for \ $\bG_T^{-1}$, \ we obtain
 \begin{align*}
  (\bG_T^{-1} \Bf_T)_2
   = \frac{T \int_0^T \frac{\dd Y_u}{Y_u}
           - \int_0^T \frac{\dd u}{Y_u} \int_0^T \dd Y_u
           + \varrho \sigma \int_0^T \frac{\dd u}{Y_u}
             \int_0^T \frac{\dd S_u - S_{u-} \, \dd L_u}{S_{u-}}
           - \varrho \sigma T
             \int_0^T \frac{\dd S_u - S_{u-} \, \dd L_u}{Y_u S_{u-}}}
          {\int_0^T Y_u \, \dd u \int_0^T \frac{\dd u}{Y_u} - T^2 } .
 \end{align*}
By the SDE \eqref{j_Heston_SDE},
 \begin{align}\label{KEY}
  \begin{aligned}
   \int_0^T \frac{\dd S_u - S_{u-} \, \dd L_u}{S_{u-}}
   = \int_0^T \bigl[\mu \, \dd u
                    + \sqrt{Y_u}
                      \bigl(\varrho \, \dd W_u
                            + \sqrt{1 - \varrho^2} \, \dd B_u\bigr)\bigr] , \\
   \int_0^T \frac{\dd S_u - S_{u-} \, \dd L_u }{Y_u S_{u-}}
   = \int_0^T
      \frac{\mu \, \dd u
            + \sqrt{Y_u}
              \bigl(\varrho \, \dd W_u + \sqrt{1 - \varrho^2} \, \dd B_u\bigr)}
           {Y_u}.
  \end{aligned}
 \end{align}
Applying again the SDE \eqref{j_Heston_SDE}, we have
 \begin{align*}
   \int_0^T \sqrt{Y_u} \, \dd W_u
   = \frac{1}{\sigma}
     \left( Y_T - y_0 - \int_0^T \kappa (\theta - Y_u) \, \dd u \right) , \qquad
   \int_0^T \frac{\dd W_u}{\sqrt{Y_u}}
   = \frac{1}{\sigma}
     \int_0^T \frac{\dd Y_u - \kappa (\theta - Y_u) \, \dd u}{Y_u} ,
 \end{align*}
 where
 \[
   \int_0^T\frac{\dd Y_u}{Y_u}
      = \log(Y_T) - \log(y_0) + \frac{\sigma^2}{2}\int_0^T\frac{\dd u}{Y_u}.
 \]
Indeed, since \ $\PP(Y_t \in (0, \infty) \ \text{for all} \ t \in [0, \infty)) = 1$, \ one can apply
 It\^o's rule to the function \ $f(x) = \log(x)$, \ $x \in (0, \infty)$, \ for which
 \ $f'(x) = 1 / x$, \ $f''(x) = - 1 / x^2$, \ $x \in (0, \infty)$, \ and we obtain
 \begin{align}\label{help_ITO}
   \log(Y_T)
   = \log(y_0) + \int_0^T \frac{\dd Y_u}{Y_u}
     - \frac{\sigma^2}{2} \int_0^T \frac{\dd u}{Y_u} , \qquad T \in [0, \infty) ,
 \end{align}
 see von Weizs\"acker and Winkler \cite[Theorem 8.1.1]{WeiWin}.
Hence, using the independence of the processes \ $Y$ \ and \ $B$, \ the conditional
 distribution of \ $(\bG_T^{-1} \Bf_T)_2$ \ given \ $(Y_u)_{u\in[0,T]}$
 \ is Gaussian and hence absolutely continuous, implying
 \[
   \PP((\bG_T^{-1} \Bf_T)_2 = 0) = \EE(\PP((\bG_T^{-1} \Bf_T)_2 = 0 \mid (Y_u)_{u\in[0,T]}))
   = 0 .
 \]
Consequently,
 \[
   H^{-1}(\bG_T^{-1} \Bf_T)
   = \argmax_{\bpsi\in\RR^3} \Lambda_T(\bpsi, \btpsi) ,
 \]
 provided that \ $\bG_T$ \ is strictly positive definite and
 \ $(\bG_T^{-1} \Bf_T)_2 \ne 0$, \ which hold almost surely, hence there
 exists a unique MLE
 \ $\hbpsi_T = (\htheta_T, \hkappa_T, \hmu_T)^\top = H^{-1}(\bG_T^{-1} \Bf_T)$ \ of
 \ $\bpsi = (\theta, \kappa, \mu)^\top$ \ on the set \ $\RR^3$ \ based on the
 observations \ $(Y_t, S_t)_{t\in[0,T]}$ \ yielding \eqref{MLE}.
Using again the explicit formula for \ $\bG_T^{-1}$, \ we obtain \eqref{MLE_coordinatewise} as well.
Note that \ $\hbpsi_T$ \ is a measurable function of the observations
 \ $(Y_t, S_t)_{t\in[0,T]}$, \ since all the integrals appearing in \ $\hbpsi_T$ \ are measurable functions
 of this process.
Indeed, in Remark \ref{observationL} we showed that for all \ $t \in [0, T]$,
 \ $L_t$ \ is a measurable function of \ $(S_u)_{u\in[0,T]}$, \ and one can use
 the arguments of Remarks \ref{observation} and \ref{Thm_MLE_cons_sigma_rho}
 together with
 Proposition I.4.44 in Jacod and Shiryaev \cite{JSh}, and Theorems 4.2.2 and
 4.2.8 in Dudley \cite{Dud}.
For example, for all \ $T\in[0,\infty)$,
 \[
   \sum_{i=1}^{\lfloor nT\rfloor}
    \frac{(S_{\frac{i}{n}} - S_{\frac{i-1}{n}})
          - S_{\frac{i-1}{n}-} (L_{\frac{i}{n}} - L_{\frac{i-1}{n}})}
         {Y_{\frac{i-1}{n}} S_{\frac{i-1}{n}-}}
   \stoch
   \int_0^T \frac{\dd S_u - S_{u-} \, \dd L_u}{Y_u S_{u-}}
   \qquad \text{as \ $n \to \infty$.}
 \]
\proofend

\begin{Rem}\label{REM_argmax_PSI}
We call the attention that later on it will turn out that \ $\hbpsi_T$ \ is a
 weakly consistent estimator of \ $\bpsi$ \ (see, Theorem \ref{Thm_MLE_cons} and
 Remark \ref{Rem_weak_const_a_szigma2}) yielding that
 \ $\PP(\hbpsi_T \in \Psi)=\PP(H^{-1}(\bG_T^{-1} \Bf_T) \in \Psi) \to 1$
 \ as \ $T \to \infty$ \ for each \ $\bpsi \in \Psi$, \ and hence
 \[
   \PP\Bigl(H^{-1}(\bG_T^{-1} \Bf_T)
             =\argmax_{\bpsi\in\Psi}
               \Lambda_T(\bpsi, \btpsi)\Bigl)
   \to 1 \qquad \text{as \ $T \to \infty$.}
 \]
Consequently, the probability that there exists a unique MLE \ $\bpsi_T^*$ \ of
 \ $\bpsi$ \ on the set \ $\Psi$ \ based on the observations
 \ $(Y_t, S_t)_{t\in[0,T]}$ \ converges to \ $1$ \ as \ $T \to \infty$, \ and
 \ $\PP(\bpsi_T^* = \hbpsi_T) \to 1$ \ as \ $T \to \infty$.
\proofend
\end{Rem}

\begin{Rem}
To make it clear, we point out that the expression for \ $(\htheta_T,\hkappa_T$) \ in
 \eqref{MLE_coordinatewise} is not the same as the MLE of \ $(\theta,\kappa)$ \ based only
 on the continuous time observation \ $(Y_t)_{t\in[0,T]}$ \ for the first coordinate process
 of \eqref{j_Heston_SDE}, see, e.g., Overbeck \cite{Ove}, because our statistical setup is different.
\proofend
\end{Rem}

\begin{Rem}\label{Luschgy}
In the literature there is another way of deriving an MLE.
S{\o}rensen \cite{SorM} defined an MLE of \ $\bpsi$ \ as a solution of
 the equation \ $\dot{\Lambda}_T(\bpsi) = 0$, \ where \ $\dot{\Lambda}_T(\bpsi)$ \ is the so-called  score vector
 given in formula (3.3) in S{\o}rensen \cite{SorM}.
Luschgy \cite{Lus2}, \cite{Lus} called this equation as an estimating equation.
With the notations of the proof of Lemma \ref{RN}, taking into account of the form of
 \ $\beta^{(\btpsi,\bpsi)}$ \ and the fact that \ $V^{(\bpsi,\btpsi)}$ \ is
 identically one, we have
 \begin{align*}
  \dot{\Lambda}_T(\bpsi)
  &:= \int_0^T
       \dot{H}(\bpsi)^\top
       A(Y_u, S_u)^\top \Sigma(Y_u, S_u)^{-1}
       \begin{bmatrix}
        \dd Y^\cont_u \\
        \dd S^\cont_u
       \end{bmatrix} \\
  &= \int_0^T
      \dot{H}(\bpsi)^\top
      A(Y_u, S_u)^\top \Sigma(Y_u, S_u)^{-1}
      \left(\begin{bmatrix}
             \dd Y_u \\
             \dd S_u
            \end{bmatrix}
            - A(Y_u, S_u) H(\bpsi) \, \dd u
            - \begin{bmatrix}
               0 \\
               S_{u-} \, \dd L_u
              \end{bmatrix}\right)
 \end{align*}
 for \ $\bpsi \in \RR^3$ \ and \ $T \in (0, \infty)$ \ with
 \[
   \dot{H}(\bpsi)
   = \begin{bmatrix}
      \kappa & \theta & 0 \\
      0 & 1 & 0 \\
      0 & 0 & 1
     \end{bmatrix} ,
   \qquad
   \bpsi = \begin{bmatrix}
            \theta \\ \kappa \\ \mu
           \end{bmatrix}
   \in \RR^3 .
 \]
Using \eqref{clmp} and the definitions of \ $\Bf_T$ \ and \ $\bG_T$, \ we obtain
 \[
   \dot{\Lambda}_T(\bpsi)
   = {\dot{H}(\bpsi)}^\top (\Bf_T - \bG_T H(\bpsi)) ,
 \]
 hence the estimating equation \ $\dot{\Lambda}_T(\bpsi) = 0$, \ $\bpsi \in \RR^3$,
 \ has a unique solution
 \ $H^{-1}(\bG_T^{-1} \Bf_T)$ \ provided that \ $\bG_T$ \ is strictly positive
 definite and \ $(\bG_T^{-1} \Bf_T)_2 \ne 0$, \ which hold almost surely.
Recall that this unique solution coincides with \ $\hbpsi_T$, \ see \eqref{MLE}.
\proofend
\end{Rem}

\section{Consistency of MLE}
\label{section_CMLE}

\begin{Thm}\label{Thm_MLE_cons}
If \ $\theta, \kappa \in (0, \infty)$ \ with
 \ $\theta \kappa \in \bigl( \frac{\sigma^2}{2}, \infty \bigr)$, \ $\mu \in \RR$,
 \ $\sigma \in (0, \infty)$, \ $\varrho \in (-1, 1)$, \ and
 \ $(Y_0, S_0) = (y_0, s_0) \in (0, \infty)^2$, \ then the MLE of
 \ $\bpsi = (\theta, \kappa, \mu)$ \ is strongly consistent, i.e.,
 \ $\hbpsi_T = \bigl(\htheta_T, \hkappa_T, \hmu_T\bigr)
    \as \bpsi = (\theta, \kappa, \mu)$
 \ as \ $T \to \infty$.
\end{Thm}

\noindent{\bf Proof.}
Obviously, it is enough to show that \ $\bG_T^{-1} \Bf_T \as H(\bpsi)$ \ as
 \ $T \to \infty$, \ since then
 \ $\hbpsi_T = H^{-1}(\bG_T^{-1} \Bf_T) \as H^{-1}(H(\bpsi)) = \bpsi$ \ as
 \ $T \to \infty$, \ using the continuity of \ $H^{-1}$ \ and that
 \ $\PP((\bG_T^{-1} \Bf_T)_2 = 0) = 0$, \ see Section \ref{section_EUMLE}.
Using the SDE \eqref{j_Heston_SDE_matrix}, we obtain
 \ $\Bf_T = \bG_T H(\bpsi) + \bh_T$, \ $T\in[0,\infty)$, \ with
 \begin{align*}
  \bh_T
  :=\int_0^T
     A(Y_u, S_u)^\top \Sigma(Y_u, S_u)^{-1} \Gamma(Y_u, S_u)
     \begin{bmatrix} \dd W_u \\ \dd B_u \end{bmatrix} .
 \end{align*}
Thus
 \begin{equation}\label{hT}
  \bG_T^{-1} \Bf_T - H(\bpsi)
  = \bG_T^{-1} \bh_T
  = (\bD_T \bG_T)^{-1} (\bD_T \bh_T) , \qquad T \in (0, \infty) ,
 \end{equation}
 with
 \[
   \bD_T
   := \begin{bmatrix}
       \frac{1}{\int_0^T \frac{\dd u}{Y_u}} & 0 & 0 \\
       0 & \frac{1}{\int_0^T Y_u \, \dd u} & 0 \\
       0 & 0 & \frac{1}{\int_0^T \frac{\dd u}{Y_u}}
      \end{bmatrix} ,
 \]
 provided that \ $\bG_T$ \ is invertible and \ $\int_0^T \frac{\dd u}{Y_u} > 0$
 \ and \ $\int_0^T Y_u \, \dd u > 0$, \ which hold almost surely, see Section
 \ref{section_EUMLE}.
We have
 \[
   \bD_T \bG_T
   = \frac{1}{(1-\varrho^2)\sigma^2}
     \begin{bmatrix}
      1 & -\frac{1}{T^{-1}\int_0^T \frac{\dd u}{Y_u}} & -\varrho\sigma \\
      -\frac{1}{T^{-1}\int_0^T Y_u \, \dd u} & 1
       & \frac{\varrho\sigma}{T^{-1}\int_0^T Y_u \, \dd u} \\
      -\varrho\sigma
       & \frac{\varrho\sigma}{T^{-1}\int_0^T \frac{\dd u}{Y_u}} & \sigma^2
     \end{bmatrix} , \qquad T \in (0, \infty) .
 \]
Using that
 \ $\Sigma(Y_u, S_u)^{-1} = (\Gamma(Y_u, S_u)^\top)^{-1} \Gamma(Y_u, S_u)^{-1}$,
 \ $u \in [0, T]$,
 \ and \eqref{bB^{-1}bA}, we obtain
 \begin{equation}\label{hTm}
  \begin{aligned}
  \bh_T
  &= \int_0^T
      A(Y_u, S_u)^\top (\Gamma(Y_u, S_u)^\top)^{-1}
      \begin{bmatrix} \dd W_u \\ \dd B_u \end{bmatrix}
   = \int_0^T
      (\Gamma(Y_u, S_u)^{-1} A(Y_u, S_u))^\top
      \begin{bmatrix} \dd W_u \\ \dd B_u \end{bmatrix} \\
  &= \frac{1}{\sigma\sqrt{1-\varrho^2}}
      \int_0^T
      \frac{1}{\sqrt{Y_u}}
      \begin{bmatrix}
       \sqrt{1-\varrho^2} \, \dd W_u - \varrho \, \dd B_u \\
       - \sqrt{1-\varrho^2} Y_u \, \dd W_u + \varrho Y_u \, \dd B_u \\
       \sigma \, \dd B_u
      \end{bmatrix} ,
  \end{aligned}
 \end{equation}
 and hence
 \[
   \bD_T \bh_T
   = \frac{1}{\sigma\sqrt{1-\varrho^2}}
     \begin{bmatrix} \oh_T^{(1)} \\ \oh_T^{(2)} \\ \oh_T^{(3)} \end{bmatrix}
 \]
 with
 \begin{gather*}
  \oh_T^{(1)}
  :=\frac{\int_0^T
           \frac{\sqrt{1-\varrho^2}\,\dd W_u - \varrho\,\dd B_u}{\sqrt{Y_u}}}
         {\int_0^T \frac{\dd u}{Y_u}} , \quad
  \oh_T^{(2)}
  :=-\frac{\int_0^T
           \sqrt{Y_u}(\sqrt{1-\varrho^2}\,\dd W_u - \varrho\,\dd B_u)}
         {\int_0^T Y_u \, \dd u} , \quad
  \oh_T^{(3)}
  :=\sigma
    \frac{\int_0^T \frac{\dd B_u}{\sqrt{Y_u}}}
         {\int_0^T \frac{\dd u}{Y_u}} .
 \end{gather*}
By part (i) of Theorem \ref{Ergodicity}, \ $\EE(Y_\infty) = \theta$ \ and
 \ $\EE\bigl(\frac{1}{Y_\infty}\bigr) = {\frac{2\kappa}{2\theta\kappa-\sigma^2}}$,
 \ and hence, part (ii) of Theorem \ref{Ergodicity} implies
 \begin{align}\label{bDbG}
  \bD_T \bG_T
  \as
  \frac{1}{(1-\varrho^2)\sigma^2}
  \begin{bmatrix}
   1 & -\frac{1}{\EE\bigl(\frac{1}{Y_\infty}\bigr)} & -\varrho\sigma \\
   -\frac{1}{\EE(Y_\infty)} & 1 & \frac{\varrho\sigma}{\EE(Y_\infty)} \\
   -\varrho\sigma
    & \frac{\varrho\sigma}{\EE\bigl(\frac{1}{Y_\infty}\bigr)} & \sigma^2
  \end{bmatrix}
  =: \frac{1}{(1-\varrho^2)\sigma^2} \bS
  \qquad \text{as \ $T \to \infty$.}
 \end{align}
We have
 \[
   \det(\bS)
   = (1-\varrho^2) \sigma^2
     \biggl(1 - \frac{1}{\EE(Y_\infty)\EE\bigl(\frac{1}{Y_\infty}\bigr)}\biggr)
   > 0
 \]
 since
 \ $\EE(Y_\infty) \EE\bigl(\frac{1}{Y_\infty}\bigr)
    = \frac{2\theta\kappa}{2\theta\kappa-\sigma^2} > 1$,
 \ hence the matrix \ $\bS$ \ is invertible, and we conclude
 \begin{equation}\label{bDbG^{-1}}
  (\bD_T \bG_T)^{-1}
  \as (1-\varrho^2) \sigma^2 \bS^{-1}
  \qquad \text{as \ $T \to \infty$.}
 \end{equation}
The aim of the following discussion is to show convergence
 \begin{equation}\label{bh}
  \bD_T \bh_T \as \bzero \qquad \text{as \ $T \to \infty$.}
 \end{equation}
The strong law of large numbers for continuous local martingales (see, e.g.,
 Theorem \ref{DDS_stoch_int}) implies \ $\oh_T^{(1)} \as 0$ \ as
 \ $T \to \infty$, \ since, by part (ii) of Theorem \ref{Ergodicity},
 \[
   \frac{1}{T} \int_0^T \frac{\dd s}{Y_s}
   \as \EE\Bigl(\frac{1}{Y_\infty}\Bigr)
       = \frac{2\kappa}{2\theta\kappa-\sigma^2} \in (0, \infty)
   \qquad \text{as \ $T \to \infty$,}
 \]
 implying
 \[
   \int_0^T \frac{\dd s}{Y_s}
   = T \cdot \frac{1}{T} \int_0^T \frac{\dd s}{Y_s}
   \as \infty
   \qquad \text{as \ $T \to \infty$.}
 \]
Convergences \ $\oh_T^{(2)} \as 0$ \ as \ $T \to \infty$, \ and
 \ $\oh_T^{(3)} \as 0$ \ as \ $T \to \infty$ \ can be proved in the
 same way, since, by part (ii) of Theorem \ref{Ergodicity},
 \[
   \frac{1}{T} \int_0^T Y_s \, \dd s \as \EE(Y_\infty) = \theta \in (0, \infty)
   \qquad \text{as \ $T \to \infty$,}
 \]
 implying
 \[
   \int_0^T Y_s \, \dd s
   = T \cdot \frac{1}{T} \int_0^T Y_s \, \dd s
   \as \infty
   \qquad \text{as \ $T \to \infty$.}
 \]
Consequently, we conclude \eqref{bh}.
By \eqref{hT}, \eqref{bDbG^{-1}} and \eqref{bh}, we obtain
 \ $\bG_T^{-1} \Bf_T \as H(\bpsi)$ \ as \ $T \to \infty$, \ hence we conclude
 the statement.
\proofend

\begin{Rem}\label{Rem_weak_const_a_szigma2}
For the case \ $\theta \kappa = \frac{\sigma^2}{2}$, \ Theorem \ref{Thm_MLE=}
 implies weak consistency of the MLE of \ $(\theta, \kappa, \mu)$.
\proofend
\end{Rem}

\section{Asymptotic behaviour of MLE}
\label{section_AMLE}

\begin{Thm}\label{Thm_MLE}
If \ $\theta, \kappa \in (0, \infty)$ \ with
 \ $\theta \kappa \in \bigl( \frac{\sigma^2}{2}, \infty \bigr)$, \ $\mu \in \RR$,
 \ $\sigma \in (0, \infty)$, \ $\varrho \in (-1, 1)$, \ and
 \ $(Y_0, S_0) = (y_0, s_0) \in (0, \infty)^2$, \ then the MLE of
 \ $\bpsi = (\theta, \kappa, \mu)$ \ is asymptotically normal, namely,
 \begin{align}\label{MLE_sub}
  T^{1/2} (\hbpsi_T - \bpsi)
  \distr \cN_3(\bzero, \bV) \qquad
  \text{as \ $T \to \infty$,}
 \end{align}
 where the matrix \ $\bV$ \ is given by
 \begin{align}\label{exp_G_infty_inverse_explicit}
  \frac{1}{2\kappa^3}
     \begin{bmatrix}
      \sigma^2(2 \theta \kappa - \varrho^2 \sigma^2)
       & - 2 (1-\varrho^2) \sigma^2 \kappa^2
       & \varrho \sigma \kappa (2 \theta \kappa - \sigma^2) \\
      - 2 (1-\varrho^2) \sigma^2 \kappa^2
       & 4 \kappa^4 (1-\varrho^2) & 0 \\
      \varrho \sigma \kappa (2 \theta \kappa - \sigma^2) & 0 &
      \kappa^2 (2 \theta \kappa - \sigma^2)
     \end{bmatrix} .
 \end{align}
With a random scaling, we have
 \begin{align}\label{MLE_subr}
  \bR_T \bQ_T (\hbpsi_T - \bpsi)
  \distr \cN_3(\bzero, \bI_3) \qquad
  \text{as \ $T \to \infty$,}
 \end{align}
 where \ $\bI_3$ \ denotes the \ $3 \times 3$ \ identity matrix, and \ $\bR_T$,
 \ $T \in (0, \infty)$, \ and \ $\bQ_T$, \ $T \in (0, \infty)$, \ are
 \ $3 \times 3$ \ (not uniquely determined) random matrices with properties
 \ $T^{-1/2} \bR_T \stoch \bC$ \ as \ $T \to \infty$ \ with some
 \ $\bC \in \RR^{3\times 3}$, \ $\bR_T^\top \bR_T = \bG_T$, \ $T \in (0, \infty)$,
 \ and \ $\bQ_T \stoch \bQ$ \ as \ $T \to \infty$, \ where
 \[
   \bQ
   := \begin{bmatrix}
       \kappa & \theta & 0 \\
       0 & 1 & 0 \\
       0 & 0 & 1
      \end{bmatrix} .
 \]
For a possible choice of \ $\bR_T$ \ and \ $\bQ_T$, $T \in (0, \infty)$, \ see
 Remark \ref{Rem_Q_R}.
\end{Thm}

\begin{Rem}\label{Rem_Q_R}
Note that the limiting covariance matrix \ $\bV$ \ in \eqref{MLE_sub} depends
 only on the unknown parameters \ $\theta$ \ and \ $\kappa$, \ but not on (the unknown)
 \ $\mu$.
\ The advantage of the random scaling is that the limiting covariance matrix in \eqref{MLE_subr}
 is the \ $3\times 3$ \ identity matrix \ $\bI_3$ \ which does not depend on any of the
 unknown parameters.
Note also that for \ $\bR_T$ \ and \ $\bQ_T$ \ one can choose, for instance,
 \[
   \bR_T
   = \frac{1}{\sigma\sqrt{1-\varrho^2}\sqrt{\int_0^T\frac{\dd u}{Y_u}}}
     \begin{bmatrix}
      \int_0^T\frac{\dd u}{Y_u} & -T & -\varrho\sigma\int_0^T\frac{\dd u}{Y_u} \\
      0 & \sqrt{\int_0^T Y_u\,\dd u \int_0^T\frac{\dd u}{Y_u}-T^2} & 0 \\
      0 & 0 & \sigma\sqrt{1-\varrho^2}\int_0^T\frac{\dd u}{Y_u}
     \end{bmatrix}
 \]
 and
 \[
   \bQ_T
   = \begin{bmatrix}
      \frac{\frac{\sigma^2}{2T}\int_0^T \frac{\dd u}{Y_u}}
           {\bigl(\frac{1}{T}\int_0^T Y_u\,\dd u\bigr)
            \bigl(\frac{1}{T}\int_0^T\frac{\dd u}{Y_u}\bigr)-1}
       & \frac{1}{T} \int_0^T Y_u \, \dd u & 0 \\
      0 & 1 & 0 \\
      0 & 0 & 1
     \end{bmatrix} .
 \]
Indeed, we have $\bR_T^\top \bR_T = \bG_T$, $T \in (0, \infty)$ (which is,
 in fact, the Cholesky factorization of $\bG_T$),
 \begin{align*}
  T^{-1/2} \bR_T
  &= \frac{1}{\sigma\sqrt{1-\varrho^2}\sqrt{\frac{1}{T}\int_0^T\frac{\dd u}{Y_u}}}
     \begin{bmatrix}
      \frac{1}{T}\int_0^T\frac{\dd u}{Y_u} & -1 &
       -\frac{\varrho\sigma}{T}\int_0^T\frac{\dd u}{Y_u} \\[2mm]
      0 & \sqrt{\bigl(\frac{1}{T}\int_0^T Y_u\,\dd u\bigr)
                \bigl(\frac{1}{T}\int_0^T\frac{\dd u}{Y_u}\bigr)-1} & 0 \\
      0 & 0 & \frac{\sigma\sqrt{1-\varrho^2}}{T}\int_0^T\frac{\dd u}{Y_u}
     \end{bmatrix} \\
  &\as \frac{1}{\sigma\sqrt{1-\varrho^2}{\sqrt{ \EE\bigl(\frac{1}{Y_\infty}\bigr)}}}
     \begin{bmatrix}
      \EE\bigl(\frac{1}{Y_\infty}\bigr) & -1 &
       -\varrho\sigma\EE\bigl(\frac{1}{Y_\infty}\bigr) \\[2mm]
      0 & \sqrt{\EE(Y_\infty)\EE\bigl(\frac{1}{Y_\infty}\bigr)-1} & 0 \\
      0 & 0 & \sigma\sqrt{1-\varrho^2}\EE\bigl(\frac{1}{Y_\infty}\bigr)
     \end{bmatrix}
 \end{align*}
 as \ $T \to \infty$, \ and \ $\bQ_T \as \bQ$ \ as \ $T \to \infty$, \ since
 \[
   \frac{\frac{\sigma^2}{2T}\int_0^T \frac{\dd u}{Y_u}}
        {\bigl(\frac{1}{T}\int_0^T Y_u\,\dd u\bigr)
         \bigl(\frac{1}{T}\int_0^T\frac{\dd u}{Y_u}\bigr)-1}
   \as \frac{\frac{\sigma^2}{2}\EE\bigl(\frac{1}{Y_\infty}\bigr)}
            {\EE(Y_\infty)\EE\bigl(\frac{1}{Y_\infty}\bigr)-1}
   = \kappa, \qquad
   \frac{1}{T} \int_0^T Y_s \, \dd s \as \EE(Y_\infty) = \theta
 \]
 as \ $T \to \infty$ \ by part (i) of Theorem \ref{Ergodicity}.
Hence then the random scaling factor has the form
 \[
   \bR_T \bQ_T
   = \frac{1}{\sigma\sqrt{1-\varrho^2}\sqrt{\int_0^T\frac{\dd u}{Y_u}}}
     \begin{bmatrix}
      \frac{\frac{\sigma^2T}{2}\bigl(\int_0^T\frac{\dd u}{Y_u}\bigr)^2}
           {\int_0^T Y_u\,\dd u \int_0^T\frac{\dd u}{Y_u}-T^2}
       & \frac{1}{T}\bigl(\int_0^T Y_u\,\dd u \int_0^T\frac{\dd u}{Y_u}-T^2\bigr)
       & -\varrho\sigma\int_0^T\frac{\dd u}{Y_u} \\[1mm]
      0 & \sqrt{\int_0^T Y_u\,\dd u \int_0^T\frac{\dd u}{Y_u}-T^2} & 0 \\
      0 & 0 & \sigma\sqrt{1-\varrho^2}\int_0^T\frac{\dd u}{Y_u}
     \end{bmatrix} .
 \]
\proofend
\end{Rem}

\noindent{\bf Proof of Theorem \ref{Thm_MLE}.}
For \eqref{MLE_sub}, it is enough to prove
 \begin{equation}\label{MLE_sub0}
  T^{1/2} (\bG_T^{-1} \Bf_T - H(\bpsi))
  \distr \cN_3(\bzero, \bV_0) \qquad
  \text{as \ $T \to \infty$,}
 \end{equation}
 where
 \[
   \bV_0
   := \frac{2\theta\kappa-\sigma^2}{2\kappa}
     \begin{bmatrix}
      \sigma^2+(1-\varrho^2)(2\theta\kappa-\sigma^2) & 2\kappa(1-\varrho^2)
       & \varrho\sigma \\
      2\kappa(1-\varrho^2) & \frac{4\kappa^2(1-\varrho^2)}{2\theta\kappa-\sigma^2} & 0 \\
      \varrho\sigma & 0 & 1
     \end{bmatrix} .
 \]
Indeed, then one can apply Lemma \ref{Lem_Kallenberg} with \ $\cS_1 = \cS_2 = \RR^3$,
 \ $\cC = \RR^3$, \ with a random vector \ $\bxi$ \ having distribution
 \ $\cN_3(\bzero, \bV_0)$, \ with
 \ $\xi_T = T^{1/2} (\bG_T^{-1} \Bf_T - H(\bpsi))$, \ $T \in (0, \infty)$, \ and
 with functions \ $F : \RR^3 \to \RR^3$ \ and
 \ $F_T : \RR^3 \to \RR^3$, \ $T \in (0, \infty)$, \ given by
 \[
   F(\bx)
   := \bQ^{-1} \bx , \qquad
   F_T(\bx)
   := \begin{cases}
       \begin{bmatrix}
        \frac{x_1-\theta x_2}{T^{-1/2}x_2+\kappa} \\
        x_2 \\
        x_3
       \end{bmatrix}
       & \text{if \ $x_2 \ne -T^{1/2} \kappa$,} \\
       \bzero & \text{if \ $x_2 = -T^{1/2} \kappa$,}
      \end{cases}
 \]
 for \ $\bx = (x_1, x_2, x_3) \in \RR^3$ \ and \ $T \in (0, \infty)$.
\ We have
 \begin{align*}
   F_T(T^{1/2}(\bG_T^{-1} \Bf_T - H(\bpsi)))
    = \begin{bmatrix}
        T^{1/2}\frac{(\bG_T^{-1} \Bf_T)_1 - \theta\kappa - \theta((\bG_T^{-1} \Bf_T)_2 - \kappa) }{(\bG_T^{-1} \Bf_T)_2} \\
       T^{1/2} ( (\bG_T^{-1} \Bf_T)_2 - \kappa ) \\
        T^{1/2} ( (\bG_T^{-1} \Bf_T)_3 - \mu ) \\
     \end{bmatrix}
     = T^{1/2}(\hbpsi_T - \bpsi),
 \end{align*}
  provided that \ $(\bG_T^{-1} \Bf_T)_2 \ne 0$, \ which holds almost surely.
Moreover, \ $F_T(\bx_T) \to F(\bx)$ \ as \ $T \to \infty$ \ if \ $\bx_T \to \bx$
 \ as \ $T \to \infty$, \ since then, for sufficiently large \ $T \in (0, \infty)$,
 \ we have \ $(\bx_T)_2 \ne -T^{1/2} \kappa$.
\ Consequently, \eqref{MLE_sub0} and Lemma \ref{Lem_Kallenberg} imply
 \[
   T^{1/2}(\hbpsi_T - \bpsi)
   = F_T(T^{1/2}(\bG_T^{-1} \Bf_T - H(\bpsi)))
   \distr F(\bxi)
   = \bQ^{-1} \bxi
   \distre \cN_3(\bzero, \bQ^{-1} \bV_0 (\bQ^{-1})^\top)
 \]
 as \ $T \to \infty$, \ where \ $\bQ^{-1} \bV_0 (\bQ^{-1})^\top = \bV$, \ hence
 we obtain \eqref{MLE_sub}.

By the first equality in \eqref{hT}, we have
 \begin{equation}\label{ThT}
  T^{1/2} (\bG_T^{-1} \Bf_T - H(\bpsi))
  = T^{1/2} \bG_T^{-1} \bh_T
  = (T^{-1} \bG_T)^{-1} (T^{-1/2} \bh_T) ,
 \end{equation}
 provided that \ $\bG_T$ \ is invertible, which holds almost surely, see Section
 \ref{section_EUMLE}.
By part (i) of Theorem \ref{Ergodicity},
 \ $\EE(Y_\infty) = \theta \in (0, \infty)$
 \ and
 \ $\EE\bigl(\frac{1}{Y_\infty}\bigr)
    = \frac{2\kappa}{2\theta\kappa-\sigma^2} \in (0, \infty)$,
 \ and hence, part (ii) of Theorem \ref{Ergodicity} and \eqref{bGt} imply
 \begin{align}\label{bG}
  T^{-1} \bG_T \as \EE(\bG_\infty) \qquad \text{as \ $T \to \infty$,}
 \end{align}
 with
 \[
   \bG_\infty
   := \frac{1}{(1-\varrho^2)\sigma^2 Y_\infty}
      \begin{bmatrix}
       1 & -Y_\infty & -\varrho \sigma \\
       -Y_\infty & Y_\infty^2 & \varrho \sigma Y_\infty \\
       -\varrho \sigma & \varrho \sigma Y_\infty & \sigma^2
      \end{bmatrix} ,
 \]
 where \ $Y_\infty$ \ has Gamma distribution with parameters
 \ $2\theta \kappa / \sigma^2$ \ and \ $2\kappa / \sigma^2$.
\ The matrix \ $\EE(\bG_\infty)$ \ is invertible, namely,
 \begin{align*}
  [\EE(\bG_\infty)]^{-1}
  &=\frac{1}
         {\bigl(\EE(Y_\infty)\EE\bigl(\frac{1}{Y_\infty}\bigr)-1\bigr)
          \EE\bigl(\frac{1}{Y_\infty}\bigr)} \\
  &\quad
    \times
    \begin{bmatrix}
     \sigma^2 \EE(Y_\infty) \EE\bigl(\frac{1}{Y_\infty}\bigr)
     - \varrho^2 \sigma^2
      & (1 - \varrho^2) \sigma^2 \EE\bigl(\frac{1}{Y_\infty}\bigr)
      & \varrho \sigma \EE(Y_\infty) \EE\bigl(\frac{1}{Y_\infty}\bigr)
        - \varrho \sigma \\
     (1 - \varrho^2) \sigma^2 \EE\bigl(\frac{1}{Y_\infty}\bigr)
      & (1 - \varrho^2) \sigma^2 \bigl(\EE\bigl(\frac{1}{Y_\infty}\bigr)\bigr)^2
      & 0 \\
     \varrho \sigma \EE(Y_\infty) \EE\bigl(\frac{1}{Y_\infty}\bigr)
     - \varrho \sigma
      & 0 & \EE(Y_\infty) \EE\bigl(\frac{1}{Y_\infty}\bigr) - 1
    \end{bmatrix} ,
 \end{align*}
 since
 \ $\EE(Y_\infty) \EE\bigl(\frac{1}{Y_\infty}\bigr)
    = \frac{2\theta\kappa}{2\theta\kappa-\sigma^2} > 1$,
 \ which yields \ $[\EE(\bG_\infty)]^{-1} = \bV_0$.
\ Whence we conclude
 \begin{equation}\label{bG^{-1}}
  (T^{-1} \bG_T)^{-1} \as \bV_0
  \qquad \text{as \ $T \to \infty$.}
 \end{equation}
By \eqref{hTm}, the process \ $(\bh_t)_{t\in[0,\infty)}$ \ is a 3-dimensional
 continuous local martingale with (predictable) quadratic variation process
 \ $\langle \bh \rangle_t = \bG_t$, \ $t \in [0,\infty)$.
\ Using \eqref{bG}, the central limit theorem for multidimensional continuous local
 martingales, see Theorem \ref{MCLT}, yields
 \ $T^{-1/2} \bh_T \distr \cN_3(\bzero, \EE(\bG_\infty))
    = \cN_3(\bzero, \bV_0^{-1})$
 \ as \ $T \to \infty$.
\ Hence, by \eqref{ThT} and \eqref{bG^{-1}},
 \[
   T^{1/2} (\bG_T^{-1} \Bf_T - H(\bpsi))
   \distr
   \cN_3(\bzero, \bV_0 \bV_0^{-1} \bV_0)
   = \cN_3(\bzero, \bV_0)
   \qquad \text{as \ $T \to \infty$,}
 \]
 thus we obtain \eqref{MLE_sub0}.

With random scaling, by \eqref{MLE_sub} and Slutsky's lemma, we obtain
 \[
   \bR_T \bQ_T (\hbpsi_T - \bpsi)
   = (T^{-1/2} \bR_T) \bQ_T \bigl[T^{1/2} (\hbpsi_T - \bpsi)\bigr]
   \distr \cN_3(\bzero, (\bC \bQ) \bV (\bC \bQ)^\top)
 \]
 as \ $T \to \infty$.
\ Moreover, by the assumptions on \ $\bR_T$, \ $T \in (0, \infty)$,
 \[
   T^{-1} \bG_T
   = (T^{-1/2} \bR_T)^\top (T^{-1/2} \bR_T)
   \stoch \bC^\top \bC \qquad
   \text{as \ $T \to \infty$.}
 \]
Thus, comparing with \eqref{bG}, we obtain
 \ $\bC^\top \bC = \EE(\bG_\infty) = \bV_0^{-1}$.
\ Using \ $\bQ^{-1} \bV_0 (\bQ^{-1})^\top = \bV$, \ we obtain
 \[
   (\bC \bQ) \bV (\bC \bQ)^\top
   = (\bC \bQ) \bQ^{-1} (\bC^\top \bC)^{-1} (\bQ^{-1})^\top (\bC \bQ)^\top
   = \bI_3 ,
 \]
 and we conclude \eqref{MLE_subr}.
\proofend

\begin{Thm}\label{Thm_MLE=}
If \ $\theta, \kappa \in (0, \infty)$ \ with
 \ $\theta \kappa = \frac{\sigma^2}{2}$, \ $\mu \in \RR$,
 \ $\sigma \in (0, \infty)$, \ $\varrho \in (-1, 1)$, \ and
 \ $(Y_0, S_0) = (y_0, s_0) \in (0, \infty)^2$, \ then
 \begin{align}\label{abm=}
  \begin{bmatrix}
   T^{1/2} (\htheta_T - \theta) \\ T^{1/2} (\hkappa_T - \kappa) \\ T (\hmu_T - \mu)
  \end{bmatrix}
  \distr
   \begin{bmatrix}
   -\frac{\sigma^2 \sqrt{1 - \varrho^2} }{\sqrt{2\kappa^3}} \, Z_1 \\[1mm]
   \sqrt{2(1 - \varrho^2)\kappa} \, Z_1 \\[1mm]
   \frac{\varrho\sigma}{\kappa\tau}
   + \frac{\sigma\sqrt{1-\varrho^2}}{\kappa\sqrt{\tau}} Z_2
  \end{bmatrix} \qquad
  \text{as \ $T \to \infty$,}
 \end{align}
 where \ $\tau := \inf\{ t \in [0, \infty) : \cW_t = 1\}$ \ with a standard
 Wiener process \ $(\cW_t)_{t\in[0,\infty)}$, \ and \ $Z_1$ \ and \ $Z_2$ \ are
 independent standard normally distributed random variables, independent from \ $\tau$.
\ With a random scaling, we have
 \begin{align}\label{abm=r}
  \begin{bmatrix}
   \frac{\sigma  T^2 }{2\sqrt{1-\varrho^2}\bigl(\int_0^T Y_u\,\dd u\bigr)^{3/2}}
   \, (\htheta_T - \theta) \\[4mm]
   \frac{1}{\sigma\sqrt{1-\varrho^2}}
   \bigl(\int_0^T Y_u \, \dd u\bigr)^{1/2}
   \, (\hkappa_T - \kappa) \\
   \frac{\sigma T^2}{2\int_0^T Y_u \, \dd u} \, (\hmu_T - \mu)
  \end{bmatrix}
  \distr
  \begin{bmatrix}
   -Z_1 \\
   Z_1 \\
   \frac{\varrho}{\tau} + \frac{ \sqrt{1-\varrho^2}}{\sqrt{\tau}} Z_2
  \end{bmatrix}
  \qquad \text{as \ $T \to \infty$.}
 \end{align}
\end{Thm}

Note that the limit distribution in Theorem \ref{Thm_MLE=} (which can be
 considered as the asymptotic error of the estimator \ $(\htheta_T, \hkappa_T, \hmu_T)$)
 \ is a mixed normal distribution.
Moreover, the first and second coordinates of the limit distributions in
 \eqref{abm=} and \eqref{abm=r} are linearly dependent.
In spite of this fact, one can give asymptotic confidence sets for \ $(\theta, \kappa)$,
 \ namely, ellipses together with their interiors and with center \ $(\htheta_T,\hkappa_T)$.
\ Indeed, the sum of the squares of the first two coordinates of the left-hand side of \eqref{abm=r},
 which one can call a normalized squared error of \ $(\theta, \kappa)$, \ converges
 weakly to \ $2Z_1^2$, \ being a chi-squared distribution of degree 1 (multiplied by 2).
Surprisingly, the mixed normal limit distributions of the third coordinate in
 \eqref{abm=} and \eqref{abm=r} are not centered.
In Appendix \ref{App_density} we derive an explicit formula for the density function of
 \ $\frac{\varrho\sigma}{\kappa\tau}
    + \frac{\sigma\sqrt{1-\varrho^2}}{\kappa\sqrt{\tau}} Z_2$,
 \ which is the limit distribution of \ $T (\hmu_T - \mu)$ \ as \ $T \to \infty$ \ in Theorem
 \ref{Thm_MLE=}.

\noindent{\bf Proof of Theorem \ref{Thm_MLE=}.}
Since \ $\PP(\tau \in (0, \infty)) = 1$, \ the limit distributions in
 \eqref{abm=} and \eqref{abm=r} are well defined.
We have again \ $\EE(Y_\infty) = \theta \in (0, \infty)$, \ implying
 \begin{align}\label{SLLN=}
  \frac{1}{T} \int_0^T Y_u \, \dd u \as \theta \qquad \text{and} \qquad
  \int_0^T Y_u \, \dd u \as \infty \qquad
  \text{as \ $T \to \infty$.}
 \end{align}
Due to Ben Alaya and Kebaier \cite[Proposition 4]{BenKeb1}, we have
 \begin{align}\label{AK}
  \frac{1}{T^2} \int_0^T \frac{\dd u}{Y_u} \distr \tau^* \qquad
  \text{as \ $T \to \infty$,}
 \end{align}
 where \ $\tau^* := \inf\{ t \in [0, \infty) : \cW_t^* = \frac{\kappa}{\sigma}\}$
 \ with a standard Wiener process \ $(\cW_t^*)_{t\in[0,\infty)}$.
\ Applying the scaling property of a standard Wiener process, we obtain
 \begin{align*}
  \tau^*
  &= \inf\Bigl\{ t \in [0, \infty) : \frac{\sigma}{\kappa} \cW_t^* = 1\Bigr\}
   \distre \inf\Bigl\{ t \in [0, \infty) : \cW_{\frac{\sigma^2}{\kappa^2}\,t} = 1\Bigr\} \\
  &= \inf\Bigl\{\frac{\kappa^2}{\sigma^2}\,s  \in [0, \infty) : \cW_s = 1\Bigr\}
   = \frac{\kappa^2}{\sigma^2} \inf\Bigl\{s  \in [0, \infty) : \cW_s = 1\Bigr\}
   = \frac{\kappa^2}{\sigma^2} \tau ,
 \end{align*}
 where \ $\distre$ \ denotes equality in distribution.
We may and do suppose that \ $\tau^* = \frac{\kappa^2}{\sigma^2} \tau$.
\ Using \ $\PP(\tau^* \in (0, \infty)) = 1$, \ we conclude
 \[
   \frac{1}{\frac{1}{T}\int_0^T \frac{\dd u}{Y_u}}
    = \frac{1}{T} \, \frac{1}{\frac{1}{T^2} \int_0^T \frac{\dd u}{Y_u}}
    \distr 0 \cdot \frac{1}{\tau^*} = 0 \qquad
   \text{as \ $T \to \infty$,}
 \]
 and hence,
 \begin{equation}\label{stoch_1/Y}
  \frac{1}{\frac{1}{T}\int_0^T \frac{\dd u}{Y_u}} \stoch 0 \qquad
  \text{as \ $T \to \infty$,}
 \end{equation}
 implying also
 \[
   \frac{1}{\int_0^T \frac{\dd u}{Y_u}}
   = \frac{1}{T} \, \frac{1}{\frac{1}{T} \int_0^T \frac{\dd u}{Y_u}}
   \stoch 0 \qquad
   \text{as \ $T \to \infty$.}
 \]
Since the function
 \ $(0, \infty) \ni T \mapsto \frac{1}{\int_0^T \frac{\dd u}{Y_u}}$ \ is
 monotone decreasing, we obtain
 \[
   \frac{1}{\int_0^T \frac{\dd u}{Y_u}} \as 0 \qquad \text{and} \qquad
   \int_0^T \frac{\dd u}{Y_u} \as \infty \qquad
   \text{as \ $T \to \infty$.}
 \]
For \eqref{abm=}, it is enough to prove
 \begin{equation}\label{abm=0}
  \bC_T (\bG_T^{-1} \Bf_T - H(\bpsi))
  \distr
  \begin{bmatrix}
   \frac{\sigma^2}{\kappa\tau} \\[1mm]
   Z_1 \sqrt{2(1 - \varrho^2)\kappa} \\[1mm]
   \frac{\varrho\sigma}{\kappa\tau}
   + \frac{Z_2\sigma\sqrt{1-\varrho^2}}{\kappa\sqrt{\tau}}
  \end{bmatrix}
  \qquad \text{as \ $T \to \infty$,}
 \end{equation}
 where
 \[
   \bC_T := \begin{bmatrix}
             T & 0 & 0 \\ 0 & T^{1/2} & 0 \\ 0 & 0 & T
            \end{bmatrix} .
 \]
Indeed, then one can apply Lemma \ref{Lem_Kallenberg} with \ $\cS_1 = \cS_2 = \RR^3$,
 \ $\cC = \RR^3$, \ with
 \[
   \xi = \begin{bmatrix}
          \frac{\sigma^2}{\kappa\tau} \\[1mm]
          Z_1 \sqrt{2(1 - \varrho^2)\kappa} \\[1mm]
          \frac{\varrho\sigma}{\kappa\tau}
          + \frac{Z_2\sigma\sqrt{1-\varrho^2}}{\kappa\sqrt{\tau}}
         \end{bmatrix} ,
 \]
 with \ $\xi_T = \bC_T (\bG_T^{-1} \Bf_T - H(\bpsi))$, \ $T \in (0, \infty)$,
 \ and with functions \ $F : \RR^3 \to \RR^3$ \ and \ $F_T : \RR^3 \to \RR^3$,
 \ $T \in (0, \infty)$, \ given by
 \[
   F(\bx)
   := \bB \bx , \qquad
   F_T(\bx)
   := \begin{cases}
       \begin{bmatrix}
        \frac{T^{-1/2} x_1 - \theta x_2}{T^{-1/2}x_2+\kappa} \\
        x_2 \\
        x_3
       \end{bmatrix}
       & \text{if \ $x_2 \ne -T^{1/2} \kappa$,} \\
       \bzero & \text{if \ $x_2 = -T^{1/2} \kappa$,}
      \end{cases}
 \]
 for \ $\bx = (x_1, x_2, x_3) \in \RR^3$ \ and \ $T \in (0, \infty)$, \ where
 \[
   \bB := \begin{bmatrix}
           0 & -\frac{\sigma^2}{2\kappa^2} & 0 \\ 0 & 1 & 0 \\ 0 & 0 & 1
          \end{bmatrix} .
 \]
We have
 \begin{align*}
  F_T(\bC_T (\bG_T^{-1} \Bf_T - H(\bpsi)))
    = \begin{bmatrix}
        \frac{T^{-1/2}T((\bG_T^{-1} \Bf_T)_1 - \theta\kappa) - \theta T^{1/2} ((\bG_T^{-1} \Bf_T)_2 - \kappa) }
             {(\bG_T^{-1} \Bf_T)_2} \\[1mm]
        T^{1/2}( (\bG_T^{-1} \Bf_T)_2 - \kappa ) \\
        T( (\bG_T^{-1} \Bf_T)_3 - \mu ) \\
      \end{bmatrix}
    = \tbC_T (\hbpsi_T - \bpsi),
 \end{align*}
 provided that \ $(\bG_T^{-1} \Bf_T)_2 \ne 0$, \ which holds almost surely, where
 \[
   \tbC_T
        :=\begin{bmatrix}
           T^{1/2} & 0 & 0 \\ 0 & T^{1/2} & 0 \\ 0 & 0 & T
          \end{bmatrix} .
 \]
Moreover, \ $F_T(\bx_T) \to F(\bx)$ \ as \ $T \to \infty$ \ if \ $\bx_T \to \bx$
 \ as \ $T \to \infty$, \ since then, for sufficiently large \ $T \in (0, \infty)$,
 \ we have \ $(\bx_T)_2 \ne -T^{1/2} \kappa$.
\ Consequently, \eqref{abm=0} and Lemma \ref{Lem_Kallenberg} imply
 \[
   \tbC_T (\hbpsi_T - \bpsi)
   = F_T(\bC_T (\bG_T^{-1} \Bf_T - H(\bpsi)))
   \distr F(\bxi)
   = \bB \bxi \qquad \text{as \ $T \to \infty$,}
 \]
 hence we obtain \eqref{abm=}.

Now we turn to prove \eqref{abm=0}.
By the first equality in \eqref{hT}, we have
 \begin{equation}\label{deco}
  \begin{aligned}
   \bC_T (\bG_T^{-1} \Bf_T - H(\bpsi))
   &= \bC_T \bG_T^{-1} \bh_T
    = (\bC_T \bG_T^{-1} \bC_T) (\bC_T^{-1} \bh_T) \\
   &= (\bC_T^{-1} \bG_T \bC_T^{-1})^{-1} (\bC_T^{-1} \bh_T) ,
  \end{aligned}
 \end{equation}
 provided that \ $\bG_T$ \ is invertible, which holds almost surely, see Section
 \ref{section_EUMLE}.
We have
 \[
   \bC_T^{-1} \bG_T \bC_T^{-1}
   = \frac{1}{(1-\varrho^2)\sigma^2}
     \begin{bmatrix}
      T^{-2} \int_0^T \frac{\dd u}{Y_u} & -T^{-1/2}
       & -\varrho \sigma T^{-2} \int_0^T \frac{\dd u}{Y_u} \\
      -T^{-1/2} & T^{-1} \int_0^T Y_u \, \dd u & \varrho \sigma T^{-1/2} \\
      -\varrho \sigma T^{-2} \int_0^T \frac{\dd u}{Y_u} & \varrho \sigma T^{-1/2}
       & \sigma^2 T^{-2} \int_0^T \frac{\dd u}{Y_u}
     \end{bmatrix},
 \]
 and, by \eqref{hTm},
 \[
   \bC_T^{-1} \bh_T
   = \frac{1}{\sigma\sqrt{1-\varrho^2}}
     \begin{bmatrix} \tildeh_T^{(1)} \\ \tildeh_T^{(2)} \\ \tildeh_T^{(3)} \end{bmatrix}
 \]
 with
 \begin{gather*}
  \tildeh_T^{(1)}
  := \frac{1}{T}
     \int_0^T \frac{\sqrt{1-\varrho^2}\,\dd W_u - \varrho\,\dd B_u}{\sqrt{Y_u}} , \qquad
  \tildeh_T^{(2)}
  := - \frac{1}{\sqrt{T}} \int_0^T \sqrt{Y_u}(\sqrt{1-\varrho^2}\,\dd W_u - \varrho\,\dd B_u) , \\
  \tildeh_T^{(3)} := \frac{\sigma}{T} \int_0^T \frac{\dd B_u}{\sqrt{Y_u}} .
 \end{gather*}
By \eqref{SLLN=}, \ $\frac{1}{T} \int_0^T Y_u \, \dd u \as \theta$ \ as \ $T \to \infty$.
\ Ben Alaya and Kebaier \cite[proof of Theorem 7]{BenKeb2} proved \ $\frac{\log (Y_T)}{T} \stoch 0$
 \ as \ $T \to \infty$.
\ Using the SDE \eqref{j_Heston_SDE_matrix}, \eqref{help_ITO} and
 \ $\theta \kappa = \frac{\sigma^2}{2}$,
 \begin{equation}\label{intWY}
  \int_0^T \frac{\sigma \, \dd W_u}{\sqrt{Y_u}}
  = \int_0^T \frac{\dd Y_u}{Y_u} - \theta \kappa \int_0^T \frac{\dd u}{Y_u} + \kappa T
  = \log(Y_T) - \log(y_0) + \kappa T ,
 \end{equation}
 thus
 \begin{equation}\label{BKconv1}
  \frac{1}{T} \int_0^T \frac{\dd W_u}{\sqrt{Y_u}} \stoch \frac{\kappa}{\sigma} \qquad
  \text{as \ $T \to \infty$.}
 \end{equation}
Consequently, \eqref{abm=0} will follow from
 \begin{equation}\label{bGbh}
  (\bC_T^{-1} \bG_T \bC_T^{-1}, \bC_T^{-1} \bh_T) \distr (\tbG_\infty, \tbh_\infty) \qquad
  \text{as \ $T \to \infty$}
 \end{equation}
 with
 \[
   \tbG_\infty
   := \frac{1}{(1-\varrho^2)\sigma^2}
      \begin{bmatrix}
       \tau^* & 0 & -\varrho \sigma \tau^* \\
       0 & \theta & 0 \\
       -\varrho \sigma \tau^* & 0 & \sigma^2 \tau^*
      \end{bmatrix} , \quad
   \tbh_\infty
   := \frac{1}{\sigma\sqrt{1-\varrho^2}}
      \begin{bmatrix}
       \frac{\kappa \sqrt{1-\varrho^2}}{\sigma} - \varrho Z_2 \sqrt{\tau^*} \\
       -Z_3\sqrt{\theta (1-\varrho^2)} + \varrho Z_4 \sqrt{\theta} \\
       \sigma Z_2 \sqrt{\tau^*}
      \end{bmatrix} ,
 \]
 where \ $Z_3$ \ and \ $Z_4$ \ are independent standard normally distributed
 random variables, independent from \ $Z_2$ \ and \ $\tau^*$.
\ Indeed, provided that \eqref{bGbh} holds, by the continuous mapping theorem,
 \[
   (\bC_T \bG_T^{-1} \bC_T, \bC_T^{-1} \bh_T)
   \distr (\tbG_\infty^{-1}, \tbh_\infty) \qquad
   \text{as \ $T \to \infty$,}
 \]
 since \ $\tbG_\infty$ \ is invertible almost surely with inverse
 \[
   \tbG_\infty^{-1}
   = \frac{1}{\tau^*}
     \begin{bmatrix}
      \sigma^2 & 0 & \varrho \sigma \\
      0 & \frac{1}{\theta} (1 - \varrho^2) \sigma^2 \tau^* & 0 \\
      \varrho \sigma & 0 & 1
     \end{bmatrix} ,
 \]
 and hence, by \eqref{deco} and the continuous mapping theorem,
 \begin{align*}
  \bC_T (\bG_T^{-1} \Bf_T - H(\bpsi))
  \distr
  \tbG_\infty^{-1} \tbh_\infty \qquad \text{as \ $T \to \infty$,}
 \end{align*}
 where, with \ $Z_1 := - \sqrt{1-\varrho^2} Z_3 + \varrho Z_4$,
 \begin{align*}
  \tbG_\infty^{-1} \tbh_\infty
  &= \frac{1}{\sigma\tau^*\sqrt{1-\varrho^2}}
     \begin{bmatrix}
      \sigma^2 \Bigl(\frac{\kappa \sqrt{1-\varrho^2}}{\sigma} - \varrho Z_2 \sqrt{\tau^*}\Bigr)
       + \varrho \sigma^2 Z_2 \sqrt{\tau^*} \\
      \frac{1}{\sqrt \theta} \, (1 - \varrho^2) \sigma^2 \tau^*
       (-Z_3 \sqrt{1-\varrho^2} + \varrho Z_4) \\
      \varrho \sigma \bigl(\frac{\kappa \sqrt{1-\varrho^2}}{\sigma} - \varrho Z_2 \sqrt{\tau^*}\bigr)
      + \sigma Z_2 \sqrt{\tau^*}
     \end{bmatrix} \\
  &= \frac{1}{\sigma\tau^*\sqrt{1-\varrho^2}}
     \begin{bmatrix}
      \sigma \kappa \sqrt{1-\varrho^2} \\
      \frac{1}{\sqrt \theta} \, (1 - \varrho^2) \sigma^2 \tau^* Z_1 \\
      \varrho \kappa \sqrt{1-\varrho^2} + \sigma (1 - \varrho^2) Z_2 \sqrt{\tau^*}
     \end{bmatrix}
  = \begin{bmatrix}
     \frac{\kappa}{\tau^*} \\
     Z_1 \sqrt{2 (1 - \varrho^2) \kappa} \\
     \frac{\varrho \kappa} {\sigma\tau^*} + \frac{Z_2\sqrt{1-\varrho^2}}{\sqrt{\tau^*}}
    \end{bmatrix} ,
 \end{align*}
 thus we obtain \eqref{abm=0} using that \ $\tau^* = \frac{\kappa^2}{\sigma^2} \tau$.

Now we turn to prove \eqref{bGbh}.
It will follow from Slutsky's lemma, continuous mapping theorem and from
 \begin{align*}
  \begin{aligned}
   &\biggl( \frac{1}{\sqrt{T}} \int_0^T \sqrt{Y_u} \, \dd W_u,
            \frac{1}{T} \int_0^T \frac{\dd W_u}{\sqrt{Y_u}},
            \frac{1}{\sqrt{T}} \int_0^T \sqrt{Y_u} \, \dd B_u,
            \frac{1}{T} \int_0^T \frac{\dd B_u}{\sqrt{Y_u}},
            \frac{1}{T} \int_0^T Y_u \, \dd u,
            \frac{1}{T^2} \int_0^T \frac{\dd u}{Y_u} \biggr) \\
   &\distr
    \Bigl( \sqrt{\theta} \, Z_3, \frac{\kappa}{\sigma}, \sqrt{\theta} Z_4,
           Z_2 \sqrt{\tau^*}, \theta, \tau^* \Bigr) \qquad
    \text{as \ $T \to \infty$,}
  \end{aligned}
 \end{align*}
 which will be a consequence of \eqref{SLLN=}, \eqref{AK}, \eqref{BKconv1}, Slutsky's
 lemma (or part (v) of Theorem 2.7 in van der Vaart \cite{Vaart}), and
 \begin{align}\label{bGbh1}
  \begin{aligned}
   &\biggl( \frac{1}{\sqrt{T}} \int_0^T \sqrt{Y_s} \, \dd W_s,
            \frac{1}{\sqrt{T}} \int_0^T \sqrt{Y_s} \, \dd B_s,
            \frac{1}{T} \int_0^T \frac{\dd B_s}{\sqrt{Y_s}},
            \frac{1}{T^2} \int_0^T \frac{\dd s}{Y_s} \biggr) \\
   &\distr
    \Bigl( \sqrt{\theta} \, Z_3, \sqrt{\theta} \, Z_4, Z_2 \sqrt{\tau^*}, \tau^* \Bigr)
    \qquad \text{as \ $T \to \infty$.}
  \end{aligned}
 \end{align}
Using the SDE \eqref{j_Heston_SDE},
 \begin{equation}\label{intYW}
  \sigma \int_0^T \sqrt{Y_s} \, \dd W_s
  = Y_T - y_0 - \theta \kappa T + \kappa \int_0^T Y_s \, \dd s ,
 \end{equation}
 consequently, \ $\int_0^T \sqrt{Y_s} \, \dd W_s$ \ is measurable with respect to the
 $\sigma$-algebra \ $\sigma((Y_s)_{s\in[0, T]})$.
\ For all \ $(u_1, u_2, u_3, v) \in \RR^4$ \ and \ $T \in (0, \infty)$, \ we
 have
 \begin{align*}
  &\EE\Biggl( \exp\Biggl\{ \frac{\ii u_1}{\sqrt{T}}
                           \int_0^T \sqrt{Y_s} \, \dd W_s
                           + \frac{\ii u_2}{\sqrt{T}}
                             \int_0^T \sqrt{Y_s} \, \dd B_s
                           + \frac{\ii u_3}{T}
                             \int_0^T \frac{\dd B_s}{\sqrt{Y_s}}
                           + \frac{\ii v}{T^2}
                             \int_0^T \frac{\dd s}{Y_s} \biggr\} \Biggr) \\
  &= \EE\Biggl[ \EE\Biggl( \exp\Biggl\{ \frac{\ii u_1}{\sqrt{T}}
                                        \int_0^T \sqrt{Y_s} \, \dd W_s
                                        + \frac{\ii u_2}{\sqrt{T}}
                                          \int_0^T \sqrt{Y_s} \, \dd B_s
                                        + \frac{\ii u_3}{T}
                                          \int_0^T \frac{\dd B_s}{\sqrt{Y_s}}
                                        + \frac{\ii v}{T^2}
                                          \int_0^T \frac{\dd s}{Y_s} \biggr\}
                           \Bigg| \, (Y_s)_{s\in[0,T]} \Biggr) \Biggr]
  \end{align*}
  \begin{align*}
  &= \EE\Biggl[\exp\Biggl\{ \frac{\ii u_1}{\sqrt{T}}
                            \int_0^T \sqrt{Y_s} \, \dd W_s
                            + \frac{\ii v}{T^2}
                              \int_0^T \frac{\dd s}{Y_s} \Biggr\} \\
  &\phantom{= \EE\Biggl[}
               \times
               \EE\Biggl( \exp\Biggl\{ \ii
                                       \int_0^T
                                        \biggl(\frac{u_2\sqrt{Y_s}}{\sqrt{T}}
                                               + \frac{u_3}{T\sqrt{Y_s}}\biggr)
                                        \dd B_s \Biggr\}
                                       \Bigg| \, (Y_s)_{s\in[0,T]} \Biggr) \Biggr] \\
  &= \EE\Biggl(\exp\Biggl\{ \frac{\ii u_1}{\sqrt{T}}
                            \int_0^T \sqrt{Y_s} \, \dd W_s
                            + \frac{\ii v}{T^2}
                              \int_0^T \frac{\dd s}{Y_s} \Biggr\}
               \exp\Biggl\{ - \frac{1}{2}
                              \int_0^T
                               \biggl( \frac{u_2^2 Y_s}{T} + \frac{u_3^2}{T^2 Y_s}
                                       + \frac{2 u_2 u_3}{T^{3/2}} \biggr)
                               \dd s \Biggr\} \Biggr) \\
  &= \exp\Biggl\{-\frac{u_2 u_3}{\sqrt{T}}\Biggr\}
     \EE\Biggl(\exp\Biggl\{ \frac{\ii u_1}{\sqrt{T}}
                            \int_0^T \sqrt{Y_s} \, \dd W_s
                            - \frac{u_2^2}{2T} \int_0^T Y_s \, \dd s
                            - \frac{u_3^2}{2T^2} \int_0^T \frac{\dd s}{Y_s}
                            + \frac{\ii v}{T^2}
                              \int_0^T \frac{\dd s}{Y_s} \Biggr\} \Biggr) ,
 \end{align*}
 where we used the independence of \ $Y$ \ and \ $B$ \ yielding that the conditional distribution
  of \ $\int_0^T \bigl(\frac{u_2\sqrt{Y_s}}{\sqrt{T}} + \frac{u_3}{T\sqrt{Y_s}}\bigr) \dd B_s$ \ given
  \ $(Y_s)_{s\in[0,T]}$ \ is normal.
We may and do suppose that \ $Z_2$ \ and \ $Z_4$ \ are independent also from
 \ $(Y_s)_{s\in[0,T]}$.
\ Then, in a similar way, for all \ $(u_1, u_2, u_3, v) \in \RR^4$ \ and
 \ $T \in (0, \infty)$, \ we have
 \begin{align*}
  &\EE\Biggl( \exp\Biggl\{ \frac{\ii u_1}{\sqrt{T}}
                           \int_0^T \sqrt{Y_s} \, \dd W_s
                           + \frac{\ii u_2 Z_4}{\sqrt{T}}
                             \biggl(\int_0^T Y_s \, \dd s\biggr)^{1/2}
                           + \frac{\ii u_3 Z_2}{T}
                             \biggl(\int_0^T \frac{\dd s}{Y_s}\biggr)^{1/2}
                           + \frac{\ii v}{T^2}
                             \int_0^T \frac{\dd s}{Y_s} \biggr\} \Biggr) \\
  &= \EE\Biggl[ \EE\Biggl( \exp\Biggl\{ \frac{\ii u_1}{\sqrt{T}}
                                        \int_0^T \sqrt{Y_s} \, \dd W_s
                                        + \frac{\ii u_2 Z_4}{\sqrt{T}}
                                          \biggl(\int_0^T Y_s \, \dd s\biggr)^{1/2}
                                        + \frac{\ii u_3 Z_2}{T}
                                          \biggl(\int_0^T \frac{\dd s}{Y_s}\biggr)^{1/2} \\
  &\phantom{= \EE\Biggl[ \EE\Biggl( \exp\Biggl\{}
                                        + \frac{\ii v}{T^2}
                                          \int_0^T \frac{\dd s}{Y_s} \biggr\}
                           \Bigg| \, (Y_s)_{s\in[0,T]} \Biggr) \Biggr] \\
  &= \EE\Biggl[\exp\Biggl\{ \frac{\ii u_1}{\sqrt{T}}
                            \int_0^T \sqrt{Y_s} \, \dd W_s
                            + \frac{\ii v}{T^2}
                              \int_0^T \frac{\dd s}{Y_s} \Biggr\} \\
  &\phantom{= \EE\Biggl[}
               \times
               \EE\Biggl( \exp\Biggl\{ \frac{\ii u_2 Z_4}{\sqrt{T}}
                                       \biggl(\int_0^T Y_s \, \dd s\biggr)^{1/2}
                                       + \frac{\ii u_3 Z_2}{T}
                                         \biggl(\int_0^T \frac{\dd s}{Y_s}\biggr)^{1/2} \Biggr\}
                                       \Bigg| \, (Y_s)_{s\in[0,T]} \Biggr) \Biggr] \\
  &= \EE\Biggl(\exp\Biggl\{ \frac{\ii u_1}{\sqrt{T}}
                            \int_0^T \sqrt{Y_s} \, \dd W_s
                            + \frac{\ii v}{T^2}
                              \int_0^T \frac{\dd s}{Y_s} \Biggr\}
               \exp\Biggl\{ - \frac{u_2^2}{2T} \int_0^T Y_s \, \dd s
                            - \frac{u_3^2}{2T^2}
                              \int_0^T \frac{\dd s}{Y_s} \Biggr\} \Biggr) ,
 \end{align*}
 which is the same as the previous expectation except the factor
 \ $\exp\{-\frac{u_2 u_3}{\sqrt{T}}\}$.
\ Ben Alaya and Kebaier \cite[proof of Theorem 7]{BenKeb2} proved
 \[
   \biggl( \frac{Y_T - \theta \kappa T + \kappa \int_0^T Y_s \, \dd s}{\sqrt{T}},
           \frac{1}{T^2} \int_0^T \frac{\dd s}{Y_s} \biggr)
   \distr \biggl(\frac{\sigma^2}{\sqrt{2\kappa}} \, Z_3, \tau^* \biggr) \qquad
   \text{as \ $T \to \infty$,}
 \]
 hence, by \eqref{intYW},
 \[
   \biggl( \frac{1}{\sqrt{T}} \int_0^T \sqrt{Y_s} \, \dd W_s,
           \frac{1}{T^2} \int_0^T \frac{\dd s}{Y_s} \biggr)
   \distr \biggl(\sqrt{\theta} \, Z_3, \tau^* \biggr) \qquad
   \text{as \ $T \to \infty$.}
 \]
By Slutsky's lemma, we obtain
 \begin{equation}\label{BKconv2}
  \biggl( \frac{1}{\sqrt{T}} \int_0^T \sqrt{Y_s} \, \dd W_s, Z_2, Z_4 ,
          \frac{1}{T} \int_0^T Y_s \, \dd s,
          \frac{1}{T^2} \int_0^T \frac{\dd s}{Y_s} \biggr)
  \distr \biggl(\sqrt{\theta} \, Z_3, Z_2, Z_4, \theta, \tau^* \biggr)
 \end{equation}
 as \ $T \to \infty$.
Using the continuity theorem, we obtain
 \begin{align*}
  &\EE\Biggl( \exp\Biggl\{ \frac{\ii u_1}{\sqrt{T}}
                           \int_0^T \sqrt{Y_s} \, \dd W_s
                           + \frac{\ii u_2}{\sqrt{T}}
                             \int_0^T \sqrt{Y_s} \, \dd B_s
                           + \frac{\ii u_3}{T}
                             \int_0^T \frac{\dd B_s}{\sqrt{Y_s}}
                           + \frac{\ii v}{T^2}
                             \int_0^T \frac{\dd s}{Y_s} \biggr\} \Biggr) \\
  &=\exp\Biggl\{-\frac{u_2 u_3}{\sqrt{T}}\Biggr\}
    \EE\Biggl(\exp\Biggl\{ \frac{\ii u_1}{\sqrt{T}}
                            \int_0^T \sqrt{Y_s} \, \dd W_s
                            - \frac{u_2^2}{2T} \int_0^T Y_s \, \dd s
                            - \frac{u_3^2}{2T^2}
                              \int_0^T \frac{\dd s}{Y_s}
                            + \frac{\ii v}{T^2}
                              \int_0^T \frac{\dd s}{Y_s} \Biggr\} \Biggr)
 \end{align*}
 \begin{align*}
  &=\exp\Biggl\{-\frac{u_2 u_3}{\sqrt{T}}\Biggr\} \\
  &\quad\times
    \EE\Biggl(\exp\Biggl\{ \frac{\ii u_1}{\sqrt{T}}
                            \int_0^T \sqrt{Y_s} \, \dd W_s
                            + \frac{\ii u_2 Z_4}{\sqrt{T}}
                              \biggl(\int_0^T Y_s \, \dd s\biggr)^{1/2}
                            + \frac{\ii u_3 Z_2}{T}
                              \biggl(\int_0^T \frac{\dd s}{Y_s}\biggr)^{1/2}
                            + \frac{\ii v}{T^2}
                              \int_0^T \frac{\dd s}{Y_s} \Biggr\} \Biggr) \\
  &\to \EE\Biggl( \exp\Biggl\{ \ii u_1 \sqrt{\theta} \, Z_3
                               + \ii u_2 \sqrt{\theta} \, Z_4
                               + \ii u_3 Z_2 \sqrt{\tau^*}
                               + \ii v \tau^* \biggr\} \Biggr) \qquad
  \text{as \ $T \to \infty$}
 \end{align*}
 for all \ $(u_1, u_2, u_3, v) \in \RR^4$.
\ By the continuity theorem, we obtain \eqref{bGbh1}.

With a random scaling, we have
 \[
   \begin{bmatrix}
    \frac{\sigma T^2 }{2\sqrt{1-\varrho^2}\,(\int_0^T Y_u\,\dd u)^{3/2}}
    \, (\htheta_T - \theta) \\[4mm]
    \frac{1}{\sigma\sqrt{1-\varrho^2}}
    \bigl(\int_0^T Y_u \, \dd u\bigr)^{1/2} \, (\hkappa_T - \kappa) \\
    \frac{\sigma T^2}{2\int_0^T Y_u \, \dd u} \, (\hmu_T - \mu)
   \end{bmatrix}
   = \ttbC_T (\hbpsi_T - \bpsi), \qquad T \in (0, \infty) ,
 \]
 with
 \[
   \ttbC_T
   := \begin{bmatrix}
       \frac{\sigma T^2}{2\sqrt{1-\varrho^2}\,(\int_0^T Y_u\,\dd u)^{3/2}}
        & 0 & 0 \\
       0 & \frac{1}{\sigma\sqrt{1-\varrho^2}}
           \bigl(\int_0^T Y_u \, \dd u\bigr)^{1/2} & 0 \\
       0 & 0 & \frac{\sigma T^2}{2\int_0^T Y_u \, \dd u}
      \end{bmatrix} .
 \]
We have
 \[
   \ttbC_T (\hbpsi_T - \bpsi)
   = (\ttbC_T \tbC_T^{-1}) \tbC_T (\hbpsi_T - \bpsi) ,
 \]
 where
 \[
   \ttbC_T \tbC_T^{-1}
   = \begin{bmatrix}
      \frac{\sigma}
           {2\sqrt{1-\varrho^2}\,(\frac{1}{T}\int_0^T Y_u\,\dd u)^{3/2}}
       & 0 & 0 \\
      0 & \frac{1}{\sigma\sqrt{1-\varrho^2}}
          \bigl(\frac{1}{T}\int_0^T Y_u \, \dd u\bigr)^{1/2} & 0 \\
      0 & 0 & \frac{\sigma}{\frac{2}{T}\int_0^T Y_u \, \dd u}
     \end{bmatrix} .
 \]
Applying \eqref{SLLN=}, \eqref{AK}, \eqref{BKconv1}, \eqref{bGbh1}, Slutsky's lemma
 (or part (v) of Theorem 2.7 in van der Vaart \cite{Vaart}) and the continuous mapping
 theorem, we obtain
 \[
    \bigl( \ttbC_T \tbC_T^{-1}, \tbC_T (\hbpsi_T - \bpsi)\bigr)
   \distr
   \left(\begin{bmatrix}
          \frac{\sqrt{2\kappa^3}}{\sigma^2\sqrt{1-\varrho^2}} & 0 & 0 \\
          0 & \frac{1}{\sqrt{2\kappa(1 - \varrho^2)}} & 0 \\
          0 & 0 & \frac{\kappa}{\sigma}
         \end{bmatrix} ,
         \begin{bmatrix}
          -\frac{\sigma^2 Z_1}{\sqrt{2\kappa^3}} \, \sqrt{1-\varrho^2} \\
          Z_1 \sqrt{2(1 - \varrho^2)\kappa} \\
          \frac{\varrho\sigma}{\kappa\tau}
          + \frac{Z_2\sigma\sqrt{1-\varrho^2}}{\kappa\sqrt{\tau}}
         \end{bmatrix}\right) \qquad \text{as \ $T \to \infty$.}
 \]
Using again the continuous mapping theorem, we obtain \eqref{abm=r}.
\proofend

\begin{Rem}
Putting formally \ $\theta\kappa = \frac{\sigma^2}{2}$ \ into the formula of \ $\bV$
 \ given in \eqref{exp_G_infty_inverse_explicit} of Theorem \ref{Thm_MLE}, one can
 observe that the joint limit distribution of the first two coordinates in
 \eqref{MLE_sub} of Theorem \ref{Thm_MLE} and in \eqref{abm=} of Theorem \ref{Thm_MLE=}
 coincide.
 \proofend
\end{Rem}

\begin{Rem}
According to Theorem 7 in Ben Alaya and Kebaier \cite{BenKeb2}, if \ $a=\frac{\sigma_1^2}{2}$ \ and \ $b\in(0,\infty)$,
 \ then, based on continuous time observations \ $(Y_t)_{t\in[0,T]}$, $T\in(0,\infty)$, \ for the MLE \ $(\ha_T,\hb_T)$ \
 of \ $(a,b)$ \ for the first coordinate process of the SDE \eqref{Heston_SDE}, we have
 \begin{align}\label{help_BA}
  \begin{bmatrix}
    T(\ha_T-a) \\
    T^{1/2}(\hb_T - b) \\
  \end{bmatrix}
  \distr
  \begin{bmatrix}
    \frac{\sigma_1^2}{b\tau} \\
      \sqrt{2b} Z_1 \\
  \end{bmatrix}
  \qquad \text{as \ $T\to\infty$,}
 \end{align}
 where \ $Z_1$ \ is a standard normally distributed random variable independent of \ $\tau$ \ introduced
 in Theorem \ref{Thm_MLE=}.
Hence, using Slutsky's lemma and that \ $\hb_T$ \ converges in probability to \ $b$ \ as \ $T\to\infty$ \
 (following from \eqref{help_BA}), we get
 \begin{align*}
  T^{1/2}\left(\frac{\ha_T}{\hb_T} - \frac{a}{b}\right)
    & = T^{1/2} \frac{b\ha_T - a\hb_T}{b \hb_T}
      = T^{1/2} \frac{b(\ha_T-a) - a (\hb_T-b)}{b \hb_T}\\
    & = \frac{T^{-1/2} bT(\ha_T-a) - a T^{1/2} (\hb_T-b)}{b \hb_T}
     \distr -\frac{a}{b^2}\sqrt{2b} Z_1
             = - \frac{\sqrt{2}a}{\sqrt{b^3}} Z_1
 \end{align*}
 as \ $T\to\infty$.
\ Let us observe that in the special case of \ $\varrho=0$, \ we have
 \ $\frac{\ha_T}{\hb_T}  = \htheta_T$ \ and \ $\hb_T = \kappa_T$, \ $T>0$ \ (for the
 explicit formulae for \ $\ha_T$ \ and \ $\hb_T$, \ see Ben Alaya and Kebaier
 \cite[Section 3.1]{BenKeb2}).
Moreover, in case of \ $a=\theta\kappa=\frac{\sigma^2}{2}$ \ and \ $b=\kappa,$ \ we have
 \ $-\frac{\sqrt{2}a}{\sqrt{b^3}} = - \frac{\sigma^2}{\sqrt{2\kappa^3}}$.
\ Hence, under the conditions of Theorem \ref{Thm_MLE=} together with \ $\varrho=0$, \
 the joint (weak) convergence of the first two coordinates of \eqref{abm=} follows from Theorem 7
 in Ben Alaya and Kebaier \cite{BenKeb2}.
 \proofend
\end{Rem}

\section{Numerical illustrations}\label{section_simulation}

We present some numerical illustrations in order to confirm our limit theorems given in Sections
 \ref{section_CMLE} and \ref{section_AMLE}.
We call the attention to the fact that our numerical illustrations using
 synthetic data can not be considered as simulations
 or a receipt for handling real data set of \ $(Y, S)$, \
 since, as it will turn out, we use the standard Wiener processes \ $(W_t)_{t\in[0,\infty)}$ \ and
 \ $(B_t)_{t\in[0,\infty)}$ \ appearing in \eqref{j_Heston_SDE} that can not be observed.
Hence the main aim of this section is to confirm the scaling factors and the limit distributions of the derived MLE in Theorems \ref{Thm_MLE} and \ref{Thm_MLE=}.
\ In order to approximate the estimator \ $\hbpsi_T$ \ given in \eqref{MLE}, one could
 generate sample paths of the model \eqref{j_Heston_SDE}, and then one could approximate the estimator \ $\hbpsi_T$ \ given in \eqref{MLE} based on the generated sample paths.
For this, it would be sufficient to simulate, for a large time \ $T>0$, \ the random variables
\begin{align*}
 &Y_T, \qquad I_{1,T} :=\int_0^T Y_u \, \dd u, \qquad I_{2,T}:=\int_0^T \frac{\dd u}{Y_u},\qquad I_{3,T}:=\int_0^T \frac{\dd Y_u}{Y_u},\\
 &I_{4,T}:=\int_0^T \frac{\dd S_u - S_{u-} \, \dd L_u}{S_{u-}},\qquad I_{5,T}:=\int_0^T \frac{\dd S_u - S_{u-} \, \dd L_u}{Y_u S_{u-}}.
\end{align*}
It is well known that the random variable \ $Y_T$ \ has a non-central chi-squared
 distribution (see, e.g., Alfonsi \cite[Proposition 1.2.11]{Alf2}) that can be simulated exactly.
Further, Broadie and Kaya \cite[Section 3.2]{BroKay} proposed an exact simulation method of \ $(Y_T ,I_{1,T})$, \ and
 more recently, Ben Alaya and Kebaier \cite[Sections 4.1 and 4.2]{BenKeb1} developed an analogous method to simulate \ $(Y_T, I_{2,T})$.
\ In the context of our current study, it would be possible to compute the Laplace transform of the couple \ $(I_{1,T}, I_{2,T})$ \ conditionally
 on \ $Y_T$, \ and using relation \eqref{help_ITO}, we could derive an exact simulation method for the random vector \ $(Y_T, I_{1,T}, I_{2,T}, I_{3,T})$.
\ However, due to the lack of an exact simulation method for the couple
 \ $(I_{4,T}, I_{5,T})$, \ we choose to approximate the quantities
 \ $(Y_T, I_{1,T}, I_{2,T}, I_{3,T}, I_{4,T}, I_{5,T})$ \ using discretization schemes, like the famous Euler one
 (see, e.g., Kloeden and Platen \cite{KloPla} or Alfonsi \cite[Chapter 2]{Alf2}).
Nevertheless, it is important to note that the discretization of the CIR process presents some troubles
 because of the square root in the diffusion coefficient. Several papers deal with this problem, see for example Alfonsi \cite{Alf}
  and Berkaoui et al. \cite{BerBosDio}.

For a given time step \ $\frac{T}{n}$ \ with \ $n\in\mathbb N$, \ we use the drift
 implicit Euler scheme introduced by Alfonsi \cite{Alf}
 to approximate the process $(Y_t)_{t\in[0,T]}$ at times $t^n_i=i\frac{T}{n}$, $i\in\{0,\dots,n\}$,
 by the following non-linear recursion, \ $Y^n_0=y_0\in(0,\infty)$ \ and
\begin{equation*}
 Y^n_{t^n_{i+1}}=\left(\frac{ \frac{\sigma}{2} (W_{t^n_{i+1}}- W_{t^n_{i}})
                      + \sqrt{Y^n_{t^n_{i}}}+\sqrt{ (\frac{\sigma}{2} (W_{t^n_{i+1}}- W_{t^n_{i}})
                      +\sqrt{ Y^n_{t^n_{i}}})^2+(1+ \frac{\kappa T}{2n})(2\theta\kappa-\frac{\sigma^2}{2})\frac{T}{n}}}
 {2+\frac{\kappa T}{n}}\right)^2,
\end{equation*}
 for \ $i\in\{0,\dots,n-1\}$.
Note that, due to Alfonsi \cite{Alf}, this scheme is well defined for
 \ $\theta, \kappa \in (0, \infty)$ and \ $\theta \kappa \in ( \frac{\sigma^2}{4}, \infty \bigr)$
 \ covering the case \ $\theta \kappa \in [ \frac{\sigma^2}{2}, \infty \bigr)$ \ as well, which ensures
 the unique existence of a MLE of \ $(\theta,\kappa,\mu)$, \ see Proposition \ref{Pro_MLE}.
Moreover, the strong convergence rate of this approximation is of order $1$ \
 in case of \ $\theta\kappa\in(\frac{\sigma^2}{2},\infty)$, \ see Alfonsi \cite{Alf} for more details.
Then, we can easily approximate \ $I_{1,T}$, $I_{2,T}$ \ and \ $I_{3,T}$ \ respectively, by
 \[
    I^n_{1,T}:=\frac{T}{n}\sum_{i=0}^{n-1}Y^n_{t^n_i},
    \qquad I^n_{2,T}:=\frac{T}{n}\sum_{i=0}^{n-1}\frac{1}{Y^n_{t^n_i}}, \qquad
    I^n_{3,T}:=\sum_{i=0}^{n-1}\frac{Y^n_{t^n_{i+1}}-Y^n_{t_{i}}}{Y^n_{t^n_i}} .
 \]
Alternatively, using the relation \eqref{help_ITO}, one can also use the approximation
 \[
   \tI^n_{3,T}:=\log(Y^n_{t^n_n}) - \log(y_0) + \frac{\sigma^2}{2} I^n_{2,T} .
 \]
Since we just would like to present some numerical illustrations of our limit theorems and
 not to provide simulations, we will not approximate the processes \ $(L_t)_{t\in[0,\infty)}$ \ and
 \ $(S_t)_{t\in[0,\infty)}$, \ instead, applying the equations \eqref{KEY}, we can use
 \begin{align*}
  I^n_{4,T} &:= \mu T + \sum_{i=0}^{n-1}
                         \sqrt{Y^n_{t^n_i}}
                         (\varrho (W_{t^n_{i+1}} - W_{t^n_i})
                          + \sqrt{1 -\varrho^2} (B_{t^n_{i+1}} - B_{t^n_i})) , \\
  I^n_{5,T} &:= \mu I^n_{2,T} + \sum_{i=0}^{n-1}
                                 \frac{\varrho (W_{t^n_{i+1}} - W_{t^n_i})
                                  + \sqrt{1 -\varrho^2} (B_{t^n_{i+1}} - B_{t^n_i})}
                                       {\sqrt{Y^n_{t^n_i}}} .
 \end{align*}
Here, we point out that \ $I^n_{4,T}$ \ and \ $I^n_{5,T}$ \ use the standard Wiener processes
 \ $(W_t)_{t\in[0,\infty)}$ \ and \ $(B_t)_{t\in[0,\infty)}$ \ appearing in \eqref{j_Heston_SDE}
 that can not be observed, so \ $I^n_{4,T}$ \ and \ $I^n_{5,T}$ \ can not be used for approximating
 \ $I_{4,T}$ \ and \ $I_{5,T}$, \ respectively, given a real dataset of \ $(Y, S)$.
\ However, the main advantage of this procedure is that it allows us to handle numerical
 illustrations involving any arbitrary purely non-Gaussian
 L\'evy process \ $(L_t)_{t\in\RR_+}$ \ with L\'evy--Khintchine representation given
 in \eqref{LK}.
Hence, by \eqref{MLE_coordinatewise}, we approximate \ $\htheta_T$, \ $\hkappa_T$ \ and \ $\hmu_T$ \ by
\begin{align*}
  \htheta^n_T
  &:= \frac{ I^n_{1,T}I^n_{2,T}I^n_{3,T}
           - T I^n_{2,T}(Y^n_T-  y_0)
           + \varrho \sigma T I^n_{2,T}I^n_{4,T}
           - \varrho \sigma T^2 I^n_{5,T}}
          {T I^n_{2,T}I^n_{3,T}
           - \bigl(I^n_{2,T})^2 (Y^n_T- y_0 )
           + \varrho \sigma \bigl(I^n_{2,T})^2I^n_{4,T}
           - \varrho \sigma T I^n_{2,T}I^n_{5,T}}, \\
  \hkappa^n_T
  &:= \frac{T I^n_{3,T}
           -  I^n_{2,T}(Y^n_T - y_0)
           + \varrho \sigma I^n_{2,T}  I^n_{4,T}
           - \varrho \sigma T  I^n_{5,T}}
         {\bigl(  I^n_{1,T} I^n_{2,T} - T^2\bigr)} , \\
  \hmu^n_T
  &:= \frac{I^n_{5,T}}
          {I^n_{2,T} } .
 \end{align*}
For the numerical implementation, we consider two case studies, one with
 \ $\theta \kappa > \sigma^2/2$, \ and another with \ $\theta \kappa = \sigma^2/2$.

\begin{figure}[ht]
 \centering
 \includegraphics[width=16cm,height=7cm]{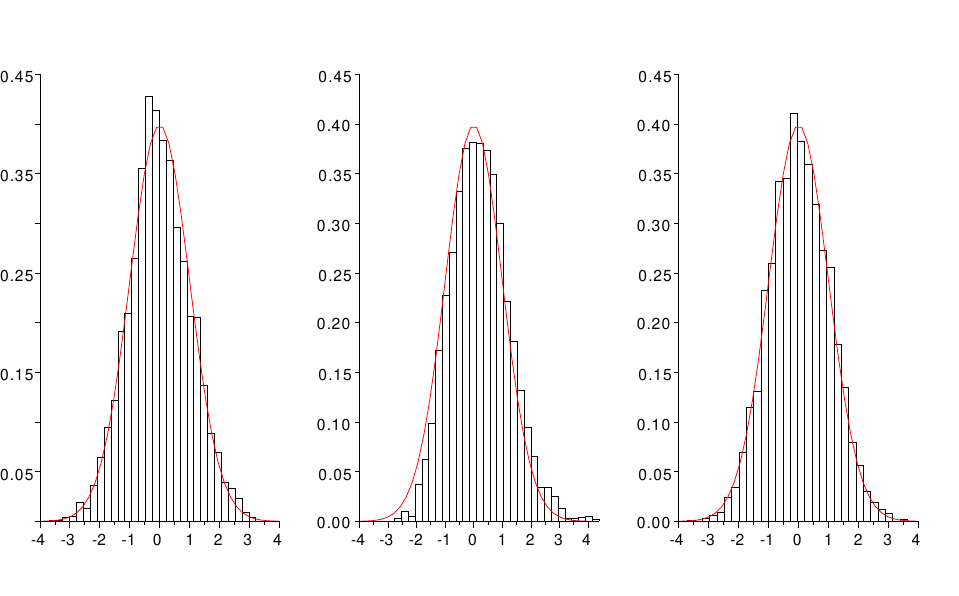}
 \caption{From the left to the right, the density histograms of the suitably scaled errors given in \eqref{error1}, \eqref{error2} and \eqref{error3}. In each case, the red line denotes the density function of the standard normal distribution.}
 \label{fig1}
\end{figure}

First we take \ $\theta=2$, \ $\kappa=0.5$,
 \ $\mu=1-\sqrt{\ee}$, \ $\sigma=0.2$, \ $\varrho=0.5$, \ $y_0=1$, \ $s_0=100$,
 \ $\frac{T}{n}=0.01$, \ and we simulate \ $M=4000$ \ independent trajectories of the normalized error \ $T^{1/2}(\hbpsi_T - \bpsi)$.
\ Note that \ $\theta\kappa > \frac{\sigma^2}{2}$ \ with this choice of parameters.
In Table \ref{Table1} we give the relative errors for \ $T \in \{10, 100, 300\}$.
\ Note that, when $T$ increases we need of course a suitable number of time steps $n$ to guarantee a good approximation.
\begin{table}[ht]
\centering
\begin{tabular}{|c|c|c|c|}
  \hline
   Relative error & $T=10$&  $T=100$ & $T=300$\\
  \hline
   $|\htheta^n_T-\theta|/\theta$   & $0.0010578$ & $0.0002387$ & $0.0000658$ \\
   $|\hkappa^n_T-\kappa|/\kappa$  & $0.2803024$ & $0.0532183$ & $0.0214441$ \\
   $|\hmu^n_T-\mu|/\mu$  & $0.0380512$ & $0.0060456$ & $0.0034771$ \\
  \hline
\end{tabular}
\caption{Relative errors.}
\label{Table1}
\end{table}
The obtained relative errors confirm the strong consistency of the estimator $\hbpsi_T$ stated in Theorem \ref{Thm_MLE_cons}.
In Figure \ref{fig1}
 we illustrate the law of each suitably scaled coordinate of the MLE
 \ $\hbpsi_T = (\htheta_T,\hkappa_T,\hmu_T)$ \ for \ $T = 300$.
As a consequence of Theorem \ref{Thm_MLE}, we have
 \begin{align}
    &\sqrt{\frac{2\kappa^3T}{\sigma^2(2\theta\kappa - \varrho^2\sigma^2)}} \,
    (\htheta_T -\theta) \distr \cN(0,1)\qquad \text{as \ $T\to\infty$,} \label{error1} \\
    &\sqrt{\frac{T}{2\kappa(1-\varrho^2)}} \, (\hkappa_T -\kappa) \distr \cN(0,1)\qquad \text{as \ $T\to\infty$,}\label{error2}\\
    &\sqrt{\frac{2\kappa T}{2\theta\kappa - \sigma^2}} \, (\hmu_T -\mu) \distr \cN(0,1)\qquad \text{as \ $T\to\infty$.}\label{error3}
 \end{align}
The obtained density histograms in Figure \ref{fig1} confirm our results in Theorem \ref{Thm_MLE}.

\begin{figure}[ht]
 \centering
 \includegraphics[width=16cm,height=7cm]{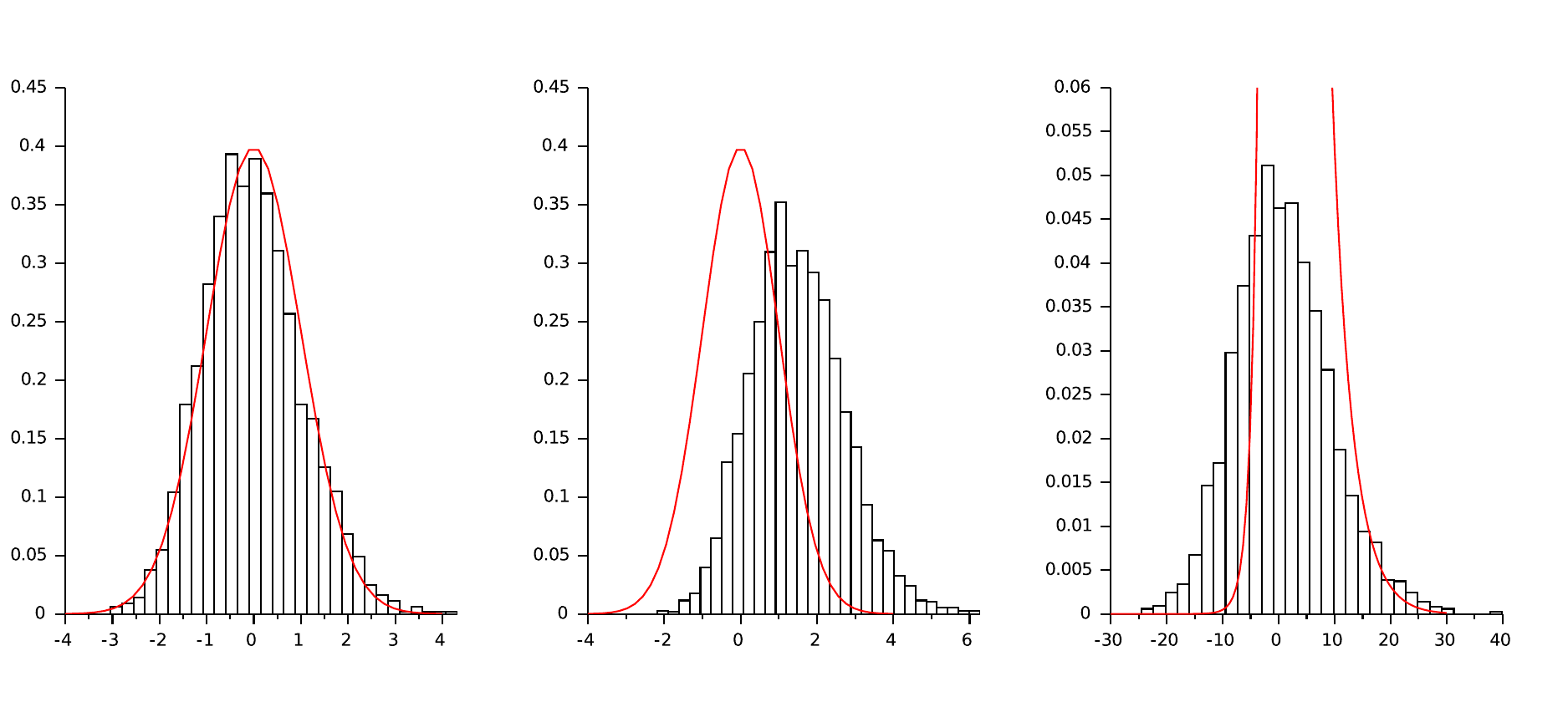}
 \caption{From the left to the right, the density histograms of the suitably scaled errors given in \eqref{error1=}, \eqref{error2=} and \eqref{error3=}. The red line denotes the density functions of corresponding limit distributions.}
 \label{fig2}
\end{figure}

Next we take \ $\theta=2$, \ $\kappa=0.5$,
 \ $\mu=1-\sqrt{\ee}$, \ $\sigma=\sqrt{2}$, \ $\varrho=0.5$, \ $y_0=1$, \ $s_0=100$,
 \ $n=30000$, \ and we simulate \ $M=4000$ \ independent trajectories of the appropriately
 normalized error \ $\hbpsi_T - \bpsi$.
\ Note that \ $\theta\kappa = \frac{\sigma^2}{2}$ \ with this choice of parameters.
In Figure \ref{fig2}
 we illustrate the law of each suitably scaled coordinate of the MLE
 \ $\hbpsi_T = (\htheta_T,\hkappa_T,\hmu_T)$ \ for \ $T = 300$.
As a consequence of Theorem \ref{Thm_MLE=}, we have
 \begin{align}
    &\sqrt{\frac{2\kappa^3T}{\sigma^4(1-\varrho^2)}} \,
    (\htheta_T -\theta) \distr \cN(0,1)\qquad \text{as \ $T\to\infty$,} \label{error1=} \\
    &\sqrt{\frac{T}{2\kappa(1-\varrho^2)}} \, (\hkappa_T -\kappa) \distr \cN(0,1)\qquad \text{as \ $T\to\infty$,}\label{error2=}\\
    &T (\hmu_T -\mu) \distr \frac{\varrho\sigma}{\kappa\tau}
   + \frac{\sigma\sqrt{1-\varrho^2}}{\kappa\sqrt{\tau}} Z_2\qquad \text{as \ $T\to\infty$.}\label{error3=}
 \end{align}
We plot the density function of the limit distribution in \eqref{error3=} using its explicit form
 given in Appendix \ref{App_density}.
Note that it is cutted at the level \ $0.06$, \ since it tends to infinity at \ $0$.
In case of the parameter \ $\kappa$, \ one can see a bias in Figure \ref{fig2}, which,
 in our opinion, is caused by the bad performance of the applied discretization scheme
 together with the approximation method of the integrals in question,
 when \ $\theta\kappa = \frac{\sigma^2}{2}$.
\ We have not been able to find any discretization scheme to explain the bias
 (we tried the truncated Euler scheme, see, e.g., Deelstra and Delbaen \cite{DeDe},
  and the symmetrized Euler scheme, see, e.g., Diop \cite{Dio} or Berkaoui
  et al.\ \cite{BerBosDio}).
Eventually, this bad performance can also be observed whenever the ratio
 \ $\frac{2\theta\kappa}{\sigma^2}$ \ is close to 1.
And to top it all, one can observe the same phenomena already in case of the MLE
 \ $(\ha_T,\hb_T)$ \ of \ $(a,b)$ \ of the first coordinate process of the
 SDE \eqref{Heston_SDE} based on continuous time observations \ $(Y_t)_{t\in[0,T]}$,
 \ $T\in(0,\infty)$, \ for both \ $\ha_T$ \ and \ $\hb_T$ \ (for an expression of
  \ $(\ha_T,\hb_T)$, \ see Overbeck \cite{Ove}).
So we conclude that the bias for \ $\kappa$ \ seen in Figure \ref{fig2} is not related
 to the fact that the model \eqref{j_Heston_SDE} contains a jump part.

\begin{figure}[ht]
 \centering
 \includegraphics[width=16cm,height=7cm]{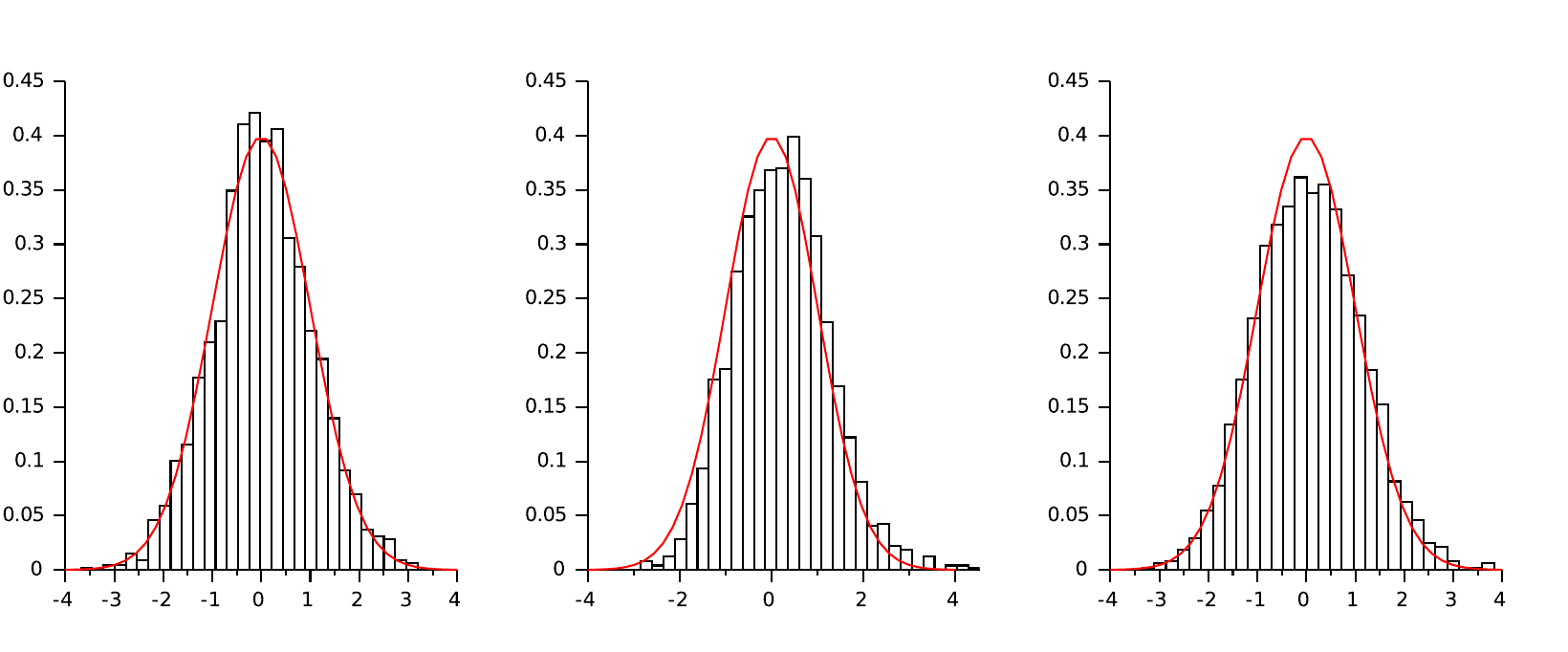}
 \caption{From the left to the right, the density histograms of the suitably scaled errors given in \eqref{error1}, \eqref{error2} and \eqref{error3}. In each case, the red line denotes the density function of the standard normal distribution.}
 \label{fig3}
\end{figure}

As we mentioned in the Introduction, the model \eqref{j_Heston_SDE} with \ $L$ \ as a
 compound Poisson process given in \eqref{L_special} is quite popular in finance.
In this special case, one can use another illustration method without applying
 the equations \eqref{KEY}, but still using the standard Wiener processes
 \ $(W_t)_{t\in[0,\infty)}$ \ and \ $(B_t)_{t\in[0,\infty)}$.
\ Namely, for all \ $0\leq s<t$, \ by \eqref{HESCPP},
 \[
  S_t = S_s\exp\biggl\{\int_s^t \Bigl(\mu - \frac{1}{2} Y_u\Bigr) \, \dd u
                    + \int_s^t
                       \sqrt{Y_u} \,
                       (\varrho\, \dd  W_u + \sqrt{1 - \varrho^2}\, \dd  B_u)
                       + \sum_{k=\pi_s+1}^{\pi_t} J_k
              \biggr\},
 \]
 hence we can approximate the price process \ $(S_t)_{t\in[0,T]}$ \ by the recursion
 \ $S^n_0=s_0\in(0,\infty)$ \ and
 \begin{align*}
 \frac{S^n_{t^n_{i+1}}}{S^n_{t^n_i}}
       = \exp\biggl\{\frac{\mu T}{n} - \frac{T}{2n} Y^n_{t^n_i}
                    + \sqrt{Y^n_{t^n_i}}(\varrho (W_{t^n_{i+1}} - W_{t^n_i})
                     + \sqrt{1 -\varrho^2} (B_{t^n_{i+1}} - B_{t^n_i})  )
                       + \sum_{k=\pi_{t^n_i}+1}^{\pi_{t^n_{i+1}}} J_k\biggr\}
 \end{align*}
 for  $i\in\{0,\dots,n-1\}$.
\ Note that the process \ $(\pi_t)_{t\in[0,\infty)}$ \ is a Poisson process with intensity
 \ $\lambda$ \ being independent of \ $(W_t,B_t)_{t\in[0,\infty)}$, \
  and it can be easily simulated.
Therefore, given independently an i.i.d. sequence of random variables
 \ $(J_k)_{k\in\NN}$, \ one can simulate at the same time the term
 \ $\sum_{k=\pi_{t^n_i}+1}^{\pi_{t^n_{i+1}}} J_k$ \ together with the increments
 \ $L_{t^n_{i+1}}-L_{t^n_i}= \sum_{k=\pi_{t^n_i}+1}^{\pi_{t^n_{i+1}}}(\ee^{J_k} - 1)$,
 \ $i\in\{0,\ldots,n-1\}$.
\ Further, one can approximate \ $I_{4,T}$ \ and \ $I_{5,T}$ \ respectively by
 \[
   \tI^n_{4,T}
   := \sum_{i=0}^{n-1}
       \left(\frac{S^n_{t^n_{i+1}}}{S^n_{t^n_{i}}} - 1\right) -L_T
    \quad \mbox{ and }\quad
   \tI^n_{5,T}
   := \sum_{i=0}^{n-1}
       \frac{1}{Y^n_{t^n_i}}
        \left(\frac{S^n_{t^n_{i+1}}}{S^n_{t^n_{i}}} - 1\right)
      - \sum_{i=0}^{n-1}\frac{L_{t^n_{i+1}}-L_{t^n_i}}{Y^n_{t^n_i}}.
 \]
We remark that \ $S^n_{t^n_i} \in (0, \infty)$, \ $n \in \NN$,
 \ $i \in \{0, 1, \ldots, n-1\}$, \ so \ $\tI^n_{4,T}$ \ and \ $\tI^n_{5,T}$ \ are
 well-defined.
Here we take again \ $\theta=2$, \ $\kappa=0.5$, \ $\mu=1-\sqrt{\ee}$, \ $\sigma=0.2$,
 \ $\varrho=0.5$, \ $y_0=1$, \ $s_0=100$, \ $n=30000$ \ and additionally, \ $\lambda = 1$ \
 and a random variable \ $J_1$ \ with standard normal distribution.
We simulate \ $M=2000$ \ independent trajectories of the normalized error \ $T^{1/2}(\hbpsi_T - \bpsi)$.
\ Note that \ $\theta\kappa > \frac{\sigma^2}{2}$ \ with this choice of parameters.
In Figure \ref{fig3} we illustrate the law of each suitably scaled coordinate of the MLE
 \ $\hbpsi_T = (\htheta_T,\hkappa_T,\hmu_T)$ \ which confirms our results in Theorem \ref{Thm_MLE}.

Finally, we note that we used the open source software Scilab for making the simulations.

\bigskip
\medskip

\appendix

\noindent{\bf\Large Appendix}

\section{Likelihood-ratio process}\label{App_LR}

Based on Jacod and Shiryaev \cite{JSh}, see also Jacod and M\'emin \cite{JM},
 S{\o}rensen \cite{SorM} and Luschgy \cite{Lus}, we recall certain sufficient
 conditions for the absolute continuity of probability measures induced by
 semimartingales together with a representation of the corresponding Radon--Nikodym
 derivative (likelihood-ratio process).

Let \ $D([0,\infty), \RR^d)$ \ denote the space of
 \ $\RR^d$-valued c\`adl\`ag functions defined on \ $[0,\infty)$.
\ Let \ $(\eta_t)_{t\in[0,\infty)}$ \ denote the canonical process
 \ $\eta_t(\omega) := \omega(t)$, \ $\omega \in D([0,\infty), \RR^d)$,
 \ $t \in [0, \infty)$.
\ Put \ $\cF_t^\eta := \sigma(\eta_s, s \in [0, t])$, \ $t \in [0, \infty)$,
 \ and
 \[
   \cD_t([0,\infty), \RR^d)
   := \bigcap_{\vare>0} \cF_{t+\vare}^\eta , \quad
   t \in [0, \infty) , \qquad
   \cD([0,\infty), \RR^d)
   := \sigma\Biggl(\,\bigcup_{t \in [0, \infty)} \cF_t^\eta\Biggr) .
 \]
Let \ $\Psi \subset \RR^k$ \ be an arbitrary non-empty set, and let \ $\PP_\bpsi$,
 \ $\bpsi \in \Psi$, \ are probability measures on the canonical space
 \ $(D([0,\infty), \RR^d), \cD([0,\infty), \RR^d))$.
\ Suppose that for each \ $\bpsi \in \Psi$, \ under \ $\PP_\bpsi$, \ the canonical
 process \ $(\eta_t)_{t\in[0,\infty)}$ \ is a semimartingale with semimartingale
 characteristics \ $(B^{(\bpsi)}, C, \nu^{(\bpsi)})$ \ associated with a
 fixed Borel measurable truncation function \ $h:\RR^d \to \RR^d$, \ see Jacod and
 Shiryaev \cite[Definition II.2.6 and Remark II.2.8]{JSh}.
Namely, \ $C_t := \langle (\eta^\cont)^{(\bpsi)} \rangle_t$,
 \ $t \in [0, \infty)$, \ where
 \ $(\langle (\eta^\cont)^{(\bpsi)} \rangle_t)_{t\in[0,\infty)}$ \ denotes
 the (predictable) quadratic variation process (with values in \ $\RR^{d\times d}$)
 \ of the continuous martingale part \ $(\eta^\cont)^{(\bpsi)}$ \ of \ $\eta$ \ under \ $\PP_\bpsi$,
 \ $\nu^{(\bpsi)}$ \ is the compensator of the integer-valued random measure
 \ $\mu^\eta$ \ on \ $[0,\infty) \times \RR^d$ \ associated with the jumps of
 \ $\eta$ \ under \ $\PP_\bpsi$ \ given by
 \begin{equation}\label{mueta}
   \mu^\eta(\omega, \dd t, \dd\bx)
   := \sum_{s\in[0,\infty)}
       \bbone_{\{\Delta\eta_s(\omega)\ne\bzero\}}
       \vare_{(s,\Delta\eta_s(\omega))}(\dd t, \dd\bx) , \qquad
   \omega \in D([0,\infty), \RR^d) ,
 \end{equation}
 where \ $\vare_{(t,\bx)}$ \ denotes the Dirac measure at the point
 \ $(t,\bx) \in [0,\infty) \times \RR^d$, \ and
 \ $\Delta\eta_t := \eta_t - \eta_{t-}$, \ $t \in (0,\infty)$,
 \ $\Delta\eta_0 := \bzero$, \ and \ $B^{(\bpsi)}$ \ is the predictable process
 (with values in \ $\RR^d$ \ having finite variation over each finite interval
 \ $[0, t]$, \ $t \in [0, \infty)$) \ appearing in the canonical decomposition
 \[
   \teta_t = \eta_0 + M^{(\bpsi)}_t + B^{(\bpsi)}_t , \qquad
   t \in [0, \infty) ,
 \]
 of the special semimartingale \ $(\teta_t)_{t\in[0,\infty)}$ \ under \ $\PP_\bpsi$
 \ given by
 \[
   \teta_t := \eta_t - \sum_{s\in[0,t]} (\eta_s - h(\Delta\eta_s)) , \qquad
   t \in [0, \infty) ,
 \]
 where \ $(M^{(\bpsi)}_t)_{t\in[0,\infty)}$ \ is a local martingale with
 \ $M^{(\bpsi)}_0 = \bzero$.
\ We call the attention that, by our assumption, the process
 \ $C= \langle (\eta^\cont)^{(\bpsi)} \rangle$ \ does not depend on \ $\bpsi$,
 \ although \ $(\eta^\cont)^{(\bpsi)}$ \ might depend on \ $\bpsi$.
\ In addition, assume that \ $\PP_\bpsi(\nu^{(\bpsi)}(\{t\} \times \RR^d) = 0)=1$
 \ for every \ $\bpsi \in \Psi$, \ $t \in [0,\infty)$, \ and
 \ $\PP_\bpsi(\eta_0 = \bx_0) = 1$ \ with some \ $\bx_0 \in \RR^d$
 \ for every \ $\bpsi \in \Psi$.
\ Note that we have the semimartingale representation
 \begin{equation}\label{semirepr}
  \begin{aligned}
   \eta_t &= \bx_0 + B^{(\bpsi)}_t + (\eta^\cont)^{(\bpsi)}_t
            + \int_0^t \int_{\RR^d} h(\bx) \,
               (\mu^\eta - \nu^{(\bpsi)})(\dd s, \dd\bx) \\
         &\quad
            + \int_0^t \int_{\RR^d} (\bx - h(\bx)) \, \mu^\eta(\dd s, \dd\bx) ,
   \qquad t \in [0, \infty) ,
  \end{aligned}
 \end{equation}
 of \ $\eta$ \ under \ $\PP_\bpsi$, \ see Jacod and Shiryaev
 \cite[Theorem II.2.34]{JSh}.
Moreover, for each \ $\bpsi \in \Psi$, \ let us choose a nondecreasing,
 continuous, adapted process \ $(F_t^{(\bpsi)})_{t\in[0,\infty)}$ \ with
 \ $F_0^{(\bpsi)} = 0$ \ and a predictable process
 \ $(c_t^{(\bpsi)})_{t\in[0,\infty)}$ \ with values in the set of all
 symmetric positive semidefinite \ $d \times d$ \ matrices such that
 \[
   C_t = \int_0^t c_s^{(\bpsi)} \, \dd F_s^{(\bpsi)}
 \]
 $\PP_\bpsi$-almost sure for every \ $t \in [0, \infty)$.
\ Due to the assumption \ $\PP_\bpsi(\nu^{(\bpsi)}(\{t\} \times \RR^d) = 0) = 1$
 \ for every \ $\bpsi \in \Psi$, \ $t \in [0, \infty)$, \ such choices of
 \ $(F_t^{(\bpsi)})_{t\in[0,\infty)}$ \ and
 \ $(c_t^{(\bpsi)})_{t\in[0,\infty)}$ \ are possible, see Jacod and
 Shiryaev \cite[Proposition II.2.9 and Corollary II.1.19]{JSh}.
Let \ $\cP$ \ denote the predictable $\sigma$-algebra on
 \ $D([0,\infty), \RR^d) \times [0,\infty)$.
\ Assume also that for every \ $\bpsi, \btpsi \in \Psi$, \ there exist a
 \ $\cP \otimes \cB(\RR^d)$-measurable function
 \ $V^{(\btpsi,\bpsi)}
    : D([0,\infty), \RR^d) \times [0, \infty) \times \RR^d \to (0, \infty)$
 \ and a predictable \ $\RR^d$-valued process \ $\beta^{(\btpsi,\bpsi)}$
 \ satisfying
 \begin{gather}
  \nu^{(\bpsi)}(\dd t, \dd\bx)
  = V^{(\btpsi,\bpsi)}(t, \bx) \nu^{(\btpsi)}(\dd t, \dd\bx) , \label{GIR1}\\
  \int_0^t \int_{\RR^d}
   \Bigl(\sqrt{V^{(\btpsi,\bpsi)}(s, \bx)} - 1\Bigr)^2 \,
   \nu^{(\btpsi)}(\dd s, \dd\bx)
  < \infty , \label{GIR2} \\
  B^{(\bpsi)}_t
  = B^{(\btpsi)}_t
    + \int_0^t
       c_s^{(\bpsi)} \beta^{(\btpsi,\bpsi)}_s \, \dd F_s^{(\bpsi)}
    + \int_0^t \int_{\RR^d}
       (V^{(\btpsi,\bpsi)}(s, \bx) - 1) h(\bx) \,
       \nu^{(\btpsi)}(\dd s, \dd\bx) , \label{GIR3} \\
  \int_0^t
   (\beta^{(\btpsi,\bpsi)}_s)^\top c_s^{(\bpsi)} \beta^{(\btpsi,\bpsi)}_s \, \dd F_s^{(\bpsi)}
  < \infty , \label{GIR4}
 \end{gather}
 $\PP_\bpsi$-almost sure for every \ $t \in [0, \infty)$.
\ Further, assume that for each \ $\bpsi \in \Psi$, \ local uniqueness holds for the
 martingale problem on the canonical space corresponding to the triplet
 \ $(B^{(\bpsi)}, C, \nu^{(\bpsi)})$ \ with the given initial value \ $\bx_0$
 \ with \ $\PP_\bpsi$ \ as its unique solution.
Then for each \ $T \in [0, \infty)$, \ $\PP_{\bpsi,T}$ \ is absolutely continuous
 with respect to \ $\PP_{\btpsi,T}$,
 \ where \ $\PP_{\bpsi,T} := \PP_{\bpsi\!\!}|_{\cD_T([0,\infty), \RR^d)}$ \ denotes the
 restriction of \ $\PP_\bpsi$ \ to \ $\cD_T([0,\infty), \RR^d)$ \ (similarly for
 \ $\PP_{\btpsi,T}$), \ and, under \ $\PP_{\btpsi,T}$, \ the corresponding
 likelihood-ratio process takes the form
 \begin{align}\label{RN_general}
  \begin{split}
  \log \frac{\dd \PP_{\bpsi,T}}{\dd \PP_{\btpsi,T}}(\eta)
  &= \int_0^T (\beta^{(\btpsi,\bpsi)}_s)^\top \, \dd(\eta^\cont)^{(\btpsi)}_s
     - \frac{1}{2}
       \int_0^T
       (\beta^{(\btpsi,\bpsi)}_s)^\top c_s^{(\bpsi)} \beta^{(\btpsi,\bpsi)}_s \, \dd F_s^{(\bpsi)} \\
  &\quad
    + \int_0^T \int_{\RR^d}
       (V^{(\btpsi,\bpsi)}(s, \bx) - 1) \,
       (\mu^\eta - \nu^{(\btpsi)})(\dd s, \dd\bx)\\
  &\quad
    + \int_0^T \int_{\RR^d}
       (\log(V^{(\btpsi,\bpsi)}(s, \bx))
        - V^{(\btpsi,\bpsi)}(s, \bx) + 1) \,
       \mu^\eta(\dd s, \dd\bx)
  \end{split}
 \end{align}
 for all \ $T \in (0, \infty)$, \ see Jacod and Shiryaev \cite[Theorem III.5.34]{JSh}.

In what follows we give a proof for \eqref{RN_general} using Jacod and Shiryaev
 \cite{JSh}, since in the literature we could not find a detailed proof.
Using the notations of Jacod and Shiryaev \cite{JSh}, under \ $\PP_\bpsi$ \ the
 triplets \ $(B^{(\bpsi)}, C, \nu^{(\bpsi)})$ \ and
 \ $(B^{(\btpsi)}, C, \nu^{(\btpsi)})$ \ satisfy III.5.5 in Jacod and Shiryaev
 \cite{JSh} with \ $Y = V^{(\btpsi,\bpsi)}$, \ $A = F^{(\bpsi)}$
 \ and \ $c = c^{(\bpsi)}$, \ and the
 filtration \ $\cD_t([0,\infty), \RR^d)$, \ $t \in [0, \infty)$, \ is generated by
 \ $(\eta_t)_{t\in[0,\infty)}$.
\ Moreover, \ $a_t = \nu^{(\bpsi)}(\{t\} \times \RR^d) = 0$ \ and
 \ $\hY_t
    = \int_{\RR^d} V^{(\btpsi,\bpsi)}(t, \bx) \, \nu^{(\bpsi)}(\{t\} \times \dd\bx)
    = 0$
 \ $\PP_\bpsi$-almost surely for all \ $t \in [0, \infty)$, \ hence
 \ $\sigma
    = \inf\{t \in [0, \infty)
            : \text{either $\hY_t > 1$, or $a_t = 1$ and $\hY_t < 1$}\} = \infty$
 \ $\PP_\bpsi$-almost surely.
Then, by \eqref{GIR2} and \eqref{GIR4}, \ we have
 \[
   H_t
   = \int_0^t
      (\beta^{(\btpsi,\bpsi)}_s)^\top c_s^{(\bpsi)} \beta^{(\btpsi,\bpsi)}_s \,
      \dd F_s^{(\bpsi)}
     + \int_0^t \int_{\RR^d}
        \Bigl(\sqrt{V^{(\btpsi,\bpsi)}(s, \bx)} - 1\Bigr)^2 \,
        \nu^{(\btpsi)}(\dd s, \dd\bx)
   < \infty
 \]
 $\PP_\bpsi$-almost sure for every \ $t \in [0, \infty)$.
\ Consequently, \ $T_n = \inf\{t \in [0, \infty) : H_t \geq n\} \to \infty$ \ as
 \ $n \to \infty$ \ $\PP_\bpsi$-almost sure (due to the fact that
 \ $(H_t)_{t\in[0,\infty)}$ \ is a nondecreasing process), and hence
 \ $\PP_\bpsi(\Delta = [0, \infty)) = 1$, \ and the Hypothesis III.5.29 in Jacod
 and Shiryaev \cite{JSh} holds.
Thus, by Theorem III.5.34 in Jacod and Shiryaev \cite{JSh}, \ $\PP_{\bpsi,T}$ \ is
 absolutely continuous with respect to \ $\PP_{\btpsi,T}$ \ for all
 \ $T\in[0,\infty)$, \ and under \ $\PP_{\btpsi,T}$, \ the density process
 (likelihood ratio process) \ $(Z_T)_{T\in[0,\infty)}$ \ takes the form
 \[
   Z_T
   = \frac{\dd \PP_{\bpsi,T}}{\dd \PP_{\btpsi,T}}(\eta)
   = \exp\biggl\{N_T
                 - \frac{1}{2}
                   \int_0^T
                    (\beta^{(\btpsi,\bpsi)}_s)^\top c_s^{(\bpsi)}
                    \beta^{(\btpsi,\bpsi)}_s \,
                    \dd F_s^{(\bpsi)}\biggr\}
     \prod_{s\in[0,T]} (1 + \Delta N_s) \ee^{-\Delta N_s}
 \]
 for \ $T \in [0, \infty)$ \ with
 \[
   N_T := \int_0^T (\beta^{(\btpsi,\bpsi)}_s)^\top \, \dd(\eta^\cont)^{(\btpsi)}_s
          + \int_0^T \int_{\RR^d}
             (V^{(\btpsi,\bpsi)}(s, \bx) - 1) \,
             (\mu^\eta - \nu^{(\btpsi)})(\dd s, \dd\bx)
 \]
 for \ $T \in [0, \infty)$.
\ Further, the density process \ $(Z_T)_{T\in[0,\infty)}$ \ satisfies
 \[
   Z_T = 1 + \int_0^T Z_{s-} \,\dd N_s \qquad
   \text{under \ $\PP_{\btpsi,T}$ \ for each \ $T \in [0, \infty)$,}
 \]
 see Jacod and Shiryaev \cite[III.5.20]{JSh}.
Taking into account the fact that
 \ $\PP_\bpsi(\nu^{(\bpsi)}(\{s\} \times \RR^d) = 0)=1$ \ for every
 \ $\bpsi \in \Psi$, \ $s \in [0,\infty)$, \ and the definition of the
 stochastic integral with respect to the random measure
 \ $\mu^\eta - \nu^{(\btpsi)}$ \ (see Jacod and Shiryaev
 \cite[Definition II.1.27]{JSh}), we obtain
 \[
   \Delta N_s(\omega)
   = (V^{(\btpsi,\bpsi)}(s, \Delta\eta_s(\omega))(\omega) - 1)
     \bbone_{\{\Delta\eta_s(\omega)\ne\bzero\}} ,
   \qquad \omega \in D([0,\infty), \RR^d) , \quad s \in [0, \infty) .
 \]
Hence, using that \ $V^{(\btpsi,\bpsi)}$ \ is positive, we have
 \ $\PP_{\btpsi}(\inf\{t\in[0,\infty) : \Delta N_t = -1 \} = \infty )=1$,
 \ and consequently, \ $\PP_{\btpsi}( Z_T \in(0,\infty))=1$ \ for all \ $T\in(0,\infty)$,
 \ see Jacod and Shiryaev \cite[Theorem I.4.61]{JSh}.
Further,
 \begin{align*}
  \prod_{s\in[0,T]} (1 + \Delta N_s) \ee^{-\Delta N_s}
  &= \prod_{s\in[0,T]}
      (V^{(\btpsi,\bpsi)}(s, \Delta\eta_s) \bbone_{\{\Delta\eta_s\ne\bzero\}}
             +  \bbone_{\{\Delta\eta_s = \bzero\}})
     \ee^{-(V^{(\btpsi,\bpsi)}(s,\Delta\eta_s)-1)\bbone_{\{\Delta\eta_s\ne\bzero\}}} \\
  &\quad
   = \exp\biggl\{\sum_{s\in[0,T]}
                  \big(\log(V^{(\btpsi,\bpsi)}(s, \Delta\eta_s))
                       - V^{(\btpsi,\bpsi)}(s, \Delta\eta_s) + 1\big)
                  \bbone_{\{\Delta\eta_s\ne\bzero\}}\biggr\} \\
  &\quad
   = \exp\biggl\{\int_0^T \int_{\RR^d}
                  \big(\log(V^{(\btpsi,\bpsi)}(s, \bx)) - V^{(\btpsi,\bpsi)}(s, \bx) + 1\big) \,
                  \mu^\eta(\dd s, \dd\bx)\biggr\} ,
 \end{align*}
 yielding \eqref{RN_general}, where the existence of the integral in the exponent above
 follows from the facts that \ $\PP_{\btpsi}( Z_T \in(0,\infty))=1$, \ assumption
 \eqref{GIR4} and \ $\PP_{\btpsi}( N_T \in\RR)=1$ \ for all \ $T\in\RR$ \ (due to Jacod
 and Shiryaev \cite[Proposition III.5.10 and III.5.12]{JSh}).

\section{Limit theorems for continuous local martingales}\label{App_MCLT}

In what follows we recall some limit theorems for continuous local martingales.
We use these limit theorems for studying the asymptotic behaviour of the MLE of
 \ $\bpsi = (\theta, \kappa, \mu)$.
\ First we recall a strong law of large numbers for continuous local martingales.

\begin{Thm}{\bf (Liptser and Shiryaev \cite[Lemma 17.4]{LipShiII})}
\label{DDS_stoch_int}
Let \ $\bigl( \Omega, \cF, (\cF_t)_{t\in[0,\infty)}, \PP \bigr)$ \ be a filtered
 probability space satisfying the usual conditions.
Let \ $(M_t)_{t\in[0,\infty)}$ \ be a square-integrable continuous local
 martingale with respect to the filtration \ $(\cF_t)_{t\in[0,\infty)}$ \ such
 that \ $\PP(M_0 = 0) = 1$.
\ Let \ $(\xi_t)_{t\in[0,\infty)}$ \ be a progressively measurable process such
 that
 \[
   \PP\left( \int_0^t \xi_u^2 \, \dd \langle M \rangle_u < \infty \right) = 1 ,
   \qquad t \in [0,\infty) ,
 \]
 and
 \begin{align}\label{SEGED_STRONG_CONSISTENCY2}
  \int_0^t \xi_u^2 \, \dd \langle M \rangle_u \as \infty \qquad
  \text{as \ $t \to \infty$,}
 \end{align}
 where \ $(\langle M \rangle_t)_{t\in[0,\infty)}$ \ denotes the (predictable) quadratic
 variation process of \ $M$.
\ Then
 \begin{align}\label{SEGED_STOCH_INT_SLLN}
  \frac{\int_0^t \xi_u \, \dd M_u}
       {\int_0^t \xi_u^2 \, \dd \langle M \rangle_u} \as 0 \qquad
  \text{as \ $t \to \infty$.}
 \end{align}
If \ $(M_t)_{t\in[0,\infty)}$ \ is a standard Wiener process, the progressive
 measurability of \ $(\xi_t)_{t\in[0,\infty)}$ \ can be relaxed to
 measurability and adaptedness to the filtration \ $(\cF_t)_{t\in[0,\infty)}$.
\end{Thm}

The next theorem is a special case of the central limit theorem for multidimensional
 square-integrable continuous local martingales, see, e.g., Jacod and Shiryaev \cite[Corollary VIII.3.24]{JSh} or van Zanten \cite[Theorem 4.1]{Zan}.

\begin{Thm}\label{MCLT}
Let \ $\bigl( \Omega, \cF, (\cF_t)_{t\in[0,\infty)}, \PP \bigr)$ \ be a filtered
 probability space satisfying the usual conditions.
Let \ $(\bM_t)_{t\in[0,\infty)}$ \ be a $d$-dimensional square-integrable
 continuous local martingale with respect to the filtration
 \ $(\cF_t)_{t\in[0,\infty)}$ \ such that \ $\PP(\bM_0 = \bzero) = 1$ \ and
 \[
   t^{-1} \langle \bM \rangle_t \stoch \bU
   \qquad \text{as \ $t \to \infty$,}
 \]
 where \ $\bU \in \RR^{d\times d}$.
\ Then
 \[
   t^{-1/2} \bM_t \distr \cN_d(\bzero, \bU) \qquad \text{as \ $t \to \infty$.}
 \]
\end{Thm}

\section{A version of the continuous mapping theorem}
\label{CMT}

The following version of continuous mapping theorem can be found for example
 in Kallenberg \cite[Theorem 3.27]{K}.

\begin{Lem}\label{Lem_Kallenberg}
Let \ $(\cS_1, d_{\cS_1})$ \ and \ $(\cS_2, d_{\cS_2})$ \ be metric spaces and
 \ $(\xi_n)_{n \in \NN}$, \ $\xi$ \ be random elements with values in \ $\cS_1$
 \ such that \ $\xi_n \distr \xi$ \ as \ $n \to \infty$.
\ Let \ $F : \cS_1 \to \cS_2$ \ and \ $F_n : \cS_1 \to \cS_2$, \ $n \in \NN$,
 \ be measurable mappings and \ $\cC \in \cB(\cS_1)$ \ such that
 \ $\PP(\xi \in \cC) = 1$ \ and
 \ $\lim_{n \to \infty} d_{\cS_2}(F_n(s_n), F(s)) = 0$ \ if
 \ $\lim_{n \to \infty} d_{\cS_1}(s_n,s) = 0$ \ and \ $s \in \cC$.
\ Then \ $F_n(\xi_n) \distr F(\xi)$ \ as \ $n \to \infty$.
\end{Lem}

\section{Explicit formula for a density function}
\label{App_density}

We show that the mixed normal but non-normal density function of
 \ $\frac{\varrho\sigma}{\kappa\tau}
    + \frac{\sigma\sqrt{1-\varrho^2}}{\kappa\sqrt{\tau}} Z_2$,
 \ which is the limit distribution of \ $T (\hmu_T - \mu)$ \ as \ $T \to \infty$ \ in Theorem
 \ref{Thm_MLE=}, has the form
 \begin{equation}\label{exp_density}
  f(x)
  = \frac{\kappa}{2\pi\sigma\sqrt{1-\varrho^2}}
    \int_0^\infty
     \frac{1}{t}
     \exp\biggl\{- \frac{1}{2t}
                 - \frac{(\kappa xt-\varrho\sigma)^2}{2\sigma^2(1-\varrho^2)t}\biggr\}
     \, \dd t , \qquad x \in \RR \setminus \{0\} ,
 \end{equation}
 and
 \begin{equation}\label{density_lim}
  \lim_{x\to0} f(x) = \infty .
 \end{equation}
It is known that the density function of \ $\tau$ \ takes the form
 \ $f_\tau(t) = (2\pi t^3)^{-1/2} \ee^{-1/(2t)} \bbone_{\RR_{++}}(t)$, \ $t \in \RR$.
\ Using the independence of \ $\tau$ \ and \ $Z_2$, \ we have
 \begin{align*}
  \PP\biggl(\frac{\varrho\sigma}{\kappa\tau}
            + \frac{\sigma\sqrt{1-\varrho^2}}{\kappa\sqrt{\tau}} Z_2 \leq x
            \,\bigg|\, \tau = t\biggr)
  &= \PP\biggl(Z_2 \leq \frac{(x-\frac{\varrho\sigma}{\kappa\tau})\kappa\sqrt{\tau}}
                             {\sigma\sqrt{1-\varrho^2}}
               \,\bigg|\, \tau = t\biggr) \\
  &= \frac{1}{\sqrt{2\pi}}
     \int_{-\infty}^{\frac{\kappa xt-\varrho\sigma}{\sigma\sqrt{(1-\varrho^2)t}}}
      \ee^{-u^2/2} \, \dd u , \qquad x \in \RR, \quad t \in \RR_{++} .
 \end{align*}
By the law of total expectation, we obtain
 \[
   \PP\biggl(\frac{\varrho\sigma}{\kappa\tau}
             + \frac{\sigma\sqrt{1-\varrho^2}}{\kappa\sqrt{\tau}} Z_2 \leq x\biggr)
   = \int_0^\infty g(x, t) \, \dd t , \qquad x \in \RR ,
 \]
 with
 \[
   g(x, t)
   := \frac{1}{\sqrt{2\pi t^3}} \ee^{-\frac{1}{2t}}
      \biggl(\frac{1}{\sqrt{2\pi}}
             \int_{-\infty}^{\frac{\kappa xt-\varrho\sigma}{\sigma\sqrt{(1-\varrho^2)t}}}
              \ee^{-u^2/2} \, \dd u\biggr) , \qquad x \in \RR, \quad t \in \RR_{++} .
 \]
The aim of the following discussion is to show that, by the dominated convergence theorem,
 \begin{align}\label{help_density_1}
   \frac{\dd}{\dd x} \int_0^\infty g(x, t) \, \dd t
   = \lim_{h\to0} \int_0^\infty \frac{g(x+h,t)-g(x,t)}{h} \, \dd t
   = f(x) , \qquad
   x \in \RR \setminus \{0\} .
 \end{align}
For all \ $x \in \RR$ \ and \ $t \in \RR_{++}$, \ we have
 \begin{align}\label{help_density_2}
   \lim_{h\to0} \frac{g(x+h,t)-g(x,t)}{h}
   = \partial_1 g(x, t)
   = \frac{\kappa}{2\pi\sigma\sqrt{1-\varrho^2}t}
     \exp\biggl\{- \frac{1}{2t}
                 - \frac{(\kappa xt-\varrho\sigma)^2}{2\sigma^2(1-\varrho^2)t}\biggr\} .
 \end{align}
Moreover, for all \ $x \in \RR \setminus \{0\}$, \ $t \in \RR_{++}$ \ and
 \ $h \in \bigl[-\frac{|x|}{2}, \frac{|x|}{2}\bigr] \setminus \{0\}$, \ we have
 \[
   \biggl|\frac{g(x + h, t) - g(x, t)}{h}\biggr|
   \leq \sup_{\xi\in[-1,1]} |\partial_1 g(x + \xi h, t)|
   \leq \sup_{\xi\in[-1,1]}
        \Bigl|\partial_1 g\Bigl(x + \xi \frac{x}{2}, t\Bigr)\Bigr| .
 \]
For \ $x \in \RR \setminus \{0\}$, \ $\xi \in [-1, 1]$ \ and \ $t \in \RR_{++}$, \ we have
 \[
   \Bigl|\partial_1 g\Bigl(x + \xi \frac{x}{2}, t\Bigr)\Bigr|
   \leq \frac{\kappa}{2\pi\sigma\sqrt{1-\varrho^2}t} \ee^{-\frac{1}{2t}}
   \leq \frac{\kappa}{\ee\pi\sigma\sqrt{1-\varrho^2}} ,
 \]
 since \ $\sup_{t\in\RR_{++}} t^{-1} \ee^{-\frac{1}{2t}} = 2 \ee^{-1}$.
\ Further, for \ $x \in \RR \setminus \{0\}$, \ $\xi \in [-1, 1]$ \ and \ $t \in \RR_{++}$,
 \ we have
 \[
   \Bigl|\kappa \Bigl(x + \xi \frac{x}{2}\Bigr) t - \varrho \sigma\Bigr|
   \geq \Bigl|\kappa\Bigl(x + \xi \frac{x}{2}\Bigr)t\Bigr| - |\varrho\sigma|
   \geq \frac{1}{2} \Bigl|\kappa \Bigl(x + \xi \frac{x}{2}\Bigr) t\Bigr|
 \]
 whenever
 \[
   \frac{1}{2} \Bigl|\kappa\Bigl(x + \xi \frac{x}{2}\Bigr)t\Bigr| \geq |\varrho\sigma| ,
 \]
 equivalently, whenever
 \[
   t \geq \frac{2|\varrho\sigma|}{\bigl|\kappa\bigl(x+\xi\frac{x}{2}\bigr)\bigr|}
       = \frac{2|\varrho|\sigma}{\kappa\bigl(1+\frac{1}{2}\xi\bigr)|x|} ,
 \]
 which holds if
 \[
   t \geq  \frac{4|\varrho|\sigma}{\kappa|x|} + 1 =: T_0(x) \in \RR_{++} .
 \]
Consequently, for \ $x \in \RR \setminus \{0\}$, \ $\xi \in [-1, 1]$ \ and
 \ $t \in (T_0(x), \infty)$, \ we have
 \[
   \Bigl|\kappa\Bigl(x + \xi \frac{x}{2}\Bigr)t - \varrho\sigma\Bigr|
   \geq \frac{1}{2} \Bigl|\kappa \Bigl(x + \xi \frac{x}{2} \Bigr) t\Bigr|
   = \frac{1}{2} \kappa \Bigl(1 + \frac{\xi}{2}\Bigr) |x| t
   \geq \frac{\kappa|x|t}{4} ,
 \]
 and hence, by the second equality in \eqref{help_density_2},
 \[
   \Bigl|\partial_1 g\Bigl(x + \xi \frac{x}{2}, t\Bigr)\Bigr|
   \leq \frac{\kappa}{2\pi\sigma\sqrt{1-\varrho^2}T_0(x)}
        \exp\biggl\{-\frac{\kappa^2x^2}{32\sigma^2(1-\varrho^2)} t \biggr\} .
 \]
We conclude that for \ $x \in \RR \setminus \{0\}$, \ $\xi \in [-1, 1]$ \ and
 \ $t \in \RR_{++}$, \ we have
 \[
   \Bigl|\partial_1 g\Bigl(x + \xi \frac{x}{2}, t\Bigr)\Bigr| \leq G(t, x)
 \]
 with
 \[
   G(t, x)
   := \frac{\kappa}{\ee\pi\sigma\sqrt{1-\varrho^2}} \bbone_{[0,T_0(x)]}(t)
      + \frac{\kappa}{2\pi\sigma\sqrt{1-\varrho^2}T_0(x)}
        \exp\biggl\{-\frac{\kappa^2x^2}{32\sigma^2(1-\varrho^2)} t \biggr\}
        \bbone_{(T_0(x),\infty)}(t)
 \]
 for \ $t \in \RR_{++}$ \ and \ $x \in \RR \setminus \{0\}$, \ and the function
 \ $\RR_+ \ni t \mapsto G(t, x)$ \ is integrable on \ $\RR_+$, \ hence the dominated
 convergence theorem can be used, and we obtain \eqref{exp_density}.

One can derive \eqref{help_density_1} in another way.
Since \ $g(x,t)\leq \frac{1}{\sqrt{2\pi t^3}} \ee^{-\frac{1}{2t}}$, \ $x\in\RR$, \ $t\in\RR_{++}$, \
 the improper integral \ $\int_0^\infty g(x,t) \, \dd t$ \ is uniformly convergent
 for \ $x\in\RR$.
\ Further, for any \ $t, a, A \in \RR_{++}$ \ and \ $x \in \RR$ \ with \ $a < |x| < A$, \ we have
 \[
   G(t,x) \leq \frac{\kappa}{\ee\pi\sigma\sqrt{1-\varrho^2}} \bbone_{[0,T_0(a)]}(t)
          + \frac{\kappa}{2\pi\sigma\sqrt{1-\varrho^2}T_0(A)}
            \exp\biggl\{-\frac{\kappa^2a^2}{32\sigma^2(1-\varrho^2)} t \biggr\}
            \bbone_{(T_0(A),\infty)}(t)
 \]
 showing that the improper integral \ $\int_0^\infty \partial_1 g(x,t) \, \dd t$ \ is uniformly
 convergent for \ $a< |x| <A$ \ with \ $a < A$, \ $a, A \in \RR_{++}$.
\ This together with the continuity of the functions
\ $\RR \times \RR_{++} \ni (x, t) \mapsto g(x,t) \in \RR$ \  and
 \ $\RR \times \RR_{++} \ni (x, t) \mapsto \partial_1 g(x,t) \in \RR$ \
 yield \eqref{help_density_1}, see, e.g., Lang \cite[pages 337-339]{Lan}.

Moreover, for all \ $x \in \RR \setminus \{0\}$, \ we have
 \[
   f(x) \geq \frac{\kappa}{2\pi\sigma\sqrt{1-\varrho^2}}
             \exp\biggl\{- \frac{1}{2}
                         + \frac{2\varrho\sigma\kappa x}{2\sigma^2(1-\varrho^2)}
                         - \frac{\varrho^2\sigma^2}{2\sigma^2(1-\varrho^2)}\biggr\}
             \int_1^\infty
              \frac{1}{t} \exp\biggl\{- \frac{\kappa^2x^2t}{2\sigma^2(1-\varrho^2)}\biggr\}
              \, \dd t ,
 \]
 where
 \[
   \lim_{x\to0} \exp\biggl\{- \frac{1}{2}
                            + \frac{2\varrho\sigma\kappa x}{2\sigma^2(1-\varrho^2)}
                            - \frac{\varrho^2\sigma^2}{2\sigma^2(1-\varrho^2)}\biggr\}
   = \exp\biggl\{- \frac{1}{2} - \frac{\varrho^2\sigma^2}{2\sigma^2(1-\varrho^2)}\biggr\}
 \]
 and
 \[
   \int_1^\infty
    \frac{1}{t} \exp\biggl\{- \frac{\kappa^2x^2t}{2\sigma^2(1-\varrho^2)}\biggr\} \, \dd t
   = \int_{x^2}^\infty
      \frac{1}{u} \exp\biggl\{- \frac{\kappa^2u}{2\sigma^2(1-\varrho^2)}\biggr\} \, \dd u
   \to \infty \qquad \text{as \ $x \to 0$,}
 \]
 hence we obtain \eqref{density_lim}.

\section*{Acknowledgements}
We would like to thank the referees for their comments that helped us to improve the
 paper.

\end{document}